\documentclass[11pt]{article}
\usepackage{amssymb,amsmath}
\usepackage{graphicx}
\usepackage{color}
\newcommand{\refe}[1]{(\ref{#1})}
\newcommand{\dst}{\displaystyle}
\newcommand{\RR}{{\mathbb R} }

\newcommand{\hlf}{\textstyle{\frac 12}}
\def\ams#1{\rm AMS classification scheme numbers: #1\par }
\def\keyw#1{\rm Keywords: #1\par }
\newcommand{\address}[1]
{\begin{center}   \rm\raggedright #1   \end{center}}
\def\ead#1{\vspace*{5pt}
\address{E-mail: \mailto{#1}}}
\def\mailto#1{{\tt #1}}
\newtheorem{theorem}{Theorem}[section]
\newtheorem{rem}[theorem]{Remark} 
 
\newtheorem{lemma}[theorem]{Lemma} 

\newenvironment{Proof}[1][Proof]{\begin{trivlist}
\item[\hskip \labelsep {\bfseries #1}]}{\end{trivlist}}

\newcommand{\qed}{\nobreak \ifvmode \relax \else
      \ifdim\lastskip<1.5em \hfill $\blacksquare$ \fi}

\numberwithin{equation}{section}

\begin{document}


\title{The orthogonality of q-classical polynomials of the Hahn class: A geometrical approach}
\author{R. \'{A}lvarez-Nodarse\dag, \ R. Sevinik-Ad\i g{\"{u}}zel\ddag\ and \ H. Ta\c{s}eli\ddag}
\date{\today}
\maketitle\address{\dag IMUS \& Departamento de An\'alisis Matem\'atico,
Universidad de Sevilla. Apdo. 1160,
E-41080  Sevilla, Spain}
\address{\ddag Department of Mathematics, Middle East Technical University 
(METU), 06531, Ankara, Turkey}

\ams{33D45, 42C05}
\medskip

\keyw{$q$-polynomials, orthogonal polynomials on $q$-linear lattices, $q$-Hahn class}

\begin{abstract}
The idea of this review article is to discuss in a unified way
the orthogonality of all  positive definite polynomial solutions of the 
$q$-hypergeometric difference equation on the $q$-linear lattice
by means of a qualitative analysis of the $q$-Pearson equation.
Therefore, our method differs from the standard ones which are based 
on the Favard theorem, the three-term recurrence relation 
and the difference equation of hypergeometric type. 
Our approach enables us to extend the orthogonality 
relations for some well-known $q$-polynomials of the Hahn class 
to a larger set of their parameters. A short version of this paper
appeared in SIGMA \textbf{8} (2012), 042, 30 pages 
\texttt{http://dx.doi.org/10.3842/SIGMA.2012.042}.
\end{abstract}

\ead{ran@us.es, sevinikrezan@gmail.com, taseli@metu.edu.tr}


\section{Introduction}
The so-called $q$-polynomials are of great interest inside the class 
of special functions since they play an important role in the 
treatment of several problems such as Eulerian series and continued fractions  
\cite{An, Fine},  $q$-algebras and quantum groups \cite {Korn1, Korn2, Vil} and $q$-oscillators
\cite{ata08,ran-ata-cos,grun}, and references therein,  among others. 

A $q$-analog of the Chebychev's discrete orthogonal polynomials is due to Markov
in 1884  \cite[page 43]{AA}, which can be regarded as the first example 
of a $q$-polynomial family. In 1949, Hahn introduced the $q$-Hahn class \cite{H} 
including the big $q$-Jacobi polynomials, on the exponential lattice although
he did not use this terminology. In fact, he did not give 
the orthogonality relations of the big $q$-Jacobi
polynomials in \cite{H} which was done by Andrew and Askey \cite{AA}.
During the last decades the $q$-polynomials have
been studied by many authors from different points of view. 
There are two most recognized approaches. The first approach, initiated by 
the work of Askey and Wilson \cite{AW} 
(see also Andrews and Askey \cite{AA}) is based on the basic hypergeometric 
series \cite{An, GR}. The second approach is due to Nikiforov and Uvarov \cite{NUpaper1, NU}
and uses the analysis of difference equations on non-uniform lattices.  
The readers are also referred to the surveys \cite{ARS,NSU,NUpaper2,sus89}. 
These approaches are associated with the so-called $q$-Askey scheme \cite{KLS} and the Nikiforov-Uvarov
scheme \cite{NUpaper2}, respectively. 
Another approach was published in \cite{MAJF} 
where the authors proved several characterizations of the $q$-polynomials 
starting from the so-called distributional $q$-Pearson equation
(for the non $q$-case see e.g. \cite{gar95,marpet} and references therein).

In particular, in  \cite{MAJF} a classification of all possible
families of orthogonal polynomials on the exponential lattice was established, and latter on 
in \cite{NM} the comparison with the $q$-Askey and Nikiforov-Uvarov schemes was done, resulting in two
new families of orthogonal polynomials. Furthermore, an important contribution to the theory of 
(orthogonal) $q$-polynomials, and in particular, to the theory of orthogonal $q$-polynomials on the 
linear exponential lattice, appeared in the recent book \cite{KLS}.
The corresponding table is generally called the $q$-Hahn tableau (see e.g., Koornwinder \cite{Korn2}).
The $q$-polynomials belonging to this class are the solutions of the 
$q$-difference equation of hypergeometric type ($q$-EHT) \cite{H}\begin{equation}\label{qEHT1}
\sigma_1(x; q)D_{q^{-1}}D_qy(x, q)+\tau(x, q)D_qy(x, q)+\lambda(q) y(x, q)=0.
\end{equation}
One way of deriving the $q$-EHT \refe{qEHT1} 
whose bounded solutions are the $q$-polynomials of the Hahn class, 
is to discretize the classical differential equation of hypergeometric type (EHT)
\begin{equation}\label{EHT}
\sigma(x)y''+\tau(x)y'+\lambda y=0,
\end{equation}
where $\sigma(x)$ and $\tau(x)$ are polynomials of at
most second and first degree, respectively, and $\lambda$ 
is a constant \cite{ran, ARS, marpet, NSU, NU}. To this end, we can use 
the approximations (see e.g. \cite[\S13, page 142]{NU})
\begin{equation*}
y'(x)\sim\frac{1}{1+q}[D_qy(x)+qD_{q^{-1}}y(x)]\quad {\rm and} \quad y''(x)\sim\frac{2q}{1+q}D_qD_{q^{-1}}y(x)
\quad {\rm as} \quad q\rightarrow1
\end{equation*}
for the derivatives in \refe{EHT}, where we use the standard notation
for the $q$ and $q^{-1}$-Jackson derivatives of $y(x)$ \cite{GR, kac}, i.e., 
\begin{equation*}
D_{\zeta}y(x)=\frac{y(x)-y(\zeta x)}{(1-\zeta)x}, \qquad \zeta\in\mathbb{C}\setminus\{0,\pm1\}
\end{equation*}%
for $x\neq0$ and $D_\zeta y(0)=y'(0)$, provided that $y'(0)$ exists.
This leads to
the $q$-EHT \refe{qEHT1} 
where
\begin{equation*}
\sigma_1(x;q):=\frac{2}{1+q}\left[\sigma(x)-\frac{1}{2}(q-1)x\tau(x)\right],\,\,\,\tau(x,q):=\tau(x),
\,\,\, \lambda(q):=\lambda, \,\,\, y(x, q):=y(x).
\end{equation*}
Notice here the relations 
$D_q=D_{q^{-1}}+(q-1)xD_qD_{q^{-1}}$ and $D_qD_{q^{-1}}=q^{-1}D_{q^{-1}}D_q$
so that \refe{qEHT1} can be rewritten in the equivalent form
\begin{equation}\label{qEHT2}
\sigma_2(x; q)D_qD_{q^{-1}}y(x, q)+\tau(x,q)D_{q^{-1}}y(x, q)
+\lambda(q)y(x, q)=0,
\end{equation}
where   
\begin{equation}\label{sigmaq12}
\sigma_2(x, q):=q\left[\sigma_1(x, q)+(1-q^{-1})x\tau(x, q)\right].
\end{equation}
It should be noted that the $q$-EHT \refe{qEHT1} and \refe{qEHT2} correspond to the second 
order linear difference equations of hypergeometric type on the linear exponential 
lattices $x(s)=c_1q^s+c_2$ and $x(s)=c_1q^{-s}+c_2$, respectively \cite{ran, ARS, NSU}. 

Notice also that \refe{qEHT1} (or \refe{qEHT2}) can be written in a very 
convenient form  \cite{NM,KLS,ks}
\begin{equation*} 
\sigma_2(x, q)D_qy(x, q)-q\sigma_1(x, q)D_{q^{-1}}y(x, q)+(q-1)x\lambda(q)y(x, q)=0,
\end{equation*}
where the coefficients $\sigma_1(x;q)$ and $\sigma_2(x;q)$ are polynomials 
of at most 2nd degree and $\tau(x,q)$ is a 1st degree polynomial in $x$.

Notice that the $q$-EHT \refe{qEHT1} can be written in the self-adjoint form
$$
D_q\left[\sigma_1(x, q)\rho(x, q)D_{q^{-1}}y(x)\right]
+q^{-1}\lambda(q)\rho(x, q)y(x)=0, 
$$
where $\rho$ is a function satisfying 
the so-called $q$-Pearson equation
$D_q\left[\sigma_1(x, q)\rho\right]=q^{-1}\tau(x, q)\rho$
that can be written as 
\begin{equation}\label{openqpearson}
\frac{\rho(qx, q)}{\rho(x, q)}=\frac{\sigma_1(x, q)
+(1-q^{-1})x\tau(x, q)}{\sigma_1(qx, q)}=\frac{q^{-1}
\sigma_2(x, q)}{\sigma_1(qx, q)},
\end{equation}
or, equivalently,
\begin{equation}\label{ide-sr}
 \sigma_2(x, q)\rho(x, q)=q \sigma_1(qx, q)\rho(qx, q).  
\end{equation}

In this paper we study, without loss of generality, the 
$q$-EHT \refe{qEHT1}, assume $0<q<1$ and take $\lambda(q)$ as
\begin{equation*}
\lambda(q):=\lambda_n(q)=-[n]_q\left[\tau'(0, q)
+\hlf[n-1]_{q^{-1}}\sigma_1''(0, q)\right], \quad n\in\mathbb{N}_0=\mathbb{N}\cup\{0\}
\end{equation*}
since we are interested only in the polynomial solutions \cite{ran, ARS, NSU}. 
For more details on the $q$-polynomials of the $q$-Hahn tableau  
we refer the readers to the works \cite{ran, ran1,  RA, NM, ARS,  DN, KLS, Korn2, MM, NSU, NUpaper1, NU, NUpaper2, sus89},
and references therein.

In this paper, we deal with the orthogonality properties of the $q$-polynomials of the 
$q$-Hahn tableau from a different viewpoint than the one used in \cite{KLS}. 
In \cite{KLS}, the authors presented a unified study of the orthogonality of $q$-polynomials based on 
the Favard Theorem. Here, the main idea is to provide a relatively simple geometrical analysis of the $q$-Pearson 
equation by taking into account every possible rational form of the polynomial coefficients of the 
$q$-difference equation. 
Roughly, our qualitative analysis is concerned with the examination of the behavior of the graphs 
of the ratio $\rho(qx, q)/\rho(x, q)$ by means of the definite right hand side (r.h.s.) of \refe{openqpearson}
in order to find out a suitable  $q$-weight function. Such a qualitative analysis implies 
all possible orthogonality relations among the polynomial solutions of the 
$q$-difference equation in question.  Moreover, it allows us to extend the orthogonality 
relations for some well-known $q$-polynomials of the Hahn class 
to a larger set of their parameters (see sections \ref{s4.1} and \ref{s4.2}).
A first attempt of using a geometrical approach 
for studying the orthogonality of $q$-polynomials of the $q$-Hahn class was presented 
in \cite{DN}. However, the study is far from being complete and only some partial results 
were obtained. We will fill this gap in this review paper.

Our main goal is to study each orthogonal polynomial system or sequence (OPS), which 
is orthogonal with respect to (w.r.t.) a $q$-weight function $\rho(x, q)>0$ 
satisfying the $q$-Pearson equation as well as certain boundary conditions (BCs)
to be introduced in Section 2. 
For each family of polynomial solutions of \refe{qEHT1} 
we search for the ones that are orthogonal in a suitable intervals 
depending on the range of the parameters coming from the coefficients of \refe{qEHT1}
and the corresponding $q$-Pearson equation. Hence, 
in Section 2, we present the candidate intervals
by inspecting the BCs as well as some preliminary results.
Theorems which help to calculate $q$-weight functions are given in Section 3.
Section 4 deals with the qualitative analysis including
the theorems stating the main results of our article. 
The last Section concludes the paper with some final remarks.

\section{The orthogonality and preliminary results}


We first introduce the so-called $q$-Jackson integrals
and afterward a well known theorem 
for the orthogonality of polynomial
solutions of \refe{qEHT1} in order to make
the article self-contained \cite{ran, cbook, NSU}.

The $q$-Jackson integrals for $q\in(0,1)$ \cite{GR, kac} are defined by
\begin{equation}\label{q-jac}
\int_0^af(x)d_qx =(1-q)a\sum_{j=0}^{\infty}q^j f(q^ja) \quad {\rm and} \quad 
\int_a^0 f(x)d_qx =(1-q)(-a) \sum_{j=0}^{\infty}q^j f(q^ja)
\end{equation}
if $a>0$ and $a<0$, respectively.
Therefore, we have
\begin{equation}\label{q-jac-a<0}
\int_a^bf(x)d_qx:=\int_0^bf(x)d_qx - \int_0^af(x)d_qx \quad {\rm and} \quad
\int_a^bf(x)d_qx:=\int_a^0f(x)d_qx + \int_0^bf(x)d_qx
\end{equation}
when $0<a<b$ and $a<0<b$, respectively. Furthermore, 
we make use of the \textit{improper} $q$-Jackson integrals
\begin{equation}\label{q-jac-imp}
\int_0^ \infty f(x)d_qx =(1-q)\sum_{j=-\infty }^{\infty}q^j f(q^j)\;\; {\rm and}\;
\int_{-\infty}^ \infty f(x)d_qx =(1-q)\sum_{j=-\infty }^{\infty}q^j [f(q^j)+f(-q^j)]
\end{equation}
where the second one is sometimes called the \textit{bilateral} $q$-integral.
The \textit{$q^{-1}$-Jackson integrals} are defined similarly. For 
instance, the improper \textit{$q^{-1}$-Jackson integral} on $(a, \infty)$
is given by
\begin{equation}\label{q-jac-var}
\int_a^{\infty}f(x)d_{q^{-1}}x =(q^{-1}-1)a\sum_{j=0}^{\infty}q^{-j}f(q^{-j}a), \quad a>0
\end{equation}
provided that $\lim_{j\to\infty}q^{-j}f(q^{-j}a)=0$ and the series is convergent.

\begin{theorem}\label{T1}
Let $\rho$ be a function satisfying the $q$-Pearson equation (\ref{openqpearson}) in such a way that
the BCs
\begin{equation}\label{bc1}
\sigma_1(x, q)\rho(x, q)x^k\bigg|_{x=a, b}=\sigma_2(q^{-1}x, q)\rho(q^{-1}x, q)x^k\bigg|_{x=a, b}=0,
\qquad k\in\mathbb{N}_0
\end{equation}
also hold. Then the sequence $\{P_n(x,q)\}$ of polynomial solutions of \refe{qEHT1} 
are orthogonal on $(a, b)$ w.r.t $\rho(x, q)$ in the sense that
\begin{equation}\label{qorthogonality}
\int_a^bP_n(x, q)P_m(x, q)\rho(x, q)d_qx=d_n^2(q)\delta_{mn},
\end{equation}
where $d_n(q)$ and $\delta_{mn}$ denote the norm of 
the polynomials $P_n$ and the Kronecker delta, respectively. 
Analogously, if the conditions
\begin{equation}\label{bc3}
\sigma_2(x, q)\rho(x, q)x^k\bigg|_{x=a, b}=\sigma_1(qx, q)\rho(qx, q)x^k\bigg|_{x=a, b}=0,
\qquad k\in\mathbb{N}_0
\end{equation}
are fulfilled, the $q$-polynomials then satisfy the relation
\begin{equation}\label{q-1orthogonality}
\int_a^bP_n(x, q)P_m(x, q)\rho(x, q)d_{q^{-1}}x=d_n^2(q)\delta_{mn}.
\end{equation}
\end{theorem}

\begin{rem}\label{rem2.2}
The relation
(\ref{qorthogonality})
means that the polynomials $P_n(x,q)$ 
are orthogonal with respect to a measure supported
on the set of points $\{q^ka\}_{k\in\mathbb{N}_0}$ and $\{q^kb\}_{k\in\mathbb{N}_0}$.
Since we are interested in the positive definite cases,
i.e., when $\rho(x,q)>0$, then,
\begin{itemize}
\item
when $a=0$, the measure is supported on the set of points 
$\{q^kb\}_{k\in\mathbb{N}_0}$ in $(0,b]$.
\item
when $a>0$, the measure should be  
supported on the finite set of points
$\{q^kb\}_{k=0}^N$ being $a=q^{N+1}b$.
\item
when $a<0$, the measure is supported on the  
set $\{q^ka\}_{k\in\mathbb{N}_0}\bigcup\{q^kb\}_{k\in\mathbb{N}_0}$ in $[a,0)\bigcup(0,b]$.
\end{itemize}
A similar analysis can be done for the relation
(\ref{q-1orthogonality}). 

\end{rem}

According to Theorem \ref{T1}, we have to determine an interval $(a, b)$ in which $\rho$
is $q$-integrable and $\rho>0$
on the lattice points of the types 
$\alpha q^{\pm k}$ and $ \beta q^{\pm k}$ for $k\in\mathbb{N}_0$ and $\alpha,\beta\in\RR$.
Such a weight function will be a solution of the $q$-Pearson equation  \refe{openqpearson}.
To this end, a qualitative analysis of the $q$-Pearson equation is presented by a detailed
inspection of the r.h.s. of \refe{openqpearson}.
Note that the r.h.s. of \refe{openqpearson} consists of
the polynomial coefficients $\sigma_1$ and $\sigma_2$ of the $q$-EHT which can be made
definite for possible forms of the coefficients. As a result, the possible behavior of $\rho$
on the left hand side (l.h.s.) of \refe{openqpearson} and the
candidate intervals can be obtained accordingly.

\subsubsection*{OPSs on finite $(a, b)$ intervals}
First assume that $(a, b)$ denotes a finite interval. 
We list the following possibilities for finding $\rho$
which obeys the BCs in \refe{bc1} or in \refe{bc3}.  

\medskip \noindent \textbf{PI.}  This is the simplest case where
$\sigma_1$ vanishes at both $x=a$ and $b$, i.e., 
$\sigma_1(a,q)=\sigma_1(b,q)=0$. Using \refe{openqpearson} rewritten of the
form 
\begin{equation}\label{openq-1pearson-equiv}
{\rho(q^{-1}x, q)}=\frac{q \sigma_1(x, q)}{\sigma_2(q^{-1}x, q)} {\rho(x, q)}
\end{equation}
we see that the function $\rho(x, q)$ becomes zero at the points 
$q^{-k}a$ and $q^{-k}b$ for $k\in\mathbb{N}$.
However, we have to take into consideration
three different situations.

\medskip \noindent (i) Let $a<0<b$. 
Since the points $q^{-k}a$ and $q^{-k}b$ lie outside the interval
$(a,b)$ and BCs are fulfilled at $x=a$ and $b$, 
there could be an OPS 
w.r.t. a measure supported on
the union of the set of points  $\{aq^{k}\}_{k\in\mathbb{N}_0}$ and $\{bq^{k}\}_{k\in\mathbb{N}_0}$
in $[a,0)\cup(0,b]$, if $\rho$ is positive.

\medskip \noindent (ii) Let $0<a<b$. In this case $\rho(x, q)$ vanishes at the points $q^{-k}a$ in $(a,b)$ and
$q^{-k}b$ out of $(a,b)$. Then, the only possibility to have an OPS on $(a,b]$ depends on the existence of $N$ such that $q^{N+1}b=a$.  This condition, however, implies 
that $bq^{k}=aq^{-(N-k)}$ and that $\rho$ vanishes at $bq^{k}$ for $k=0,1,\ldots,N$, and, 
therefore, it must be rejected. The similar statement is true when $a<b<0$, which 
can be obtained by a simple linear scaling transformation so that it 
does not represent an independent case.

\medskip \noindent (iii) Let $a=0<b$ (or, $a<b=0$).  
This case is much more involved. First of all, if $a=0$ is a zero of
$\sigma_1(x,q)$ then it is a zero of $\sigma_2(x,q)$ as well, both containing a factor $x$. 
Therefore, the r.h.s. of $q$-Pearson equation \refe{openqpearson}
can be simplified and {\bf PI}(i)-(ii)
are not valid anymore. In fact, in this case an 
OPS w.r.t. a measure supported on the set of points $\{bq^{k}\}_{k\in\mathbb{N}_0}$ in $(0,b]$ can be defined.

\medskip \noindent \textbf{PII.} 
The relation in \refe{ide-sr}
suggest an alternative possibility to define an OPS on $(a,b)$.
Namely, if $q^{-1}a$ and $q^{-1}b$ are both zeros of $\sigma_2(x, q)$, 
by using \refe{openqpearson} rewritten of the form
\begin{equation}\label{openqpearson-equiv}
{\rho(qx, q)}=\frac{q^{-1}
\sigma_2(x, q)}{\sigma_1(qx, q)}{\rho(x, q)},
\end{equation}
it follows that  $\rho(x, q)$ vanishes at the points $q^{k}a$ and $q^{k}b$ for $k\in\mathbb{N}_0$.
Then two different situations appear depending on whether $a<0<b$ or $0<a<b$. In the first case,
 $\rho(x, q)=0$ at the points $q^{k}a$ and $q^{k}b$ for $k\in\mathbb{N}_0$ in $[a,b]$, which is not interesting.
In the second case, the $q^kb$ are in $[a, b]$ whereas the $q^ka$ remain out of $[a,b]$, 
so that we could have an OPS if  
there exists $N$ such that $q^{-N-1}a=q^{-1}b$. However, since $q^{-k}a=q^{N-k}b$, 
$\rho$ vanishes at the $q^{-k}a$ which are in $[a, b]$ as well.

\medskip
\noindent \textbf{PIII.} Let $q^{-1}a$ and $b$ be the roots of
$\sigma_2$ and $\sigma_1$, respectively. 
Then we see, from (\ref{openq-1pearson-equiv})
and (\ref{openqpearson-equiv}), that $\rho=0$ at
$x=q^{-k}b$ for $k\in\mathbb{N}$ and at $x=q^ka$ for $k\in\mathbb{N}_0$. That is, if  $a<0<b$, 
$\rho=0$ on $x\in(a, 0)$ 
and, therefore, an OPS can not be constructed on $(a, b)$ unless $a\rightarrow0^-$.
In this limiting case of $x\in(0, b]$, it can be possible to introduce
a desired weight function supported on the set  
$\{bq^{k}\}_{k\in\mathbb{N}_0}$.
If $0<a<b$, on the other hand, $\rho$ vanishes for $x<a$
and $x>b$. Thus
there could be an OPS  
w.r.t. a measure supported on the finite set of 
points  $\{q^{k}b\}_{k=0}^N$
provided that $q^{N+1}b=a$ for some finite $N$ integer. 
Alternatively, we can define an equivalent OPS w.r.t. a measure
supported on the equivalent finite set of points $\{q^{-k}a\}_{k=0}^N$ 
provided now that $q^{-N-1}a=b$, where $N$ is a finite integer. 
Note that in the limiting case of $a\rightarrow0^+$
the set of points $\{q^{k}b\}_{k\in\mathbb{N}_0}$ becomes infinity.

\medskip
\noindent \textbf{PIV.} Assume that $a$ and $q^{-1}b$ are the roots of
$\sigma_1(x, q)$ and $\sigma_2(x, q)$, respectively. Then, from (\ref{openq-1pearson-equiv})
and (\ref{openqpearson-equiv}), it follows that $\rho(x, q)$ vanishes at the points
$q^{-k}a$, $k\in\mathbb{N}$ and $q^kb$, $k\in\mathbb{N}_0$. 
So, if $a<0<b$, it is not 
possible to find a  weight function satisfying the BCs. Nevertheless, as in {\bf PIII}, 
in the limiting case of $b\rightarrow0^+$
an OPS w.r.t. a measure supported at the points $q^{k}a$, $k\in\mathbb{N}_0$ in $[a, 0)$ can be constructed.
If $0<a<b$, there is no possibility to introduce an OPS. Note that when $a=0<b$,
an OPS also does not exist.  

\medskip
\noindent 

\subsubsection*{OPSs on infinite intervals}

Assume now that $(a, b)$ is an infinite interval. 
Without any loss of generality, 
let $a$ be a finite number and 
$b\rightarrow\infty$.
In fact, the system on the
infinite interval $(-\infty, b)$
is not independent which may be transformed
into $(a, \infty)$ on replacing $x$ by $-x$.
Obviously one BC in (\ref{bc1}) reads as
$$
\lim_{b\rightarrow\infty}\sigma_1(b, q)\rho(b, q)b^k=0\quad \mbox{or}\quad
\lim_{b\rightarrow\infty}\sigma_2(b, q)\rho(b, q)b^k=0,
\quad k\in\mathbb{N}_0,
$$
and there are additional cases  for $x=a$.

\medskip
\noindent \textbf{PV.}  If $x=a\neq0$ is root of $\sigma_1(x, q)$ then, 
from (\ref{openq-1pearson-equiv}), $\rho(x, q)$ vanishes at the points 
$q^{-k}a$ for $k\in\mathbb{N}$ which are interior points of $(a, \infty)$ 
when $a>0$. Therefore there is no OPS on $(a, \infty)$ for $a>0$. 
If $a<0$ we can find a $q$-weight function in $[a,0)\cup(0,\infty)$ supported 
on the union of the sets $\{q^ka\}_{k\in\mathbb{N}_0}$ and $\{q^{\pm k}\alpha\}_{k\in\mathbb{N}_0}$ for arbitrary $\alpha>0$  
where $\alpha$ can be taken as unity.
If $a=0$, on the other hand, then a weight function 
in $(0,\infty)$ can be defined at the points 
$q^{\pm k}\alpha$ for arbitrary $\alpha>0$ and $k\in\mathbb{N}_0$.

\medskip
\noindent \textbf{PVI.} If $x=q^{-1}a$ is a root of $\sigma_2(x, q)$, as we have already 
discussed, $\rho$ is zero at $q^ka$ for $k\in\mathbb{N}_0$. 
Therefore, for $a>0$ a $q$-weight function can exist
in $(q^{-1}a, \infty)$ supported on the set of points $\{q^{-k}a\}_{k\in\mathbb{N}}$. An OPS does not exist if $a<0$. 
Finally, if $a=0$ it is possible to find a $\rho$ on $(0, \infty)$ supported at the points 
$q^{\pm k}\alpha$ for arbitrary $\alpha>0$ and $k\in\mathbb{N}_0$. 

\medskip
\noindent \textbf{PVII.}
Finally, we consider the possibility of satisfying the BC 
$$
\lim_{a\rightarrow -\infty}\sigma_1(a, q)\rho(a, q)a^k=0
$$
in the limiting case as $a\rightarrow-\infty$. If this condition
holds a weight function and, hence, an OPS 
w.r.t. a measure supported on the set of points $\{\pm q^{\pm k}\alpha\}_{k\in\mathbb{N}_0}$, for arbitrary 
$\alpha>0$, can be defined.

\medskip

The aforementioned considerations are expressible as a theorem.   

\begin{theorem}\label{thm*} 
Let $a_1(q)$, $b_1(q)$ and $a_2(q)$, $b_2(q)$
denote the zeros of $\sigma_1(x,q)$ and $\sigma_2(x,q)$, respectively.
Let $\rho$ be a bounded and non-negative function satisfying the 
$q$-Pearson equation \refe{openqpearson} as well as the BCs \refe{bc1} or \refe{bc3}.
Then $\rho$ is a desired weight function for the polynomial 
solutions $P_n(x,q)$ of \refe{qEHT1} only in the following cases:

\begin{itemize}
 \item[\textbf{1.}]  Let $a<0<b$, where $a=a_1(q)$ and $b=b_1(q)$. 
Then $\rho$ is supported on  $\{q^ka\}_{k\in\mathbb{N}_0}\bigcup\{q^kb\}_{k\in\mathbb{N}_0}$
and
\begin{equation}\label{qortho1}
\int_{a_1(q)}^{b_1(q)}P_n(x, q)P_m(x, q)\rho(x, q)d_qx=d_n^2(q)\delta_{mn}.
\end{equation}
where the $q$-Jackson integral is of type \refe{q-jac-a<0}.
\item[\textbf{2.}] Let $a=0<b$, where $b=a_1(q)$. 
Then $\rho$ is supported on the set of points 
$\{q^kb\}_{k\in\mathbb{N}_0}$ in $(0,b]$ and
\begin{equation}\label{qortho2}
\int_{0}^{a_1(q)}P_n(x, q)P_m(x, q)\rho(x, q)d_qx=d_n^2(q)\delta_{mn},
\end{equation}
where the $q$-Jackson integral is of type \refe{q-jac}.
 
\item[\textbf{3.}] Let $0<a<b$, where $a=a_2(q)$ and $b=q^{-1}a_1(q)$. 
Then $\rho$ is supported on the finite set of points $\{q^{-k}a\}_{k=0}^N$
when $q^{-N-1}a=b$ and
\begin{equation}\label{qortho4-1}
\int_{a_2(q)}^{q^{-1}a_1(q)=q^{-N-1}a_2(q)}P_n(x, q)P_m(x, q)\rho(x, q)d_{q^{-1}}x=d_n^2(q)\delta_{mn},
\end{equation}
which is
the finite sum of the form
\begin{eqnarray*}
\int_{a_2(q)}^{q^{-N-1}a_2(q)}[\cdot] d_{q^{-1}}x
=(1-q^{-1})a_2(q)\sum_{k=0}^{N}
P_n(q^{-k}a_2(q), q)P_m(q^{-k}a_2(q), q)\rho(q^{-k}a_2(q), q).
\end{eqnarray*}


\item[\textbf{4.}] Let $a=a_1(q)<0$ and $b\rightarrow\infty$. 
Then $\rho$ is supported on the set $\{q^ka\}_{k\in\mathbb{N}_0}\bigcup 
\{q^{\mp k}\alpha\}_{k\in\mathbb{N}_0}$ for arbitrary $\alpha>0$ and 
\begin{equation}\label{qortho7}
\int_{a_1(q)}^{\infty}P_n(x, q)P_m(x, q)\rho(x, q)d_qx:= 
\int_{a_1(q)}^{0}[\cdot]d_qx+
\int_{0}^{\infty}[\cdot]d_qx=d_n^2(q)\delta_{mn},
\end{equation}
where the first $q$-Jackson integral is of type \refe{q-jac} 
and the second one is of type \refe{q-jac-imp}, respectively.

\item[\textbf{5.}] Let $a=a_2(q)>0$ and $b\rightarrow\infty$.
Then $\rho$ is supported on the set of points $\{q^{-k}a\}_{k\in\mathbb{N}_0}$ 
in $[a, \infty)$ and
\begin{equation}\label{qortho8}
\int_{a_2(q)}^{\infty}P_n(x, q)P_m(x, q)\rho(x, q)d_{q^{-1}}x=d_n^2(q)\delta_{mn},
\end{equation}
where the $q^{-1}$-Jackson integral is of type \refe{q-jac-var}.

\item[\textbf{6.}] Let $a=0$ and $b\rightarrow\infty$. 
Then $\rho$ is supported on the set of points $\{q^{\pm k}\alpha\}_{k\in\mathbb{N}_0}$
for arbitrary $\alpha>0$ and
\begin{equation}\label{qortho9}
\int_0^{\infty}P_n(x, q)P_m(x, q)\rho(x, q)d_qx=d_n^2(q)\delta_{mn},
\end{equation}
where the $q$-Jackson integral is of type \refe{q-jac-imp}.

\item[\textbf{7.}] Let $a\rightarrow-\infty$ and $b\rightarrow\infty$.
Then $\rho$ is supported on the set of points $\{\mp q^{\pm k}\alpha\}_{k\in\mathbb{N}_0}$ 
for arbitrary $\alpha>0$ and
\begin{equation}\label{qortho10}
 \int_{-\infty}^{\infty}P_n(x, q)P_m(x, q)\rho(x, q)d_qx=d_n^2(q)\delta_{mn},
\end{equation}
where the bilateral $q$-Jackson integral is of type \refe{q-jac-imp}.  

\end{itemize}
\end{theorem}

Before starting our analysis, let us mention that in accordance with \cite{NM, MAJF, NUpaper2} 
the $q$-polynomials can be classified by means
of the degrees of the polynomial coefficients $\sigma_1$ and $\sigma_2$ and the fact 
that either $\sigma_1(0, q) \sigma_2(0, q)\not=0$ or $\sigma_1(0, q)=\sigma_2(0, q)=0$.
Therefore, we can define two classes, namely, the non-zero ($\emptyset$) and zero ($0$) classes
corresponding to the cases where $\sigma_1(0, q) \sigma_2(0, q)\neq0$ and 
$\sigma_1(0, q)=\sigma_2(0, q)=0$, respectively (this is a consequence of the fact that
$\sigma_2(0, q)=q\sigma_1(0, q)$, i.e., $\sigma_1$ and $\sigma_2$ both have the same constant terms).
In each class we 
consider all possible degrees of the polynomial coefficients $\sigma_1(x, q)$ and $\sigma_2(x, q)$ as shown
in \cite[page 182]{MAJF}. We follow the notation introduced in \cite{NM, MAJF}, i.e.,  
the statement $\emptyset$-Laguerre/Jacobi implies that deg$\sigma_2=1$, deg$\sigma_1=2$, and 
$\sigma_1(0, q)\sigma_2(0, q)\neq0$ and the statement $0$-Jacobi/Laguerre indicates that deg$\sigma_2=2$,
deg$\sigma_1=1$ and $\sigma_1(0, q)=\sigma_2(0, q)=0$.

In the following we use the Taylor polynomial expansion for the coefficients
\begin{eqnarray} \label{generalsigmaq12}
\tau(x, q)&=&\tau'(0, q)x+\tau(0, q),\, \tau'(0, q)\neq0,\nonumber\\
\sigma_1(x, q)&=&\hlf\sigma_1''(0, q)x^2+\sigma_1'(0, q)x+\sigma_1(0, q)=\hlf\sigma_1''(0, q)[x-a_1(q)][x-b_1(q)], \\
\sigma_2(x, q) &=&\hlf\sigma_2''(0, q)x^2+
\sigma_2'(0, q)x+\sigma_2(0, q)=\hlf\sigma_2''(0, q)[x-a_2(q)][x-b_2(q)].\nonumber 
\end{eqnarray}

\begin{theorem}[Classification of the OPS of the $q$-Hahn class \cite{NM, MAJF}]
All orthogonal polynomial solutions of the $q$-difference equations \refe{qEHT1} and \refe{qEHT2}  
can be classified as follows:
\begin{enumerate}
\item \textbf{$\emptyset$-Jacobi/Jacobi} polynomials
where $\sigma_2''(0,q)\neq0$ and $\sigma_1''(0,q)\neq0$
with $\sigma_1(0,q)\sigma_2(0,q)\neq0$.
\item \textbf{$\emptyset$-Jacobi/Laguerre} polynomials
where $\sigma_2''(0,q)\neq0$ and $\sigma_1''(0,q)=0$, $\sigma_1'(0,q)\neq0$
with \break $\sigma_1(0,q) \sigma_2(0,q)\neq0$.
\item \textbf{$\emptyset$-Jacobi/Hermite} polynomials
where $\sigma_2''(0,q)\neq0$, $\sigma_1''(0,q)=0$, $\sigma_1'(0,q)=0$ and $\sigma_1(0,q)\neq0$
with $\sigma_1(0,q)\sigma_2(0,q)\neq0$.
\item \textbf{$\emptyset$-Laguerre/Jacobi} polynomials
where $\sigma_2''(0,q)=0$, $\sigma_2'(0,q)\neq0$ and $\sigma_1''(0,q)\neq0$
with  \break $\sigma_1(0,q) \sigma_2(0,q)\neq0$.
\item \textbf{$\emptyset$-Hermite/Jacobi} polynomials
where $\sigma_2''(0,q)=0$, $\sigma_2'(0,q)=0$, $\sigma_2(0,q)\neq0$
and $\sigma_1''(0,q)\neq0$ with $\sigma_1(0,q)\sigma_2(0,q)\neq0$.
\item \textbf{$0$-Jacobi/Jacobi} polynomials
where $\sigma_2''(0,q)\neq0$, $\sigma_2'(0,q)\neq0$ and $\sigma_1''(0,q)\neq0$,
$\sigma_1'(0,q)\neq0$ with $\sigma_2(0,q)=\sigma_1(0,q)=0$.
\item \textbf{$0$-Jacobi/Laguerre} polynomials
where $\sigma_2''(0,q)\neq0$, $\sigma_2'(0,q)\neq0$ and $\sigma_1''(0,q)=0$, $\sigma_1'(0,q)\neq0$ with 
$\sigma_2(0,q)=\sigma_1(0,q)=0$.
\item \textbf{$0$-Bessel/Jacobi} polynomials
where $\sigma_2''(0,q)\neq0$, $\sigma_2'(0,q)=0$ and $\sigma_1''(0,q)\neq0$,
$\sigma_1'(0,q)\neq0$ with $\sigma_2(0,q)=\sigma_1(0,q)=0$.
\item \textbf{$0$-Bessel/Laguerre} polynomials
where $\sigma_2''(0,q)\neq0$, $\sigma_2'(0,q)=0$ and $\sigma_1''(0,q)=0$, $\sigma_1'(0,q)\neq0$ with $\sigma_2(0,q)=\sigma_1(0,q)=0$.
\item \textbf{$0$-Laguerre/Jacobi} polynomials
where $\sigma_2''(0,q)=0$, $\sigma_2'(0,q)\neq0$ and
$\sigma_1''(0,q)\neq0$, $\sigma_1'(0,q)\neq0$ with $\sigma_2(0,q)=\sigma_1(0,q)=0$.

\end{enumerate}

\end{theorem}


\section{The $q$-weight function}
 In the following sections we will discuss the solutions of the 
$q$-Pearson equation (\ref{openqpearson}) defined on the $q$-linear lattices  
enumerated in Remark \ref{rem2.2}. Since it is a linear difference 
equation on a given lattice, its solution can be uniquely determined
by the equation  (\ref{openqpearson}) and the boundary conditions
(for more details on the general theory of linear difference equations see e.g. 
 \cite[\S1.2]{Ela05}). In fact, the explicit form of a $q$-weight function 
 can be deduced by means of Theorem \ref{lemma1}. 

\begin{theorem}\label{lemma1}
Let $f$ satisfy the difference equation
\begin{equation}
\label{lemma1-help}
\frac{f(qx;q)}{f(x;q)}=\frac{a(x;q)}{b(x;q)},
\end{equation}
such that the limits
$\dst\lim_{x\to0}f(x;q)=f(0,q)$ and $\dst\lim_{x\to\infty }f(x;q)=f(\infty,q)$
exist, where $a(x;q)$ and $b(x;q)$ are definite functions. 
Then $f(x;q)$ admits the two $q$-integral representations 
\begin{equation}\label{explicit-f}
f(x, q)=f(0, q)\exp\left[{\dst\int_0^x\frac{1}{(q-1)t}\ln\left[\frac{a(t, q)}
{b(t, q)}\right]d_qt}\right]
\end{equation}
and
\begin{equation}
\label{explicit-f1}
f(x, q)=f(\infty, q)\exp\left[{\dst\int_x^{\infty}
\frac{1}{(1-q^{-1})t}\ln\left[\frac{a(t, q)}{b(t, q)}\right]d_{q^{-1}}t}\right]
\end{equation}
provided that the integrals are convergent.
\end{theorem}
\begin{Proof}
Taking the logarithms of both sides of \refe{lemma1-help}, multiplying by  
$1/(q-1)t$ and then integrating from $0$ to $x$, we have 
\begin{equation}\nonumber
\int_0^x\frac{1}{(q-1)t}\ln\left[\frac{f(qt, q)}{f(t, q)}\right]d_qt
=\int_0^x\frac{1}{(q-1)t}\ln\left[\frac{a(t, q)}{b(t, q)} \right]d_qt.
\end{equation}
The l.h.s. is expressible as
\begin{eqnarray*}
\int_0^x\frac{1}{(q-1)t}\ln\left[\frac{f(qt, q)}{f(t, q)}\right]d_qt
\!\!&=&\!\!\lim_{n\rightarrow \infty}\sum_{j=0}^n\left[\ln\left(f(q^jx, q)\right)
-\ln\left(f(q^{j+1}x, q)\right)\right] 
\\
&=&\ln\left[f(x, q)\right]-\ln\left[f(0, q)\right],
\end{eqnarray*}
which completes the proof, on using the fact that
$f(q^{n+1}x, q)\rightarrow f(0, q)$ as $n\rightarrow\infty$ for
$0<q<1$. The second representation \refe{explicit-f1} can be proven in a similar way. \qed
\end{Proof}

This theorem can be used to derive the $q$-weight functions for every
$\sigma_1$ and $\sigma_2$. However, here we  take into account the quadratic
coefficients leading to $\emptyset$-Jacobi/Jacobi and $0$-Jacobi/Jacobi cases. 
The results, some of which can be found in \cite {NM}, are stated by the next theorem. 

\begin{theorem}
In the $\emptyset$-Jacobi/Jacobi case, let 
$\sigma_1(x, q)$ and $\sigma_2(x, q)$ be of forms
\refe{generalsigmaq12}
in which $\sigma_1''(0, q)a_1(q)b_1(q)\neq0$ and 
$\sigma_2''(0, q)a_2(q)b_2(q)\neq0$. And let, in $0$-Jacobi/Jacobi case, 
$b_1(q)=b_2(q)=0$, $\sigma_1''(0, q)\neq0$ and $\sigma_2''(0, q)\neq0$.
Then a solution $\rho(x, q)$ of $q$-Pearson equation \refe{openqpearson}  
is expressible in the equivalent forms
shown in Table \ref{tab1}.
\end{theorem}

\begin{small}
\begin{table}[!ht]
\centering 
\caption {Expressions for the $q$-weight function $\rho(x, q)$ in the Jacobi/Jacobi cases}
\medskip
\label{tab1}
\renewcommand{\tabcolsep}{0.2cm}
\renewcommand{\arraystretch}{1.4}
{\small
\begin{tabular}{|l|}
\hline 
$\emptyset$-Jacobi/Jacobi case \\
\hline \\
1.\quad$\dst\frac{(a_1^{-1}qx,\, b_1^{-1}qx; q )_{\infty}}
{(a_2^{-1}x,\, b_2^{-1}x; q )_{\infty}}$   \\
2.\quad $\dst
\left|x\right|^{\alpha}\frac{(b_1^{-1}qx, a_2q/x; q )_{\infty}}{(a_1/x, b_2^{-1}x; q )_{\infty}}$ where 
$q^{\alpha}=\dst
\frac{q^{-2}\sigma_2''(0, q)b_2}{\sigma_1''(0, q)b_1}$ \\
\hline\hline
 $0$-Jacobi/Jacobi case\\
\hline\\
3.\quad  $\dst\left|x\right|^{\alpha}\frac{(a_1^{-1}qx; q)_{\infty}}
{(a_2^{-1}x; q)_{\infty}}$ where
$q^{\alpha}=\dst\frac{q^{-2}\sigma_2''(0, q)a_2} 
{\sigma_1''(0, q)a_1}$\\
4.\quad$\dst\left|x\right|^{\alpha}\sqrt{x^{\log_qx-1}}(qa_1^{-1}x, qa_2/x; q)_{\infty}$ where
$q^{\alpha}=\dst{\frac{q^{-2}\sigma_2''(0, q)}{\sigma_1''(0, q)a_1}}$\\
\hline
\end{tabular}}
\end{table}
\end{small}

\begin{Proof}
We start proving the first expression in Table \ref{tab1}. Keeping in mind that
$q^{-1}\sigma_2(0, q)=\sigma_1(0, q)$ and that $q^{-1}
\sigma_2''(0, q)a_2(q)b_2(q)=\sigma_1''(0, q)a_1(q)b_1(q)$
we have from (\ref{openqpearson}) 
\begin{equation}\label{rho1}
\frac{\rho(qx, q)}{\rho(x, q)}=\frac{q^{-1}\sigma_2''(0, q)
[x-a_2(q)][x-b_2(q)]}{\sigma_1''(0, q)[qx-a_1(q)][qx-b_1(q)]}=
\frac{[1-a_2^{-1}(q)x][1-b_2^{-1}(q)x]}{[1-a_1^{-1}(q)qx][1-b_1^{-1}(q)qx]}
\end{equation}
which gives
\[\rho(x, q)=\rho(0, q)\exp\left\{\int_0^x\frac{1}{(q-1)t}
\Big[\ln(1-a_2^{-1}t)+\ln(1-b_2^{-1}t)
-\ln(1-a_1^{-1}qt)-\ln(1-b_1^{-1}qt)\Big]d_qt\right\}\]
on using (\ref{explicit-f}).
By definition (\ref{q-jac}) of the $q$-integral, we first obtain
\[\rho(x, q)=\rho(0, q)\exp\left\{\ln{\prod_{k=0}^{\infty}(1-a_1^{-1}q^{k+1}x)
(1-b_1^{-1}q^{k+1}x)}-\ln{\prod_{k=0}^{\infty}(1-a_2^{-1}q^{k}x)(1-b_2^{-1}q^{k}x)}\right\}\]
and, therefore,
\begin{equation}\label{quadraticsigmas}
\rho(x, q)=\rho(0, q)\frac{(a_1^{-1}qx, b_1^{-1}qx; q)_{\infty}}{(a_2^{-1}x, b_2^{-1}x; q)_{\infty}},
\qquad \rho(0, q)\neq0.
\end{equation}
This implies that $a_1(q)q^{-1-k}$ and $b_1(q)q^{-1-k}$ for $k\geq0$ are zeros of $\rho$. 
Furthermore, $a_2(q)q^{-j}$ and $b_2(q)q^{-j}$ for $j\geq0$ stand for the simple poles of $\rho$.
Note here that $\rho(0, q)$ can be made unity,
and $a_1(q)$, $b_1(q)$, $a_2(q)$ and $b_2(q)$ are non-zero constants.
Therefore the solution in (\ref{quadraticsigmas}) is continuous
everywhere except at the simple poles.

To show 4., we rewrite 
the $q$-Pearson equation in the form
\begin{equation}\label{rho4}	
\frac{\rho(qx, q)}{\rho(x, q)}=\frac{ax[1-a_2(q)/x]}
{[1-a_1^{-1}(q)qx]}, 
\qquad a=\frac{q^{-2}\sigma_2''(0, q)}{\sigma_1''(0, q)a_1(q)}
\end{equation}
and assume that $\rho$ is a product of three functions 
$\rho(x, q)=f(x, q)g(x, q)h(x, q)$. Hence, if $f$, $g$ and $h$ are
solutions of 
\begin{equation}\label{fgh}
\frac{f(qx, q)}{f(x, q)}=ax, \quad
\frac{g(qx, q)}{g(x, q)}=\frac{1}{[1-a_1^{-1}(q)qx]}
\;\; {\rm and}\;\;
\frac{h(qx, q)}{h(x, q)}=\left[1-\frac{a_2(q)}{x}\right]
\end{equation}
respectively, then $\rho=fgh$ is a solution of \refe{rho4}. 
A solution of \refe{fgh} for $f(x,q)$ is of the form
$f(x, q)=\dst\left|x\right|^{\alpha}H^{(1)}(x)$,
which can be verified by direct substitution. Here, the function 
$H^{(\beta)}(x)=\sqrt{x^{\log_q^{x^{\beta}}-\beta}}$ with $\beta\neq0$
was first defined in \cite{NM}, and $\alpha\neq0$ is such that $q^{\alpha}=a$. 
The equation in \refe{fgh} for $g(x, q)$ can be solved in a way similar to that of \refe{rho1}.
So we find that
$g(x,q)=g(0, q)(a_1^{-1}qx; q)_{\infty}$,
where $g(0,q)=1$.
The expression \refe{explicit-f} is not suitable in finding
$h(x, q)$ which gives a divergent infinite product. Instead, 
we employ \refe{explicit-f1} so that the equation for $h(x, q)$ becomes
$$
\frac{h(q^{-1}x, q)}{h(x, q)}=\frac{1}
{\left[1-qa_2(q)/x\right]}
$$ 
whose solution is of the form
$h(x, q)=h(\infty, q)(qa_2/x; q)_{\infty}$,
where $h(\infty, q)$ can be taken again as unity without 
loss of generality. 
Clearly $h(x, q)$ is uniformly convergent in any compact subset 
of the complex plane that does not contain the point at the origin. Moreover, the product  
converges to an arbitrary constant $c$, which has been set to unity, as $x\rightarrow\infty$.
Thus,
\begin{equation*}
\rho(x,q)=f(x,q)g(x,q)h(x, q)=\left|x\right|^{\alpha}
\sqrt{x^{\log_qx-1}}(qa_1^{-1}x, qa_2/x; q )_{\infty}.
\end{equation*}

In order to obtain the expressions 2 and 3 in Table \ref{tab1} for the weight function 
we use the same procedure as before, but starting from the $q$-Pearson equation written in 
the forms
\begin{equation}\label{rho2}
\frac{\rho(qx, q)}{\rho(x, q)}=a\frac{[1-a_2(q)/x][1-b_2^{-1}(q)x]}{[1-a_1(q)q^{-1}/x][1-b_1^{-1}(q)qx]}, \quad
a=\frac{q^{-2}\sigma_2''(0, q)b_2(q)}{\sigma_1''(0, q)b_1(q)}
\end{equation}
and
\begin{equation}\label{rho3}
\frac{\rho(qx, q)}{\rho(x, q)}=a\frac{[1-a_2^{-1}(q)x]}{[1-a_1^{-1}(q)qx]}, \quad
a=\frac{q^{-2}\sigma_2''(0, q)a_2(q)}{\sigma_1''(0, q)a_1(q)},
\end{equation}
respectively. This completes the proof.
 \qed
\end{Proof}

\begin{rem} Notice that for getting the expressions of the weight function we have  
used the $q$-Pearson equation rewritten in different forms, namely  \refe{rho1},
\refe{rho4}, \refe{rho2} and \refe{rho3}, and different solution procedure in each case,
therefore, it is not surprising that $\rho$ has several equivalent representations 
displayed in Table \ref{tab1}. However, they all satisfy the same linear equation
and, therefore, they differ only by a multiplicative constant. 
\end{rem}

For the sake of the completeness, let us obtain the analytic representations of
$q$-weight functions satisfying \refe{openqpearson} for the other cases.
 
\begin{theorem}
Let $\sigma_1$ and $\sigma_2$ be polynomials of at most 2nd degree
in $x$ as the form \refe{generalsigmaq12}.
Then a solution $\rho(x, q)$ of $q$-Pearson equation \refe{openqpearson} 
for each $\emptyset$-Jacobi/Laguerre, 
$\emptyset$-Jacobi/Hermite, 
$\emptyset$-Laguerre/Ja\-co\-bi, 
$\emptyset$-Her\-mi\-te/Ja\-co\-bi, 
$0$-Jacobi/Laguerre, $0$-Bessel/Jacobi,
$0$-Bessel/Laguerre and $0$-Laguerre/Jacobi
case is expressible in the equivalent forms
shown in Table \ref{tab2}.

\begin{table}[!ht]
\centering 
\caption {Expressions for the $q$-weight function $\rho(x, q)$ for the other cases}\medskip
\label{tab2}
\renewcommand{\tabcolsep}{0.1cm}
\renewcommand{\arraystretch}{1}
{\small
\begin{tabular}{|l|l|}
\hline    
&  \\
$\emptyset$-Jacobi/Laguerre & $1.\, \dst\frac{(a_1^{-1}qx; q )_{\infty}}
{(a_2^{-1}x, \, b_2^{-1}x; q)_{\infty}}$ 
$2.\, \left|x\right|^{\alpha}\sqrt{x^{\log_qx-1}}\dst\frac{(a_2q/x, b_2q/x; q )_{\infty}}{(a_1/x; q)_{\infty}}$,
$q^{\alpha}=\dst\frac{\hlf\sigma_2''(0, q)q^{-2}}{\sigma_1'(0, q)}$\\
& $3.\, \left|x\right|^{\alpha}x^{\log_qx}(qa_1^{-1}x, qa_2/x, qb_2/x; q)_{\infty}, \quad
q^{\alpha}=-\frac{q^{-2}\hlf \sigma_2''(0, q)}{\sigma_1'(0, q)a_1}$\\
$\emptyset$-Jacobi/Hermite & $1. \,\dst\frac{1}{(a_2^{-1}x, b_2^{-1}x; q )_{\infty}}$
 $2.\, \left|x\right|^{\alpha}x^{\log_qx-1}\dst(a_2q/x,\, b_2q/x; q )_{\infty}$,
$q^{\alpha}=\dst\frac{\hlf\sigma_2''(0, q)q^{-1}}{\sigma_1(0, q)}$\\
$\emptyset$-Laguerre/Jacobi & $1. \, \dst\frac{(a_1^{-1}qx,\, b_1^{-1}qx; q )_{\infty}}
{(a_2^{-1}x; q )_{\infty}}$  $2.\, \left|x\right|^{\alpha}\dst\frac{(qa_2/x, qb_1^{-1}x; q)_{\infty}}{(a_1/x; q)_{\infty}}$,
$q^{\alpha}=-\frac{q^{-2}\sigma_2'(0, q)}{\hlf\sigma_1''(0, q)b_1}$\\
$\emptyset$-Hermite/Jacobi & $\quad(a_1^{-1}qx, b_1^{-1}qx; q )_{\infty}$\\
&  \\
\hline 
&  \\
$0$-Jacobi/Laguerre & $1.\, \dst\left|x\right|^{\alpha}\frac{1}{(a_2^{-1}x; q)_{\infty}}$,
$q^{\alpha}=-\dst\frac{q^{-2}\hlf\sigma_2''(0, q)a_2}{\sigma_1'(0, q)}$ 
$2.\, \dst\left|x\right|^{\alpha}\sqrt{x^{\log_qx-1}}(qa_2/x; q)_{\infty}$, 
$\dst{\frac{q^{-2}\hlf\sigma_2''(0, q)}{\sigma_1'(0, q)}}$\\
$0$-Bessel/Jacobi & $\quad\left|x\right|^{\alpha}\sqrt{x^{\log_qx-1}}(a_1^{-1}qx; q)_{\infty}$,
$q^{\alpha}=-\dst\frac{q^{-2}\hlf\sigma_2''(0, q)}{\hlf\sigma_1''(0, q)a_1}$ \\
$0$-Bessel/Laguerre & $\quad\left|x\right|^{\alpha}\sqrt{x^{\log_qx-1}}$,
$q^{\alpha}=\dst\frac{q^{-2}\hlf\sigma_2''(0, q)}{\sigma_1'(0, q)}$ \\
$0$-Laguerre/Jacobi & $\quad\left|x\right|^{\alpha}
(a_1^{-1}qx; q)_{\infty}$, $q^{\alpha}=-\dst\frac{q^{-2}\sigma_2'(0, q)}{\hlf\sigma_1''(0, q)a_1}$ \\
&  \\
\hline
\end{tabular}}
\end{table}
\end{theorem}

\begin{Proof}
The proof is similar to the previous one. That is, to obtain the 
second formula for the $\emptyset$-Jacobi/Laguerre family we rewrite the 
$q$-Pearson equation \refe{openqpearson} in the form
$$
\frac{\rho(qx, q)}{\rho(x, q)}=ax\frac{[1-a_2/x][1-b_2/x]}{[1-a_1q^{-1}/x]}, \quad
a=\frac{q^{-2}\hlf\sigma_2''(0, q)}{\sigma_1'(0, q)}
$$
and then apply the same procedure described in the proof of the previous theorem.  \qed
\end{Proof}


\section{The orthogonality of $q$-polynomials: the Jacobi/Jacobi cases}

The rational function on the r.h.s. of the $q$-Pearson equation
\refe{openqpearson} has been examined in detail. Since it is 
the ratio of two polynomials $\sigma_1$ and $\sigma_2$ of at most 
second degree, we deal with a definite rational function having 
at most two zeros and two poles.  
In the analysis of the unknown quantity $\rho(qx, q)/\rho(x, q)$ 
on the l.h.s. of \refe{openqpearson}, 
we sketch roughly its graph by using every possible form 
of the definite rational function in question. 
In particular, we split the $x$-interval into subintervals 
according to whether $\rho(qx, q)/\rho(x, q)<1$ or $\rho(qx, q)/\rho(x, q)>1$,
which yields valuable information about the monotonicity of $\rho(x, q)$.
Other significant properties of $\rho$ are provided by the 
asymptotes, if there exist any, of $\rho(qx, q)/\rho(x, q)$. 
A full analysis along these lines is sufficient for a complete
characterization of the orthogonal $q$-polynomials.
A similar characterization is made in a very recent book 
\cite{KLS} based on the three-term recursion and the Favard theorem. 

Here, in this section, we discuss only the cases in which both
$\sigma_1$ and $\sigma_2$ are of second degree, i.e., the $\emptyset$- and $0$-Jacobi/Jacobi
cases.

\subsection{The non-zero case}\label{s4.1}

Let the coefficients $\sigma_1$ and $\sigma_2$ be quadratic polynomials in $x$
such that $\sigma_1(0, q)\sigma_2(0, q)\neq0$. If $\sigma_1$ is
written in terms of its roots, i.e., $\sigma_1(x, q)=\hlf\sigma''_1(0, q)[x-a_1(q)][x-b_1(q)]$
then, from (\ref{sigmaq12}),
\[\sigma_2(x, q)=\Big[\hlf\sigma_1''(0, q)
+(1-q^{-1})\tau'(0, q)\Big]qx^2-\Big[\hlf\sigma_1''(0, q)(a_1+b_1)
-(1-q^{-1})\tau(0, q)\Big]qx+\hlf q\sigma_1''(0, q)a_1b_1\]
where $\hlf\sigma_1''(0, q)+(1-q^{-1})\tau'(0, q)\neq0$ by hypothesis.
Then $q$-Pearson equation (\ref{openqpearson}) takes the form

\begin{equation}\label{constantqpearson}
f(x,q):=\frac{\rho(qx, q)}{\rho(x, q)}=\frac{q^{-1}
\sigma_2(x, q)}{\sigma_1(qx, q)}=
\left[1+\frac{(1-q^{-1})\tau'(0, q)}{\hlf\sigma''_1(0, q)}\right]  
\frac{[x-a_2(q)][x-b_2(q)]}{[qx-a_1(q)][qx-b_1(q)]}
\end{equation}
provided that the discriminant denoted by $\Delta_q$,
\[ \Delta_q:=\Big[a_1(q)+b_1(q)-\frac{(1-q^{-1})\tau(0, q)}{\hlf\sigma''_1(0, q)}\Big]^2-
4a_1(q)b_1(q)\Big[1+\frac{(1-q^{-1})\tau'(0, q)}{\hlf\sigma_1''(0, q)}\Big], \]
of the quadratic polynomial in the nominator of 
$f(x,q)$ in \refe{constantqpearson} is non-zero.
Here $x=a_2$ and $x=b_2$ denote the zeros of $f$,
and they are constant multiples of the roots of $\sigma_2(x,q)$.
 
We see that the lines $x=q^{-1}a_1$ and $x=q^{-1}b_1$ stand for 
the vertical asymptotes of $f(x,q)$ and 
the point $y=1$ is always its $y$-intercept since $\sigma_2(0, q)=q\sigma_1(0, q)$.
Moreover, the locations of the zeros 
of $f$ are determined by the straightforward lemma.

\begin{lemma}
Define the parameter
\begin{equation}\label{Lambda}
 \Lambda_q=\frac{1}{q^2}\left[1+\frac{(1-q^{-1})\tau'(0, q)}{\hlf\sigma''_1(0, q)}\right]\neq0
\end{equation}
so that the line $y=\Lambda_q$ denotes the horizontal asymptote of $f(x, q)$.
Then we encounter the following cases for the roots of the equation $f(x, q)=0$.

\medskip \noindent {\bf Case 1.} If $\Lambda_q>0$ and $a_1(q)<0<b_1(q)$, 
$f$ has two real and distinct roots with opposite signs. 

\medskip \noindent {\bf Case 2.} If $\Lambda_q>0$ and $0<a_1(q)<b_1(q)$, there exist three possibilities
 
\medskip {\bf(a)}  if $\Delta_q>0$, $f$ has two real roots with the same signs

\medskip {\bf(b)}  if $\Delta_q=0$, $f$  has a double root

\medskip {\bf(c)} if $\Delta_q<0$, $f$ has a pair of complex conjugate roots.

\medskip \noindent {\bf Case 3.} If $\Lambda_q<0$ and $a_1(q)<0<b_1(q)$, there exist three possibilities 

\medskip{\bf(a)} if $\Delta_q>0$, $f$ has two real roots with the same signs

\medskip {\bf(b)}  if $\Delta_q=0$, $f$  has a double root

\medskip {\bf(c)} if $\Delta_q<0$, $f$ has a pair of complex conjugate roots.

\medskip \noindent {\bf Case 4.}  If $\Lambda_q<0$ and  $0<a_1(q)<b_1(q)$, $f$ has two real distinct roots
 with opposite signs.
\end{lemma}

From \refe{constantqpearson} it is clear that we need to consider the cases 
$\Lambda_q>1$ and $0<\Lambda_q<1$ separately. 
Now, our strategy consists of sketching first the graphs of $f(x,q)$ 
depending on all possible relative positions of the zeros of $\sigma_1$ and $\sigma_2$. 
To obtain the behaviours of $q$-weight functions $\rho$ from the graphs of $f(x,q)=\rho(qx, q)/\rho(x, q)$, 
we divide the real line into subintervals in which $\rho$ is either
monotonic decreasing or increasing. We take into consideration only the subintervals
where $\rho>0$. Note that if $\rho$ is initially positive then we have $\rho>0$
everywhere in an interval where $\rho(qx, q)/\rho(x, q)>0$.
Then we find suitable intervals  
in cooperation with Theorem \ref{thm*}.

\begin{figure}[!htp]
\centering
\includegraphics{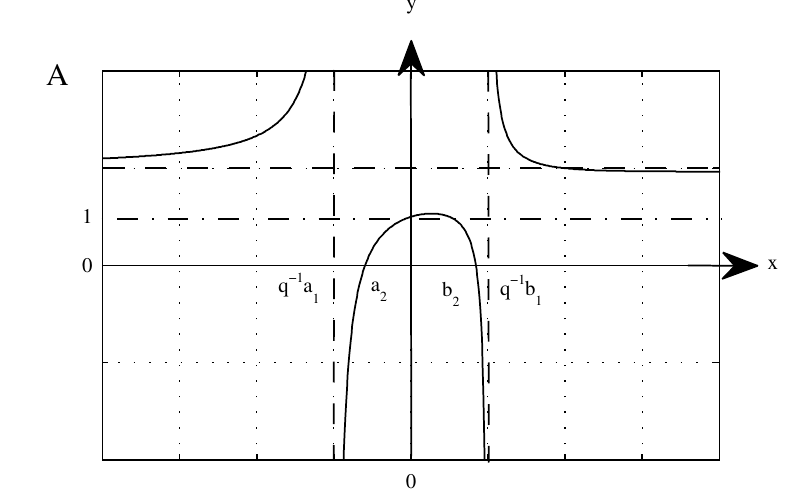}
\hfill
\includegraphics{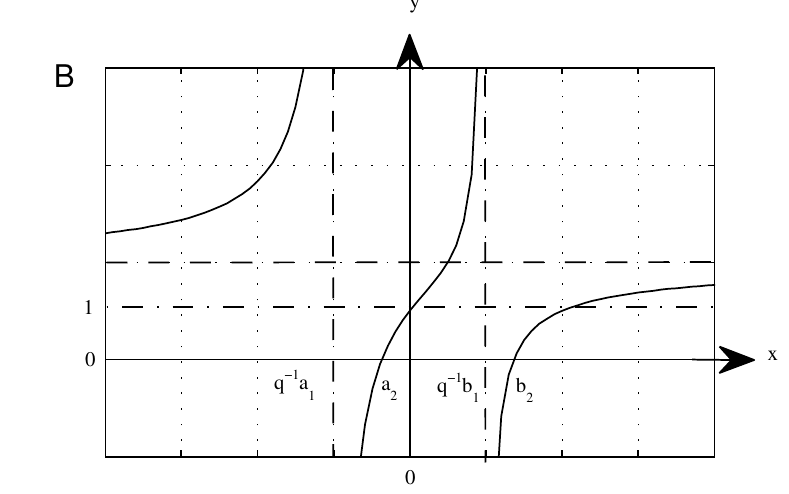}
\caption{The graph of $f(x, q)$ in {\bf Case 1} with $\Lambda_q>1$.
In A, the zeros are in order $q^{-1}a_1<a_2<0<b_2<q^{-1}b_1$,
and in B, $q^{-1}a_1<a_2<0<q^{-1}b_1<b_2$.}
\label{case1i}
\end{figure}

In Figure \ref{case1i}A, the intervals $(q^{-1}a_1, a_2)$ 
and $(b_2, q^{-1}b_1)$ are rejected immediately since 
$f$ is negative there. The subinterval $(a_2, b_2)$ should also be rejected
in which $\rho=0$ by {\bf PII}. For the same reason $(q^{-1}b_1, \infty)$ 
and $(-\infty, q^{-1}a_1)$, by symmetry, are not suitable by 
{\bf PV}. Therefore, an OPS fails to exist.

Let us analyse the problem presented in Figure \ref{case1i}B. The positivity of $\rho$ implies that
the intervals  $(q^{-1}a_1, a_2)$ and $(q^{-1}b_1, b_2)$ should be eliminated.
With the transformation $x=-t$, we eliminate also $(-\infty, q^{-1}a_1)$ by {\bf PV}.
The interval $(a_2, q^{-1}b_1)$ is not suitable too, by {\bf PIII}.
So it remains only $(b_2, \infty)$ to examine which coincides 
with the 5{\it th} case in Theorem \ref{thm*}.  
Since  $\rho(qx, q)/\rho(x, q)=1$ at $x_0=-\tau(0, q)/\tau'(0, q)>b_2(q)$,
then $\rho$ is increasing on $(b_2, x_0)$ and decreasing 
on $(x_0, \infty)$. 
As is shown from the figure $f$ has a finite limit 
as $x\rightarrow+\infty$ so that we could have the case 
$\rho\to0$ as $x\to\infty$. Even if $\rho\to0$ as $x\to\infty$,
we must show also that  
$\sigma_1(x, q)\rho(x, q)x^k\to0$ as $x\to\infty$ to satisfy the
BC. In fact, instead of the usual $q$-Pearson equation 
we have to consider the equation
\begin{equation}\label{openexpqpearson1}
g(x, q):=\frac{\sigma_1(qx, q)\rho(qx, q)(qx)^k}{\sigma_1(x, q)\rho(x, q)x^k}=q^k\frac{\sigma_1(x, q)
+(1-q^{-1})x\tau(x, q)}{\sigma_1(x, q)}=q^k\frac{q^{-1}
\sigma_2(x, q)}{\sigma_1(x, q)}
\end{equation} 
in case of an infinite interval, what we call it here the \textit{extended} $q$-Pearson equation
to determine the behaviour of the quantity
$\sigma_1(x, q)\rho(x, q)x^k$ as $x\to\infty$, which has been
easily derived from \refe{openqpearson}. It is obvious that 
the extended $q$-Pearson equation is the difference equation not for
the weight function $\rho(x,q)$ but for $\sigma_1(x, q)\rho(x, q)x^k$.

\begin{figure}[!htp]
\centering
\includegraphics{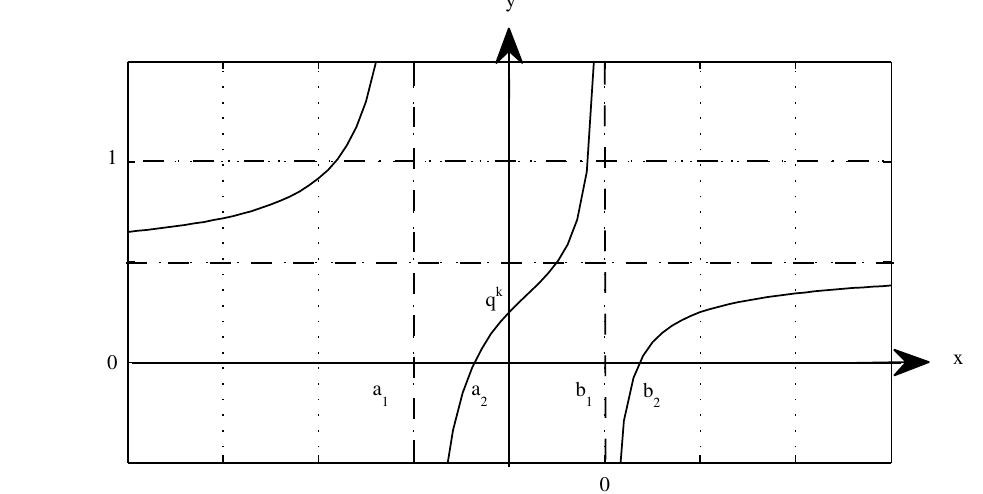}
\caption{The graph of $g(x, q)$ corresponding to Figure \ref{case1i}B.}
\label{case1iBexpqpearson}
\end{figure}

In Figure \ref{case1iBexpqpearson} we draw the graph of a typical
$g$ for some $0<q<1$, where $k$ is large enough.
From this figure we see that
$g<1$ for $x>b_2$ so that $\sigma_1(x, q)\rho(x, q)x^k$ does not vanish at $\infty$
since it is increasing as $x$ increases.
Thus we cannot find a weight function $\rho$ on $(b_2, \infty)$.

\begin{figure}[!htp]
\centering
\includegraphics{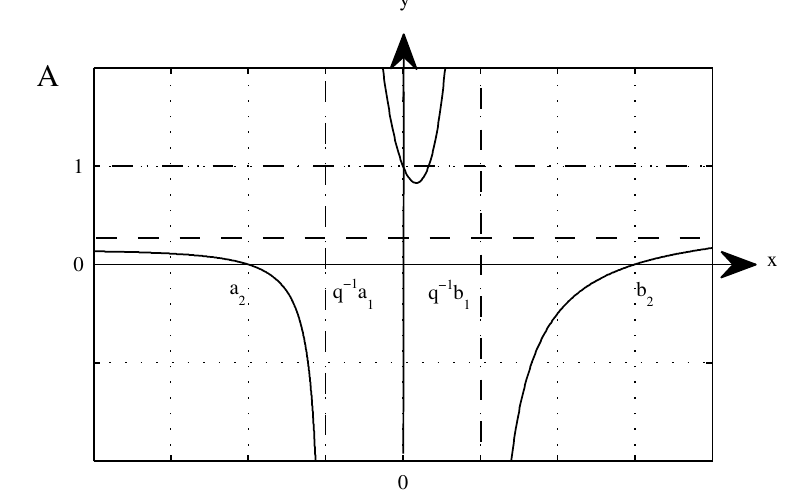}
\hfill
\includegraphics{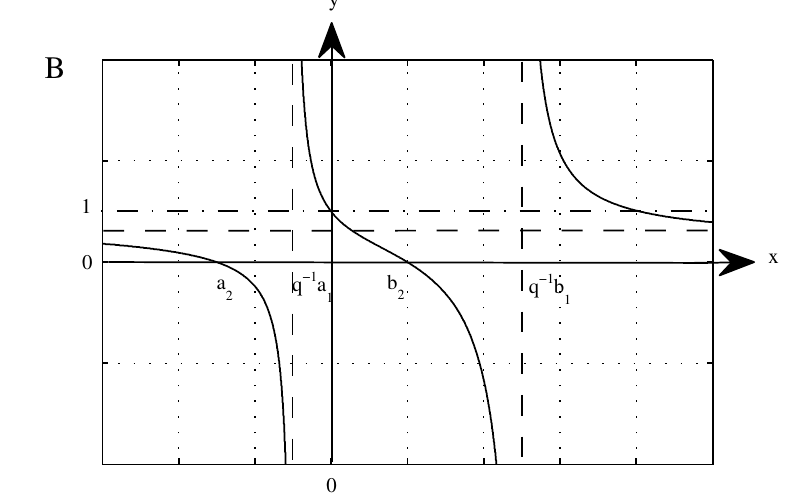}
\caption{The graph of $f(x, q)$ in {\bf Case 1} with $0<\Lambda_q<1$.
In A, the zeros are in order $a_2<q^{-1}a_1<0<q^{-1}b_1<b_2$
and in B, $a_2<q^{-1}a_1<0<b_2<q^{-1}b_1$.}
\label{case1ii}
\end{figure}

\begin{figure}[!htp]
\centering
\includegraphics{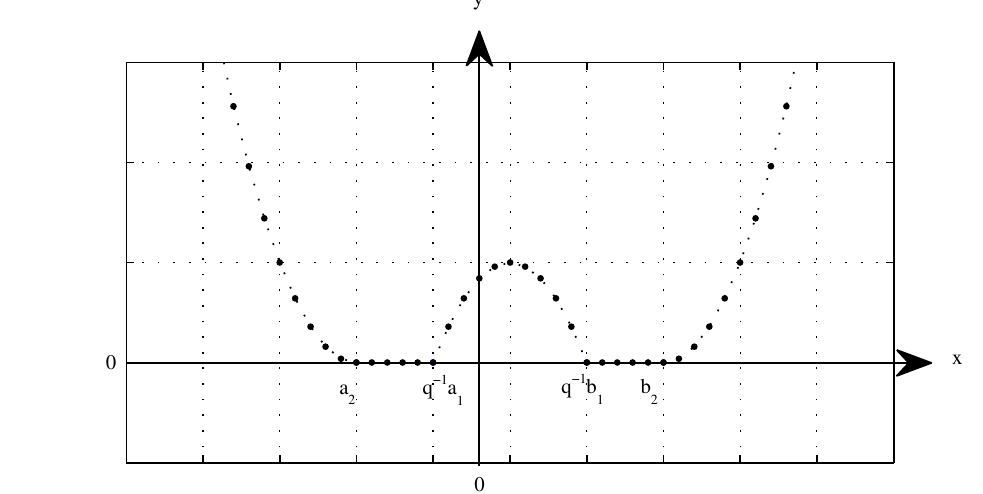}
\caption{The graph of $\rho(x, q)$ associated with the case in Figure \ref{case1ii}A.}
\label{case1iiArho}
\end{figure}

From Figure \ref{case1ii}A,
we first eliminate the intervals $(a_2, q^{-1}a_1)$ and $(q^{-1}b_1, b_2)$
because of the positivity of $\rho$.
The interval $(b_2, \infty)$ coincides again with the 5{\it th} case in
Theorem \ref{thm*}. However, $f(x, q)<1$ on this interval so that $\rho$ is
increasing on $(b_2, \infty)$ which implies that $\rho$ can not vanish as $x\to\infty$.
Thus $\sigma_1(x, q)\rho(x, q)x^k$ is never zero as $x\to\infty$ for some $k\in\mathbb{N}_0$.
The same is true for $(-\infty, a_2)$ by symmetry.
For the last subinterval $(q^{-1}a_1, q^{-1}b_1)$, we face 
the 1{\it th} case in Theorem \ref{thm*}.
Since  $\rho(qx, q)/\rho(x, q)=1$ at $q^{-1}a_1<x_0=-\tau(0, q)/\tau'(0, q)<q^{-1}b_1$,
then $\rho$ is increasing on $(q^{-1}a_1, x_0)$ and decreasing on 
$(x_0, q^{-1}b_1)$. Furthermore, $\rho(qx, q)/\rho(x, q)\to\infty$,
and hence $\rho\to0$, as $x\to q^{-1}a_1^+$ and
$x\to q^{-1}b_1^-$. As a result, the typical shape of $\rho$
is shown in Figure \ref{case1iiArho}
assuming a positive initial value of $\rho$ in each subinterval.
Then, an OPS
with such a weight function
in Figure \ref{case1iiArho} supported on
the union of set of points $\{q^ka_1(q)\}_{k\in\mathbb{N}_0}$ and $\{q^kb_1(q)\}_{k\in\mathbb{N}_0}$ exists
(see Theorem \ref{thm*}-\textbf{1}).
This OPS can be stated in the Theorem \ref{thmbigqjacobi}.

\begin{theorem}\label{thmbigqjacobi}
Consider the case where $a_2<a_1<0<b_1<b_2$ and $0<q^2\Lambda_q<1$.
Let $a=a_1(q)$ and $b=b_1(q)$ be zeros of $\sigma_1(x,q)$. 
Then there exists a sequence  of polynomials  $\{P_n\}$ for
$n\in\mathbb{N}_0$ orthogonal  
w.r.t. the weight function (see Eq. 1 in Table \ref{tab1})
\begin{equation}\label{rho-thmbigqjacobi}
\rho(x, q)=\frac{(qa^{-1}x, qb^{-1}x; q)_{\infty}}{(a_2^{-1}x, b_2^{-1}x; q)_{\infty}},
\end{equation}
supported on $\{q^ka\}_{k\in\mathbb{N}_0}\bigcup \{q^kb\}_{k\in\mathbb{N}_0}$,
(see \refe{qortho1} of Theorem \ref{thm*}-\textbf{1}).
\end{theorem}
The OPS in Theorem \ref{thmbigqjacobi} coincides with the 
case VIIa1 in Chapter 10 of  \cite[pages 292 and 318]{KLS}.  
In fact, a typical example of this family is the big $q$-Jacobi 
polynomials $P_n(x;a, b, c; q)$ satisfying the $q$-EHT with the coefficients
\[
\sigma_1(x, q)=q^{-2}(x-a_1)(x-b_1), \quad \sigma_2(x, q)=abq(x-a_2)(x-b_2),
\]
\begin{equation}\label{coeffbigqjacobi}
\tau(x, q)=\frac{1-abq^2}{(1-q)q}x+\frac{a(bq-1)+c(aq-1)}{1-q} \quad {\rm and} \quad
\lambda_n(q)=q^{-n}[n]_q\frac{1-abq^{n+1}}{q-1}
\end{equation}
where
$a_1=cq$, $b_1=aq$, $a_2=b^{-1}c$ and $b_2=1$. 
The conditions $0<q^2\Lambda_q<1$ and $a_2<a_1<0<b_1<b_2$ give the 
known constrains $c<0$, $0<b<q^{-1}$ and $0<a<q^{-1}$ on the parameters
of $P_n(x;a, b, c; q)$ with 
orthogonality on $\{cq,cq^2,cq^3,...\}\bigcup \{...,aq^3,aq^2,aq\}$ in the sense
\refe{qortho1} with
$$
d_n^2=(a-c)q(1-q)
\frac{(q, abq^2, a^{-1}cq, ac^{-1}q; q)_{\infty}}{(aq, bq, cq, abc^{-1}q; q)_{\infty}}
\dst\frac{(q, abq; q)_n}{(abq, abq^2; q)_{2n}}(aq, bq, cq, abc^{-1}q; q)_n(-ac)^nq^{n(n+3)/2}.
$$
It should be noted that the difference between these conditions
and those of Figure \ref{case1ii} comes from the fact that
we have considered not only the conditions
on $\rho$ but also on $\sigma_1\rho$ in Theorem \ref{thmbigqjacobi}.
Finally, the analysis of the case in Figure \ref{case1ii}B
does not yield an OPS.

\begin{figure}[!htp]
\centering
\includegraphics{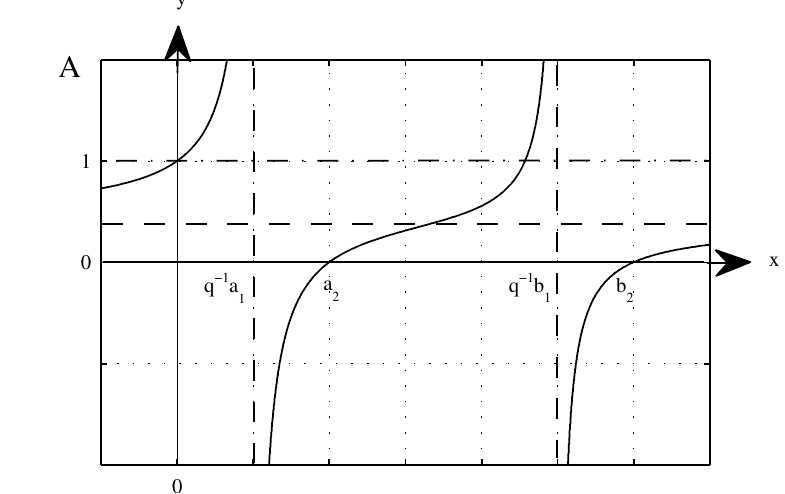}
\hfill
\includegraphics{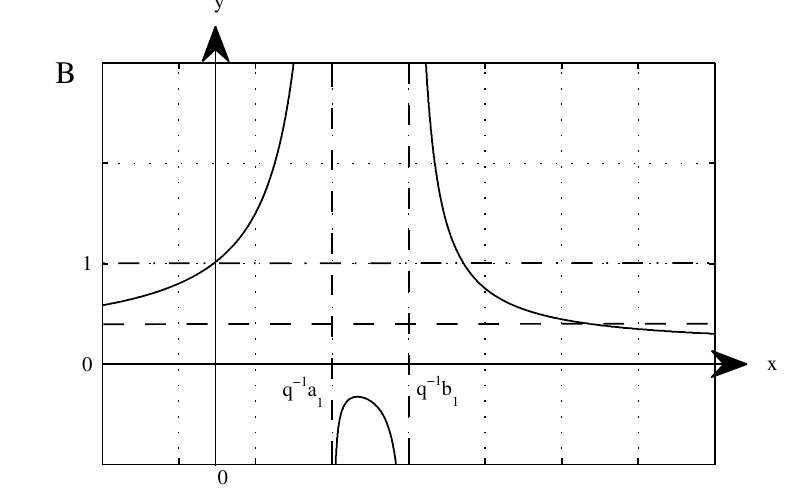}
\caption{The graph of $f(x, q)$ in {\bf Case 2} with $0<\Lambda_q<1$.
In A, we have {\bf Case 2(a)} with $0<q^{-1}a_1<a_2<q^{-1}b_1<b_2$. In B, 
we have {\bf Case 2(c)} with $0<q^{-1}a_1<q^{-1}b_1$ and $a_2, b_2\in\mathbb{C}$.}
\label{case2aciiEF}
\end{figure}

The case in Figure \ref{case2aciiEF}B is inappropriate to define an OPS.
On the other hand, in Figure \ref{case2aciiEF}A, the intervals
$(q^{-1}a_1, a_2)$ and $(q^{-1}b_1, b_2)$ are rejected by the positivity of $\rho$.
The intervals $(-\infty, q^{-1}a_1)$ and $(b_2, \infty)$ coincide with 4{\it th}, by symmetry,
and 5{\it th} cases in Theorem \ref{thm*}. However, $f(x, q)<1$ on $(-\infty, 0)$ 
and on $(b_2, \infty)$ so that $\rho$ is
decreasing on $(-\infty, 0)$ and increasing on $(b_2, \infty)$ 
which implies that $\rho$ can not vanish as $x\to-\infty$ and $x\to\infty$.
Thus $\sigma_1(x, q)\rho(x, q)x^k$ is never zero as $x\to-\infty$ and $x\to\infty$ 
for some $k\in\mathbb{N}_0$.
The only possible interval for the case in Figure \ref{case2aciiEF}A
is $(a_2, q^{-1}b_1)$
which corresponds to the 3{\it th} case of Theorem \ref{thm*}. Note that
$\rho(qx, q)/\rho(x, q)=1$ at $a_2<x_0=-\tau(0, q)/\tau'(0, q)<q^{-1}b_1$
and $\rho$ is increasing on $(a_2, x_0)$ and decreasing on 
$(x_0, q^{-1}b_1)$. Furthermore, $\rho(qa_2, q)=0$ and $\rho\to0$ as $x\to q^{-1}b_1^-$
since $\rho(qa_2, q)/\rho(a_2, q)$=0 and 
$\rho(qx, q)/\rho(x, q)\to\infty$ as $x\to q^{-1}b_1^-$, respectively.
It is clear that  
the BCs are satisfied at $x=a_2$ and $x=q^{-1}b_1$.  
Thus we can find a suitable $\rho$
on $[a_2, q^{-1}b_1)$ 
supported at the points $q^{-k}a_2$, $k=0,1,...,N$ where $q^{-N-1}a_2=q^{-1}b_1$.
Therefore, we state the following theorem.

\begin{theorem}\label{theoremcase2aiiE}
Consider the case where $0<a_1<a_2<b_1<b_2$ and $0<q^2\Lambda_q<1$.
Let $a=a_2(q)$ and $b=q^{-1}b_1(q)$ be zeros of $\sigma_2(x,q)$ 
and  $\sigma_1(qx,q)$, respectively.  
Then there exists a finite family of polynomials  $\{P_n\}$ 
orthogonal 
w.r.t. the weight function (see Eq. 2 in Table \ref{tab1}) 

\begin{equation}\label{rho-theoremcase2aiiE}
\rho(x, q)=\left|x\right|^{\iota}\frac{(\frac{qa}{x}, b^{-1}x; q)_{\infty}}{(\frac{a_1}{x}, b_2^{-1}x; q)_{\infty}},
\quad q^{\iota}=\frac{q^{-3} \sigma_2''(0, q)b_2}{\sigma_1''(0, q)b}
\end{equation}
supported on the set of points $\{q^{-k}a\}_{k=0}^N$ 
where $q^{-N-1}a=b$ (see \refe{qortho4-1} of Theorem \ref{thm*}-\textbf{3}).
\end{theorem}
The OPS in Theorem \ref{theoremcase2aiiE} corresponds to the 
case IIIb9 in Chapter 11 of \cite[page 366]{KLS}.  
A well known example of this family is the 
$q$-Hahn polynomials satisfying the $q$-EHT with the coefficients
\[
\sigma_1(x, q)=q^{-2}(x-a_1)(x-b_1), \quad \sigma_2(x, q)=\alpha\beta q(x-a_2)(x-b_2),
\]
\begin{equation}\label{coeffqhahn}
\tau(x, q)=\frac{1-\alpha\beta q^2}{(1-q)q}x+\frac{\alpha q^{-N}+\alpha\beta q-\alpha-q^{-N-1}}{1-q} \quad
{\rm and} \quad \lambda_n(q)=-q^{-n}[n]_q\frac{1-\alpha\beta q^{n+1}}{1-q}
\end{equation}
where $a_1=\alpha q$, $b_1=q^{-N}$, $a_2=1$ and $b_2=\beta^{-1}q^{-N-1}$.
The conditions $0<q^2\Lambda_q<1$ and $0<a_1<a_2<b_1<b_2$ give
the conditions $0<\alpha <q^{-1}$ and $0<\beta<q^{-1}$ on the parameters
of $Q_n(x; \alpha, \beta, N|q)$ with orthogonality on $\{1,q^{-1},q^{-2},...,q^{-N}\}$ in the sense
\refe{qortho4-1} where
\begin{equation}\label{dnqhahn}
d_n^2=\frac{(q, q^{N+1}, \beta^{-1}, \alpha^{-1}\beta^{-1}q^{-N-1}; q)_{\infty}}{(\alpha q,\beta q^{N+1}, \beta^{-1}q^{-N}, \alpha^{-1}\beta^{-1}q^{-1}; q)_{\infty}}\frac{(q, \alpha\beta q, \alpha q, q^{-N}, \beta q, \alpha\beta q^{N+2}; q)_n}{(\alpha\beta q, \alpha\beta q^2; q)_{2n}}(-\alpha q^{-N})^nq^{n(n+1)/2}(q^{-1}-1).
\end{equation}
In the literature, this relation is usually written as a finite sum \cite[page 367]{KLS}.

\begin{figure}[!htp]
\centering
\includegraphics{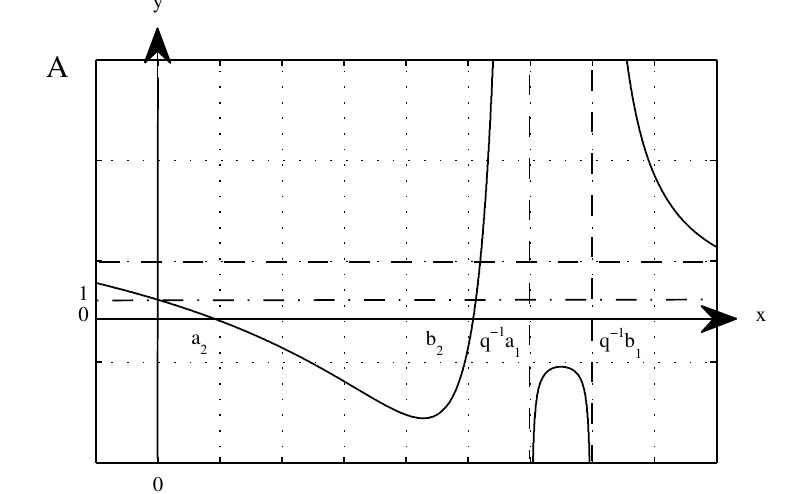}
\hfill
\includegraphics{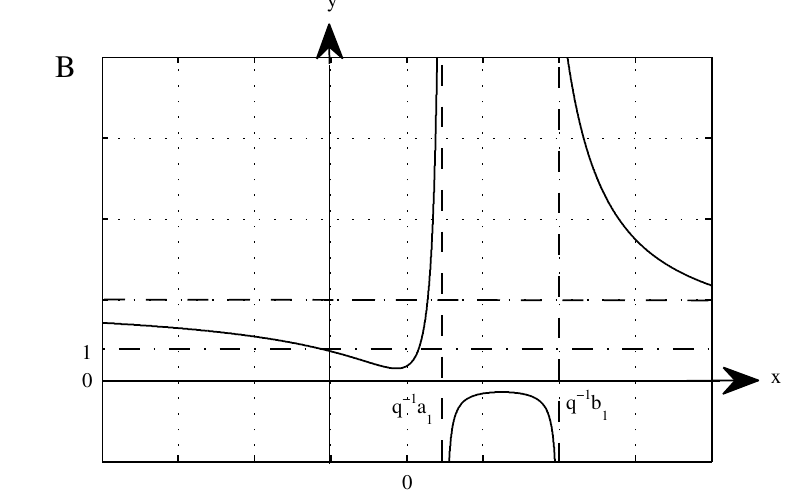}
\caption{The graph of $f(x, q)$ in {\bf Case 2} with $\Lambda_q>1$.
In A, we have {\bf Case 2(a)} with $0<a_2<b_2<q^{-1}a_1<q^{-1}b_1$. In B, 
we have {\bf Case 2(c)} with $0<q^{-1}a_1<q^{-1}b_1$ and
$a_2, b_2\in\mathbb{C}$.}
\label{case2aciEF}
\end{figure}

In Figure \ref{case2aciEF}A, the intervals $(a_2, b_2)$ 
and  $(q^{-1}a_1, q^{-1}b_1)$ are rejected by the 
positivity of $\rho$. We also eliminate the intervals 
$(-\infty, a_2)$ and $(q^{-1}b_1, \infty)$
due to {\bf PVI}, by symmetry, and  
{\bf PV}, respectively.
The last interval $[b_2, q^{-1}a_1)$ coincides with the 3{\it th} case
in Theorem \ref{thm*}. 
Notice that $f(x, q)=1$ at $b_2<x_0=-\tau'(0, q)/\tau(0, q)<q^{-1}a_1$,
then $\rho$ is increasing on $(b_2, x_0)$ 
with $\rho(qb_2, q)=0$ since $\rho(qb_2, q)/\rho(b_2, q)=0$ and decreasing on
$(x_0, q^{-1}a_1)$ with $\rho\to0$ as $x\to q^{-1}a_1^-$ 
since $\rho(qx, q)/\rho(x, q)\rightarrow\infty$. As a result, 
$[b_2, q^{-1}a_1)$ is an interval in which a desired $\rho$ is defined 
on the supporting points $q^{-k}b_2(q)$ for $k=1,2,...,N$ such that 
$q^{-N-1}b_2=q^{-1}a_1$.
Notice that the BC (\ref{bc3}) holds since $q^{-1}a_1$  and $b_2$
are one of the roots of $\sigma_1(qx, q)$ and $\sigma_2(x, q)$, respectively. 
As a significant remark, observe that the analysis is valid
in the limiting cases $a_1\rightarrow b_1$ and $a_2\rightarrow b_2$ as well.
Hence, the resulting OPS is presented in Theorem \ref{theoremcase2aiE}.
However, the case in  Figure \ref{case2aciEF}B does not give any OPS. 

\begin{theorem}\label{theoremcase2aiE}
Consider the case where $0<a_2\leq b_2<a_1\leq b_1$ and $q^2\Lambda_q>1$.
Let $a=b_2(q)$ and $b=q^{-1}a_1(q)$ be one of the zeros of 
$\sigma_2(x,q)$ and $\sigma_1(qx,q)$, respectively.  
Then there exists a finite family of polynomials $\{P_n\}$ 
orthogonal 
w.r.t. the weight function (see Eq. 2 in Table \ref{tab1})
\begin{equation*}
\rho(x, q)=\left|x\right|^{\iota}\frac{(b^{-1}x, \frac{qa}{x}; q )_{\infty}}{(\frac{b_1}{x}, a_2^{-1}x; q )_{\infty}}, 
\quad q^{\iota}=\dst
\frac{q^{-3}\sigma_2''(0, q)a_2}{\sigma_1''(0, q)b}
\end{equation*}
supported on the set of points $\{q^{-k}a\}_{k=0}^N$ 
where $q^{-N-1}a=b$ (see \refe{qortho4-1} of Theorem \ref{thm*}-\textbf{3}).
\end{theorem}
An example of this family is again the $q$-Hahn polynomials $Q_n(x; \alpha, \beta, N|q)$ 
with orthogonality on $\{1,q^{-1},q^{-2},...,q^{-N}\}$. 
They satisfy \refe{qEHT1} and \refe{qEHT2} with the coefficients 
$\sigma_1(x, q)=q^{-2}(x-a_1)(x-b_1)$, $\sigma_2(x, q)=\alpha\beta q(x-a_2)(x-b_2)$,
as in \refe{coeffqhahn}
but where now $a_1=q^{-N}$, $b_1=\alpha q$, $a_2=\beta^{-1}q^{-N-1}$ and $b_2=1$. 
These polynomials satisfy the orthogonality relation \refe{dnqhahn} with the same $d_n^2$
but with a different choice of the parameters 
$\alpha\geq q^{-N-1}$ and $\beta\geq q^{-N-1}$ which comes from the conditions
$0<a_2\leq b_2<a_1\leq b_1$ and $q^2\Lambda_q>1$. The authors did not mention this
different set of the parameters for the $q$-Hahn polynomials in \cite{KLS}.
However it is given in \cite[page 76]{ks}.

\begin{figure}[!htp]
\centering
\includegraphics{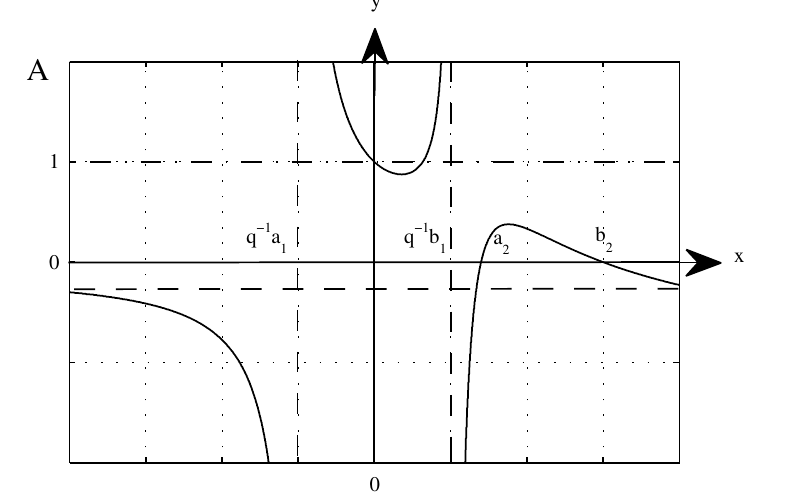}
\hfill
\includegraphics{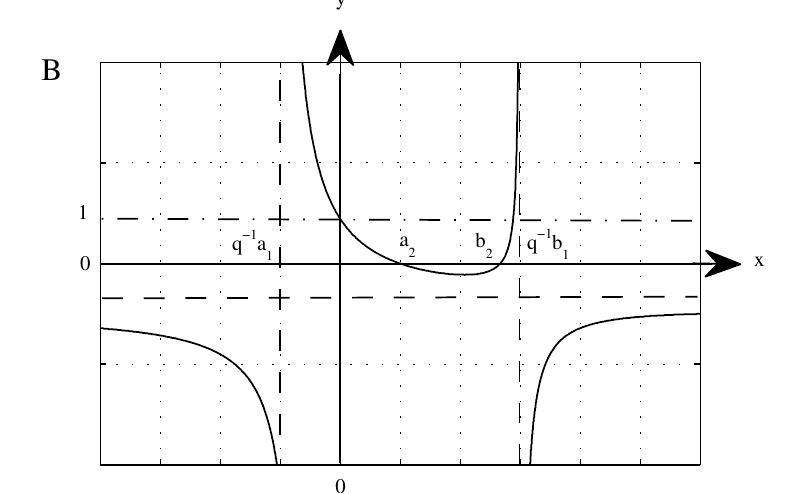}
\caption{The graph of $f(x, q)$ in {\bf Case 3(a)} with $\Lambda_q<0$.
In A, the zeros are in order $q^{-1}a_1<0<q^{-1}b_1<a_2<b_2$, and
in B, $q^{-1}a_1<0<a_2<b_2<q^{-1}b_1$.}
\label{case3aAB}
\end{figure}

In Figure \ref{case3aAB}A, 
the only suitable interval is  $(q^{-1}a_1, q^{-1}b_1)$ 
which coincides with the 1{\it st} case in Theorem \ref{thm*}.
In fact, $\rho(qx, q)/\rho(x, q)=1$ at $q^{-1}a_1<x_0=-\tau(0, q)/\tau'(0, q)<q^{-1}b_1$,
then it follows that $\rho$ is increasing on $(q^{-1}a_1, x_0)$ and decreasing on 
$(x_0, q^{-1}b_1)$ with $\rho\to0$ as $x\to q^{-1}a_1^+$ and $x\to q^{-1}b_1^-$
since $\rho(qx, q)/\rho(x, q)\to\infty$.
Notice that the BCs \refe{bc1} hold at 
$x=a_1$ and $x=b_1$. Then, there exists 
an OPS w.r.t. a $\rho$ 
supported on the set of points $\{q^ka_1\}_{k\in\mathbb{N}_0}\bigcup\{q^kb_1\}_{k\in\mathbb{N}_0}$. 
Thus, we have the following result.

\begin{theorem}\label{case3aA}
Consider the case $a_1<0<b_1<a_2\leq b_2$, and 
$q^2\Lambda_q<0$. Let $a=a_1(q)$ and $b=b_1(q)$ be zeros of $\sigma_1(x,q)$.  
Then there exists a sequence  of polynomials  $\{P_n\}$ for $n\in\mathbb{N}_0$ orthogonal  w.r.t. the weight function \refe{rho-thmbigqjacobi} 
supported on  $\{q^ka\}_{k\in\mathbb{N}_0}\bigcup\{q^kb\}_{k\in\mathbb{N}_0}$ (see \refe{qortho1} of Theorem \ref{thm*}-\textbf{1}).
\end{theorem}

An example of this family is again the big $q$-Jacobi polynomials which are
orthogonal on the set $\{cq,cq^{2},cq^{3},...\}\bigcup\{...,aq^3,aq^2,aq\}$. 
They satisfy the $q$-EHT with the coefficients in \refe{coeffbigqjacobi} where
$a_1=cq$, $b_1=aq$, $a_2=b^{-1}c$ and $b_2=1$. 
This case corresponds to the case VIIa1 in Chapter 10 of \cite[pages 292 and 318]{KLS}. 
However, notice that the conditions,  $a_1<0<b_1<a_2\leq b_2$ and 
$q^2\Lambda_q<0$, lead to the new constrains $c<0$, $b<0$, $abc^{-1}q\leq 1$ and $0<a<q^{-1}$,
which give a larger set of parameters for the orthogonality of the big $q$-Jacobi
polynomials than the one reported in \cite[page 319]{KLS}.

In Figure \ref{case3aAB}B, 
the only possible interval is $[b_2, q^{-1}b_1)$ 
which corresponds to the 3{\it th} case in Theorem \ref{thm*}.
In this case $\rho(qx, q)/\rho(x, q)=1$ at $b_2<x_0=-\tau(0, q)/\tau'(0, q)<q^{-1}b_1$,
then it follows that $\rho$ is increasing on $[b_2, x_0)$ and decreasing on 
$(x_0, q^{-1}b_1)$. Furthermore, $\rho(qb_2, q)=0$ and $\rho(x, q)\to0$ as $x\to q^{-1}b_1^-$
since $\rho(qb_2, q)/\rho(b_2, q)=0$ and $\rho(qx, q)/\rho(x, q)\to\infty$ as $x\to q^{-1}b_1^-$.
Thus there is a suitable $\rho$ 
supported on the set of points $\{q^{-k}b_2\}_{k\in\mathbb{N}_0}$
where $q^{-N-1}b_2=q^{-1}b_1$ (see Theorem \ref{thm*}-\textbf{3}). 
Hence we state the following theorem.
 
\begin{theorem}\label{case3aB}
Consider the case where $a_1<0<a_2\leq  b_2<b_1$ and 
$q^2\Lambda_q<0$.
Let $a=b_2(q)$ and $b=q^{-1}b_1(q)$ be zeros of $\sigma_2(x,q)$ and $\sigma_1(qx,q)$, respectively.  
Then there exists a finite family of polynomials $\{P_n\}$ 
 orthogonal w.r.t. the weight function (see Eq. 2 in Table \ref{tab1})
\begin{equation*}
\rho(x, q)=\left|x\right|^{\iota}\frac{(\frac{qa}{x}, b^{-1}x; q)_{\infty}}{(\frac{a_1}{x}, a_2^{-1}x; q)_{\infty}},
\quad q^{\iota}=\frac{q^{-3} \sigma_2''(0, q)a_2}{\sigma_1''(0, q)b}
\end{equation*}
supported on the set of points $\{q^{-k}a\}_{k\in\mathbb{N}_0}$ 
where $q^{-N-1}a=b$ (see \refe{qortho4-1} of Theorem \ref{thm*}-\textbf{3}).
\end{theorem}

A typical example of this family is again the $q$-Hahn polynomials 
orthogonal on $\{1, q^{-1},q^{-2},...,q^{-N}\}$. 
They satisfy \refe{qEHT1} and \refe{qEHT2} given by \refe{coeffqhahn}, 
being $a_1=\alpha q$, $b_1=q^{-N}$, $a_2=\beta^{-1}q^{-N-1}$ and $b_2=1$. The conditions $a_1<0<a_2\leq  b_2<b_1$ and 
$q^2\Lambda_q<0$
lead to the orthogonality relation for the $q$-Hahn polynomials that is valid
in a larger  set of the parameters, $\alpha<0$ and $\beta\geq q^{-N-1}$.
This new parameter set is not mentioned in \cite{KLS}.

\begin{figure}[!htp]
\centering
\includegraphics{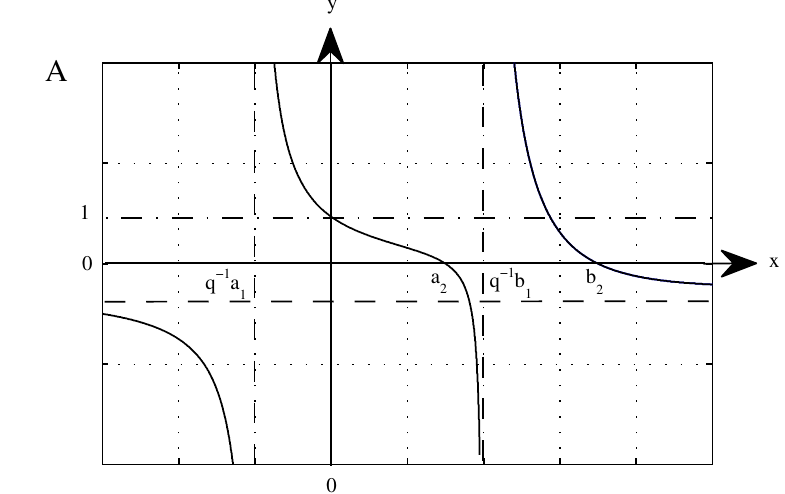}
\hfill
\includegraphics{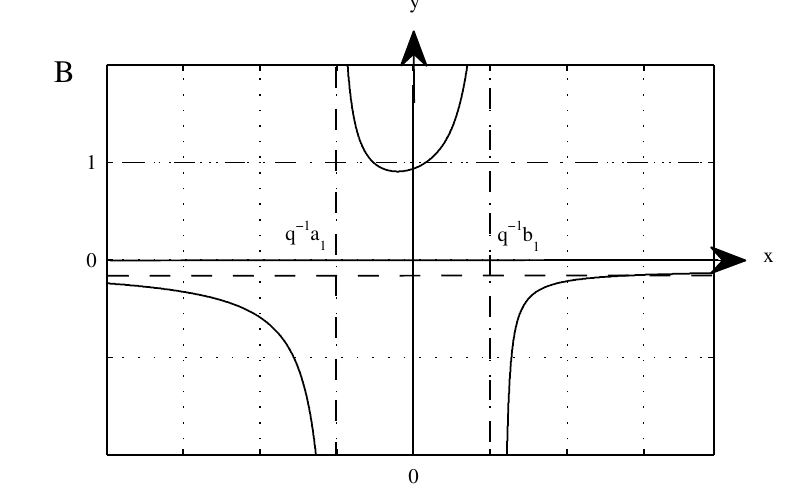}
\caption{The graph of $f(x, q)$ in {\bf Case 3} with $\Lambda_q<0$.
In A, we have {\bf Case 3(a)} with $q^{-1}a_1<0<a_2<q^{-1}b_1<b_2$. In B, 
we have {\bf Case 3(c)} with $q^{-1}a_1<0<q^{-1}b_1$ and $a_2, b_2\in\mathbb{C}$.}
\label{case3acCD}
\end{figure}

We could not find an OPS in case of Figure \ref{case3acCD}A.
In Figure \ref{case3acCD}B, the interval $(q^{-1}a_1, q^{-1}b_1)$ 
coincides with the 1{\it st} case in Theorem \ref{thm*}. Notice that $\rho(qx, q)/\rho(x, q)=1$ 
at $q^{-1}a_1<x_0=-\tau(0, q)/\tau'(0, q)<q^{-1}b_1$, then 
$\rho$ is increasing on $(q^{-1}a_1, x_0)$ and decreasing on 
$(x_0, q^{-1}b_1)$ with $\rho(x, q)\to0$ as $x\to q^{-1}a_1^+$ and $x\to q^{-1}b_1^-$
since $\rho(qx, q)/\rho(x, q)\to\infty$.
Then there
is a suitable $\rho$ supported on 
the set $\{q^ka_1\}_{k\in\mathbb{N}_0}\bigcup \{q^kb_1\}_{k\in\mathbb{N}_0}$.
Therefore,  the following theorem holds:

\begin{theorem}\label{theoremcase3cD}
Consider the case $a_1<0<b_1$, $a_2$, $b_2\in\mathbb{C}$ and 
$q^2\Lambda_q<0$.
Let $a=a_1(q)$ and $b=b_1(q)$ be zeros of $\sigma_1(x,q)$.
Then there exists a sequence  of polynomials  $\{P_n\}$ for $n\in\mathbb{N}_0$ orthogonal 
w.r.t. the weight function \refe{rho-thmbigqjacobi},
supported on $\{q^ka\}_{k\in\mathbb{N}_0}\bigcup \{q^kb\}_{k\in\mathbb{N}_0}$,
(see \refe{qortho1} of Theorem \ref{thm*}-\textbf{1}) with
\begin{eqnarray*}
d_n^2&=&\left(b_1-a_1\right)(1-q)q^{n(n-1)/2}\left(-a_1b_1\right)^n
\frac{(q, q^{-1}a_2^{-1}b_2^{-1}a_1b_1; q)_n}
{(q^{-1}a_2^{-1}b_2^{-1}a_1b_1, a_2^{-1}b_2^{-1}a_1b_1; q)_{2n}}\nonumber\\
&\times&(a_2^{-1}a_1, a_2^{-1}b_1, b_2^{-1}a_1, b_2^{-1}b_1; q)_n\frac{(q, qb_1a_1^{-1}, qa_1b_1^{-1}, a_2^{-1}b_2^{-1}a_1b_1; q)_{\infty}}
{(a_2^{-1}a_1, a_2^{-1}b_1, b_2^{-1}a_1, b_2^{-1}b_1; q)_{\infty}}
\end{eqnarray*}
where $q^2\Lambda_q=a_1b_1a_2^{-1}b_2^{-1}$, $a_2=i\alpha$, $b_2=\overline{a_2}=-i\alpha$, $\alpha\in\mathbb{R}$.
\end{theorem}
This case is included in the case VIIa1 in Chapter 10 of \cite[pages 292 and 318]{KLS}
(with $\gamma_2=\overline{\gamma_1}$) but it is not mentioned there. 
In fact, this case is similar to 
the big $q$-Jacobi polynomials studied in (\textbf{Cases 1} in Figure \ref{case1ii}A 
and \textbf{Case 3(a)} in Figure \ref{case3aAB}A). The difference 
is that the roots $a_2(q)$ and $b_2(q)$ are complex.

\begin{figure}[!htp]
\centering
\includegraphics{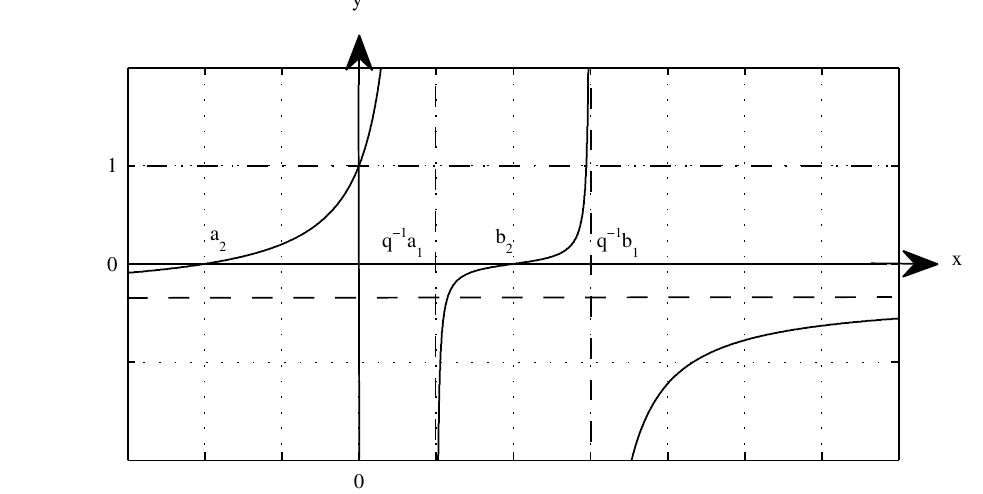}
\caption{The graph of $f(x, q)$ in {\bf Case 4} with $\Lambda_q<0$.
The zeros are in order $a_2<0<q^{-1}a_1<b_2<q^{-1}b_1$.}
\label{case4C}
\end{figure}

In Figure \ref{case4C},
the only possible interval is $[b_2, q^{-1}b_1)$ which corresponds to the one described 
in Theorem \ref{thm*}-\textbf{3}. A similar analysis shows that  there exists a $q$-weight function defined on the interval $[b_2, q^{-1}b_1)$ supported at the points
$q^{-k}b_2$ for $k=0,1,...,N$ where $q^{-N-1}b_2=q^{-1}b_1$ 
which lead to the following theorem:

\begin{theorem}\label{theoremcase4C}
Consider the case where $a_2<0<a_1<b_2<b_1$ and 
$q^2\Lambda_q<0$.
Let $a=b_2(q)$ and $b=q^{-1}b_1(q)$ be zeros of $\sigma_2(x,q)$ and $\sigma_1(qx,q)$, respectively.
Then there exists a finite family of polynomials $\{P_n\}$ 
 orthogonal w.r.t. the weight function \refe{rho-theoremcase2aiiE} in Theorem 
\refe{theoremcase2aiiE}
supported on the set of points $\{q^{-k}a\}_{k=0}^N$ 
where $q^{-N-1}a=b$ (see \refe{qortho4-1} of Theorem \ref{thm*}-\textbf{3}).
\end{theorem}

An example of this family is again the $q$-Hahn polynomials 
orthogonal on $\{1, q^{-1}, q^{-2},..., q^{-N}\}$. 
They satisfy \refe{qEHT1} and \refe{qEHT2} with
the coefficients \refe{coeffqhahn}
where $a_1=\alpha q$, $b_1=q^{-N}$, $a_2=\beta^{-1}q^{-N-1}$ and $b_2=1$.
The conditions $a_2<0<a_1<b_2<b_1$ and 
$q^2\Lambda_q<0$ lead to
another new constrains $0<\alpha <q^{-1}$ and $\beta<0$ on the parameters of the $q$-Hahn polynomials which extend the orthogonality relation for
the $q$-Hahn polynomials and it has been not reported in \cite{KLS}.

For completeness, we have also examined the cases listed below
for which an OPS fails to exist.

\medskip \noindent\textbf{Case 2(a)} with $0<a_2<q^{-1}a_1<b_2<q^{-1}b_1$ and $\Lambda_q>1$.

\medskip \noindent\textbf{Case 2(a)} with $0<a_2<q^{-1}a_1<q^{-1}b_1<b_2$ and $\Lambda_q>1$.

\medskip \noindent \textbf{Case 2(a)} with $0<q^{-1}a_1<a_2<b_2<q^{-1}b_1$ and $\Lambda_q>1$.

\medskip \noindent \textbf{Case 2(a)} with $a_2<b_2<0<q^{-1}a_1<q^{-1}b_1$ and $\Lambda_q>1$. 

\medskip \noindent \textbf{Case 2(a)} with $0<q^{-1}a_1<a_2<b_2<q^{-1}b_1$ and $0<\Lambda_q<1$.

\medskip \noindent \textbf{Case 2(a)} with $0<a_2<q^{-1}a_1<q^{-1}b_1<b_2$ and $0<\Lambda_q<1$.

\medskip \noindent \textbf{Case 2(a)} with $0<q^{-1}a_1<q^{-1}b_1<a_2<b_2$ and $0<\Lambda_q<1$.

\medskip \noindent \textbf{Case 2(a)} with $a_2<b_2<0<q^{-1}a_1<q^{-1}b_1$ and $0<\Lambda_q<1$.

\medskip \noindent \textbf{Case 4} with $a_2<0<q^{-1}a_1<q^{-1}b_1<b_2$ and $\Lambda_q<0$.

\medskip \noindent \textbf{Case 4} with $a_2<0<b_2<q^{-1}a_1<q^{-1}b_1$ and $\Lambda_q<0$.

\subsection{The zero case}
We make a similar analysis here with the same notations.
Let the coefficients $\sigma_1$ and $\sigma_2$ be quadratic polynomials in $x$ such that
$\sigma_1(0, q)=\sigma_2(0, q)=0$. If $\sigma_1$ is written as $\sigma_1(x, q)=
\hlf\sigma''_1(0, q)x[x-a_1(q)]$ 
then, from (\ref{sigmaq12}) $\sigma_2(x, q)=\hlf\sigma_2''(0, q)x^2+\sigma_2'(0, q)x$ where
$$
\hlf\sigma_2''(0, q)=q\left[\hlf\sigma_1''(0, q)+(1-q^{-1})\tau'(0, q)\right]\neq0 \quad
{\rm and} \quad \sigma_2'(0, q)=q(1-q^{-1})\tau(0, q)-\hlf\sigma_1''(0, q)a_1(q).
$$
Then it follows from (\ref{openqpearson}) that
\begin{eqnarray}\label{zerojacobiqpearson}
f(x, q):=\frac{\rho(qx, q)}{\rho(x, q)}&=&
\left[1+\frac{(1-q^{-1})\tau'(0, q)}{\hlf\sigma_1''(0, q)}\right]\frac{[x-a_2(q)]}{q[qx-a_1(q)]}, 
\qquad x\neq0
\end{eqnarray}
provided that $\left[1+\frac{(1-q^{-1})\tau'(0, q)}{\hlf\sigma_1''(0, q)}\right]a_2(q)=
\left[a_1(q)-\frac{(1-q^{-1})\tau(0, q)}{\hlf\sigma_1''(0, q)}\right]$.
Let us point out that $\Lambda_q$ defined in \refe{Lambda} is also
horizontal asymptote of $f(x, q)$ in \refe{zerojacobiqpearson}. Moreover, 
$f$ intersects the $y$-axis at the point 
$$
y:=y_0=q^{-1}\left[1-\frac{(1-q^{-1})}{a_1(q)}\frac{\tau(0, q)}{\hlf \sigma_1''(0,q)}\right].
$$
In the zero cases notice that one of the boundary of $(a, b)$ interval could be zero.
Therefore for such a case we need to know the behaviour of $\rho$ at the origin. 

\begin{lemma}\label{zeropoint} If $0<y_0<1$, then  $\rho(z, q)\to0$ as $z\to0$. Otherwise it 
diverges to $\mp\infty$.
\end{lemma}
\begin{Proof}
From \refe{openqpearson-equiv}, we write
$\rho(q^kx, q)=\rho(x, q)\prod_{i=0}^{k-1}\frac{q^{-1}\sigma_2(q^ix, q)}{\sigma_1(q^{i+1}x, q)}$
from which
\begin{equation}
\rho(q^kx, q)=q^{-k}\left[1-\frac{(1-q^{-1})}{a_1(q)}\frac{\tau(0, q)}
{\hlf\sigma_1''(0, q)}\right]^k\frac{(x/a_2(q); q)_k}{(qx/a_1(q); q)}_k\rho(x, q)
\end{equation}
is obtained. Taking $k\to\infty$ the result follows.
\qed
\end{Proof}

In a similar fashion, we introduce the two additional cases
which include all possibilities.

\medskip \noindent \textbf{Case 5.} $\Lambda_q>0$ with
{\bf(a)} $0<y_0<1$ or
{\bf(b)} $y_0>1$ or
{\bf(c)} $y_0<0$.

\medskip \noindent \textbf{Case 6.} $\Lambda_q<0$ with
{\bf(a)} $0<y_0<1$ or
{\bf(b)} $y_0>1$ or
{\bf(c)} $y_0<0$.

\begin{figure}[!htp]
\centering
\includegraphics{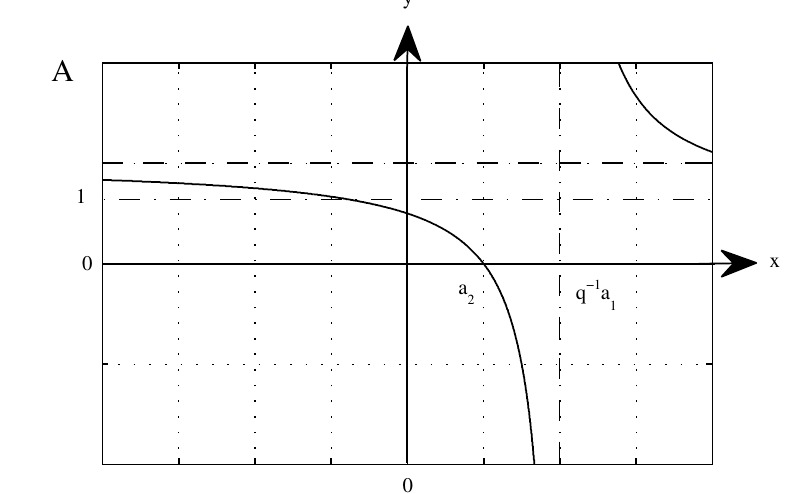}
\hfill
\includegraphics{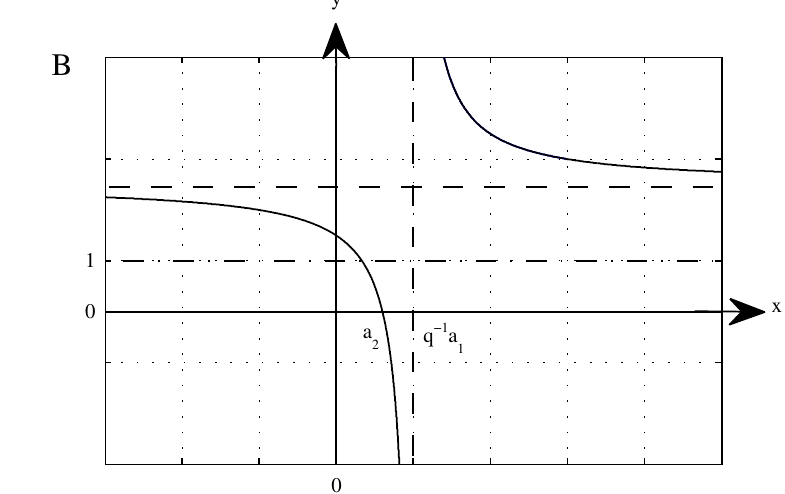}
\caption{The graph of $f(x, q)$ in \textbf{Case 5}
with $\Lambda_q>1$ and
$0<a_2<q^{-1}a_1$. In A, we have \textbf{Case 5(a)}. In B, we have \textbf{Case 5(b)}.}
\label{zerocase1iAB}
\end{figure}

In Figure \ref{zerocase1iAB}A,
we consider all possible intervals in which we can have a suitable 
$q$-weight function $\rho$.
By the positivity of $\rho$, the interval
$(a_2, q^{-1}a_1)$ is not suitable. The other intervals $(0, a_2)$
and $(q^{-1}a_1, \infty)$ are both eliminated due to the {\bf PIV}
and {\bf PV}, respectively. The interval $(-\infty, 0)$ is the one 
described in Theorem \ref{thm*}-\textbf{6} by symmetry. Since
$\rho(qx, q)/\rho(x, q)=1$ at $x_0=-\tau(0, q)/\tau'(0, q)<0$,
$\rho$ is increasing on $(-\infty, x_0)$ and decreasing on 
$(x_0, 0)$. Moreover, since $0<y_0<1$ $\rho\rightarrow0$ as $x\rightarrow0^-$ 
according to Lemma \ref{zeropoint}. On the other hand, since $\rho(qx, q)/\rho(x, q)$ has a finite limit 
as $x\rightarrow-\infty$, we may have $\rho\to0$ as $x\to-\infty$, but 
we should check that 
 $\sigma_1(x, q)\rho(x, q)x^k\to0$ as $x\to-\infty$ by using  the
extended $q$-Pearson equation (\ref{openexpqpearson1}). However,
the graph of $g$ \refe{openexpqpearson1} looks like the one represented 
in Figure \ref{zerocase1iAB}A
 together with the property that $g(x, q)<1$ on $(-\infty, 0)$ for large $k$ 
 which leads to that  $\sigma_1(x, q)\rho(x, q)x^k$ is decreasing on $(-\infty, 0)$
 with  $\sigma_1(x, q)\rho(x, q)x^k\not\to0$ as $x\to-\infty$. Therefore, this case 
 does not lead to any suitable $\rho$ and, therefore, OPS. The same result is valid 
for the case in Figure \ref{zerocase1iAB}B.

\begin{figure}[!htp]
\centering
\includegraphics{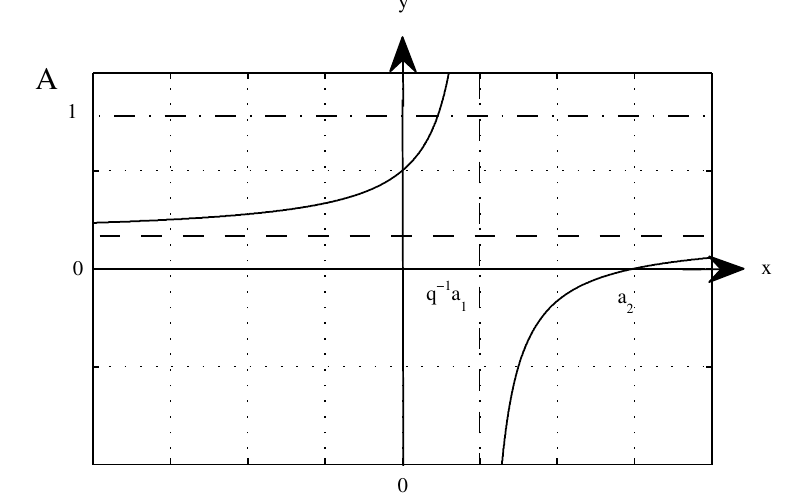}
\hfill
\includegraphics{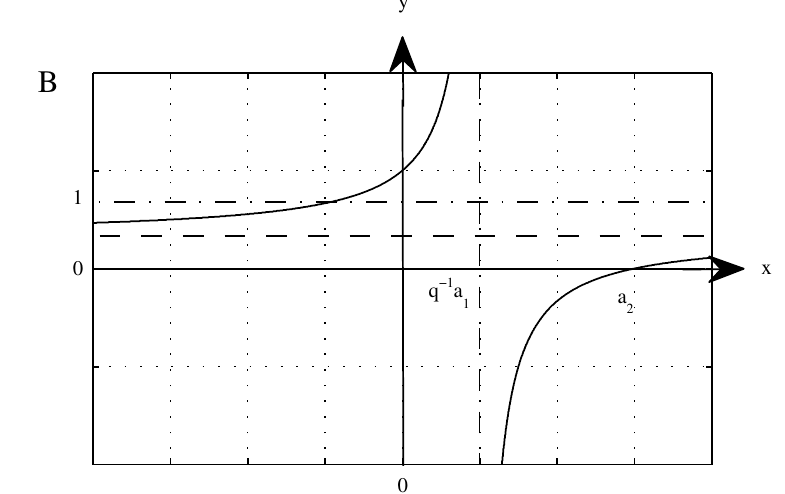}
\caption{The graph of $f(x, q)$ in \textbf{Case 5}
with $0<\Lambda_q<1$ and
$0<q^{-1}a_1<a_2$. In A, we have \textbf{Case 5(a)}. In B, we have \textbf{Case 5(b)}.}
\label{zerocase1iiGH}
\end{figure}

In Figure \ref{zerocase1iiGH}A,
the only suitable interval is $(0, q^{-1}a_1)$ which coincides 
with the 2{\it nd} case in Theorem \ref{thm*}.
Notice that $\rho(qx, q)/\rho(x, q)=1$ at $0<x_0=-\tau(0, q)/\tau'(0, q)<q^{-1}a_1$,
then $\rho$ is increasing on $(0, x_0)$ with $\rho\to0$ as $x\to0^+$ since
$0<y_0<1$ and decreasing on $(x_0, q^{-1}a_1)$ with $\rho(x, q)\to0$ as $x\to q^{-1}a_1^-$. 
Thus there is a a $q$-weight function supported at the points $q^ka_1$ for $k\in\mathbb{N}_0$ 
Hence, the resulting OPS is introduced in Theorem \ref{theoremzerocase1iiG}.
But, the case in Figure \ref{zerocase1iiGH}B does not yield any OPS.

\begin{theorem}\label{theoremzerocase1iiG}
Consider the case where $0<a_1<a_2$, $0<qy_0<1$ and $0<q^2\Lambda_q<1$.
Let $a=0$ and $b=a_1(q)$ be the zeros of $\sigma_1(x,q)$.
Then there exists a sequence  of polynomials  $\{P_n\}$ 
for $n\in\mathbb{N}_0$ orthogonal 
w.r.t. the weight function (see Eq. 3 in Table \ref{tab1})
\begin{equation*}
\rho(x, q)=\left|x\right|^{\alpha}\frac{(qb^{-1}x; q)_{\infty}}{(a_2^{-1}x; q)_{\infty}},
\quad q^{\alpha}=\frac{q^{-2} \sigma_2''(0, q)a_2}{ \sigma_1''(0, q)b}
\end{equation*}
supported on $\{q^kb\}_{k\in\mathbb{N}_0}$ (see \refe{qortho2} of Theorem \ref{thm*}-\textbf{2}).
\end{theorem}
The OPS in Theorem \ref{theoremzerocase1iiG} 
corresponds to the case IVa3 in Chapter 10 of \cite[pages 277 and 311]{KLS}.  
In fact, a typical example of this family is the little 
$q$-Jacobi polynomials $P_n(x; a, b|q)$  
satisfying the $q$-EHT with the coefficients
\[
\sigma_1(x, q)=q^{-2}x(x-a_1), \quad \sigma_2(x, q)=abqx(x-a_2),
\]
\begin{equation}\label{coefflittleqjacobi}
\tau(x, q)=\frac{1-abq^2}{(1-q)q}x+\frac{aq-1}{(1-q)q}\quad {\rm and} \quad \lambda_n(q)=-q^{-n}[n]_q\frac{1-abq^{n+1}}{1-q}
\end{equation}
where $a_1=1$ and $a_2=b^{-1}q^{-1}$.
The conditions $0<q^2\Lambda_q<1$, $0<qy_0<1$ and $0<a_1<a_2$ yield the 
restrictions $0<a<q^{-1}$ and $0<b<q^{-1}$ on the parameters of
$P_n(x; a, b|q)$ with orthogonality 
on $\{...,q^2,q,1\}$ in the sense \refe{qortho2} 
where
\begin{equation}\label{dnlittleqjacobi}
d_n^2=a^nq^{n^2}(1-q)\frac{(q, abq; q)_n}{(abq, abq^2; q)_{2n}}(aq, bq; q)_n
\frac{(q, abq^2; q)_{\infty}}{(aq, bq; q)_{\infty}}.
\end{equation}
In the literature, this relation can be found as an infinite sum \cite [Page 312]{KLS}.

\begin{figure}[!htp]
\centering
\includegraphics{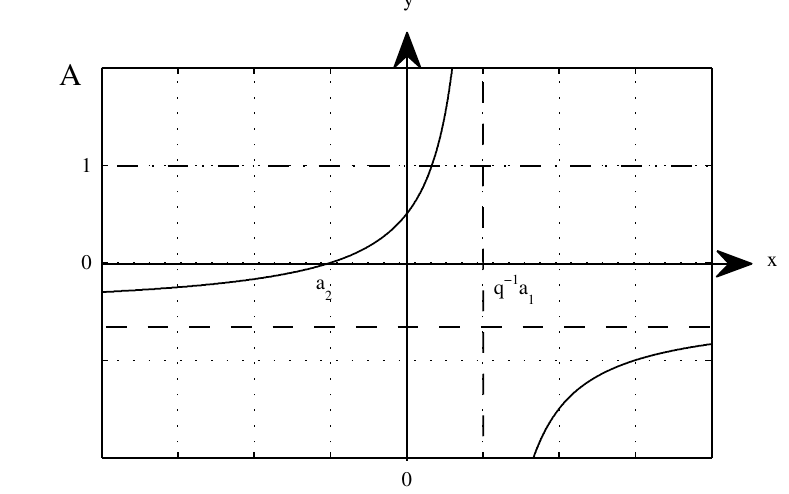}
\hfill
\includegraphics{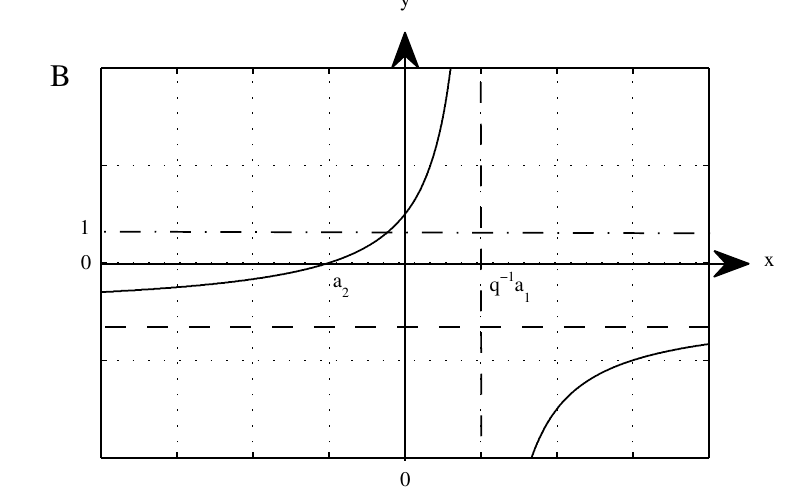}
\caption{The graph of $f(x, q)$ in \textbf{Case 6}
with $\Lambda_q<0$ and
$a_2<0<q^{-1}a_1$. In A, we have \textbf{Case 6(a)}. In B, we have \textbf{Case 6(b)}.}
\label{zerocase2IJ}
\end{figure}

In Figure \ref{zerocase2IJ}A,
the only interval is $(0, q^{-1}a_1)$ which corresponds to the interval described in 
Theorem \ref{thm*}-\textbf{2}. Notice that $\rho(qx, q)/\rho(x, q)=1$ at $0<x_0=-\tau'(0, q)/\tau(0, q)<q^{-1}a_1$
then $\rho$ is increasing on $(0, x_0)$ and decreasing on $(x_0, q^{-1}a_1)$.
Furthermore, since $0<y_0<1$, $\rho\to0$ 
as $x\to0^+$ according to Lemma \ref{zeropoint} 
and $\rho\to0$ as $x\to q^{-1}a_1^-$ since $\rho(qx, q)/\rho(x, q)\to+\infty$
as $x\to q^{-1}a_1^-$. Then, there exists a suitable
$\rho$ on $(0, a_1]$ supported at the points $q^ka_1$ for $k\in\mathbb{N}_0$
and, therefore, an OPS exists which is stated in the following theorem.
However, the case analysed in Figure \ref{zerocase2IJ}B does not give
any OPS.

\begin{theorem}\label{theoremzerocase2I}
Consider the case where $a_2<0<a_1$, $0<qy_0<1$ and $q^2\Lambda_q<0$.
Let $a=0$ and $b=a_1(q)$ be zeros of $\sigma_1(x,q)$.
Then there exists a sequence  of polynomials  $\{P_n\}$ for 
$n\in\mathbb{N}_0$ orthogonal on $(a,b]$ 
w.r.t. the weight function  in Theorem \ref{theoremzerocase1iiG},
supported on $\{q^kb\}_{k\in\mathbb{N}_0}$ (see \refe{qortho2} of Teorem \ref{thm*}-\textbf{2}).
\end{theorem}
The OPS in Theorem \ref{theoremzerocase2I} coincides with 
the case IVa4 in Chapter 10 of \cite[pages 278 and 312]{KLS}. 
An example of this family is again the little $q$-Jacobi polynomials 
with orthogonality on the set of points $\{...,q^2,q,1\}$. 
They satisfy the $q$-EHT with the coefficients in \refe{coefflittleqjacobi}
where $a_1=1$ and $a_2=b^{-1}q^{-1}$.
These polynomials have the same orthogonality property with the same $d_n^2$
as in \refe{dnlittleqjacobi} but the
constrains $0<a<q^{-1}$ and $b<0$ on the parameters are different
due to the conditions.
This extends the orthogonality relation
of the little $q$-Jacobi polynomials for $0<a<q^{-1}$ and $0<b<q^{-1}$
to a larger set of the parameters $0<a<q^{-1}$ and $b<0$.
Notice that combining this with the previous Case 5(a) one can obtain the
orthogonality relation of the little $q$-Jacobi polynomials 
for $0<aq<1, bq<1$.

\begin{figure}[!htp]
\centering
\includegraphics{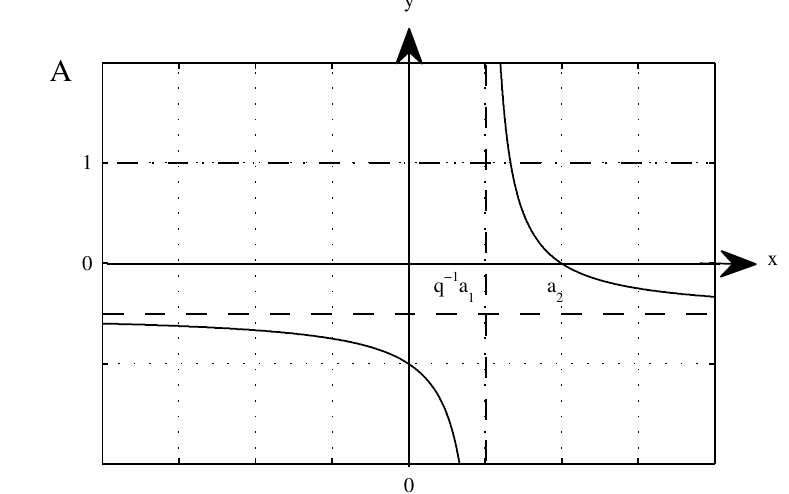}
\hfill
\includegraphics{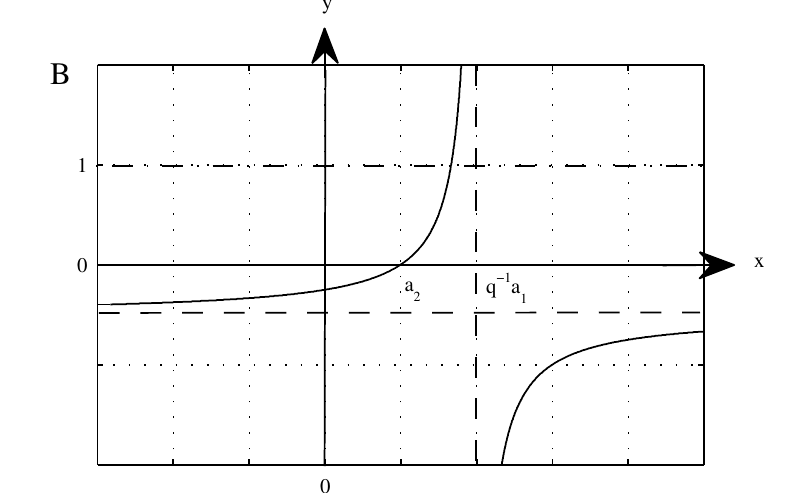}
\caption{The graph of $f(x, q)$ in \textbf{Case 6(c)}
with $\Lambda_q<0$. In A, the zeros are in order $0<q^{-1}a_1<a_2$. 
In B, the zeros are in order $0<a_2<q^{-1}a_1$.} 
\label{zerocase2KL}
\end{figure}

The case in Figure \ref{zerocase2KL}A does not yield any OPS.
On the other hand, in Figure \ref{zerocase2KL}B, the only possible
interval is $[a_2, q^{-1}a_1)$ which is 3{\it th} case in Theorem \ref{thm*}. 
Note that $\rho(qx, q)/\rho(x, q)=1$ at $a_2<x_0=-\tau'(0, q)/\tau(0, q)<q^{-1}a_1$
and that $\rho$ is increasing on $(a_2, x_0)$ and decreasing
on $(x_0, q^{-1}a_1)$. Furthermore,
$\rho(qa_2, q)=0$ since $\rho(qa_2, q)/\rho(a_2, q)=0$  and 
$\rho\to0$ as $x\to q^{-1}a_1^-$ since $\rho(qx, q)/\rho(x, q)\to\infty$
as $x\to q^{-1}a_1^-$.
Then there is an OPS on $\{q^{-k}a_2\}_{k=0}^N$ where $q^{-N-1}a_2=q^{-1}a_1$. 
Therefore, we have the following theorem.

\begin{theorem}\label{theoremzerocase2L}
Consider the case where $0<a_2<a_1$, $qy_0<0$ and $q^2\Lambda_q<0$.
Let $a=a_2(q)$ and $b=q^{-1}a_1(q)$ be zeros of $\sigma_2(x, q)$ and $\sigma_1(qx,q)$, respectively.
Then there exists a finite family of polynomials  $\{P_n\}$ orthogonal 
w.r.t. the weight function (see Eq. 4 in Table \ref{tab1})
\begin{equation*}
\rho(x, q)=\left|x\right|^{\alpha}\sqrt{x^{\log_qx-1}}(\frac{qa}{x}, b^{-1}x; q)_{\infty}, 
\quad q^{\alpha}=-\frac{q^{-3} \sigma_2''(0, q)}{ \sigma_1''(0, q)b}
\end{equation*}
supported on the set of points $\{q^{-k}a\}_{k=0}^N$ where $q^{-N-1}a=b$ 
(see \refe{qortho4-1} of Teorem \ref{thm*}-\textbf{3}).
\end{theorem}

An example of this family is the $q$-Kravchuk 
polynomials $K_n(x; p, N; q)$ satisfying the $q$-EHT with the coefficients
$$
\sigma_1(x, q)=q^{-2}x(x-a_1), \quad \sigma_2(x, q)=-px(x-a_2),
$$
$$
\tau(x, q)=\frac{1+pq}{(1-q)q}x-\frac{p+q^{-N-1}}{1-q}\quad {\rm and} \quad \lambda_n(q)=-q^{-n}[n]_q\frac{1+pq^{n}}{1-q}
$$
where $a_1=q^{-N}$ and $a_2=1$.
The conditions $q^2\Lambda_q<0$, $qy_0<0$ and $0<a_2<a_1$ lead to the condition $p>0$ on the parameter
of $K_n(x; p, N; q)$ with orthogonality
on $\{1, q^{-1},q^{-2},...,q^{-N}\}$ 
in the sense \refe{qortho4-1} where
\begin{equation}
d_n^2=(q^{-1}-1)
p^{-N}q^{-(^{N+1}_2)}(-q^{-N}p)^nq^{n^2}\frac{1+p}{1+pq^{2n}}(-pq; q)_N(q, q^{N+1}; q)_{\infty}\frac{(q, -pq^{N+1}; q)_n}{(-p, q^{-N}; q)_n}.
\end{equation}
In the literature, this relation is usually written as a finite sum \cite[Page 98]{ks}.
This case is not mentioned in \cite{KLS} for $0<q<1$. 
However the $q$-Kravchuk polynomials with this set of parameters
are described in \cite[page 98]{ks}.

In the following independent cases we fail to define an OPS.

 \medskip \noindent \textbf{Case 5(a)} with $0<\Lambda_q<1$ and $0<a_2<q^{-1}a_1$.

\medskip \noindent \textbf{Case 5(b)} with $\Lambda_q>1$ and $0<q^{-1}a_1<a_2$.

\medskip \noindent \textbf{Case 5(c)} with $\Lambda_q>1$ and $q^{-1}a_1<0<a_2$.

\medskip \noindent \textbf{Case 5(c)} with $0<\Lambda_q<1$ and $a_2<0<q^{-1}a_1$.


\section{The orthogonality of the $q$-polynomials: other cases}

This section includes the main analysis of the 
other families 
by taking into account the rational function on the r.h.s. of the $q$-Pearson equation
\refe{openqpearson} along the same lines with the $\emptyset$-Jacobi/Jacobi
and $0$-Jacobi/Jacobi cases handled in the previous section.

\subsection{$q$-Classical $\emptyset$-Jacobi/Laguerre Polynomials}\label{s4.2}
Let the coefficients $\sigma_2$ and $\sigma_1$ be quadratic
and linear polynomials in $x$, respectively, such that $\sigma_1(0, q)\sigma_2(0, q)\neq0$. If $\sigma_1$ is written in terms of 
its root, i.e., $\sigma_1(x, q)=\sigma_1'(0, q)[x-a_1(q)]$, 
$a_1(q)=-\frac{\sigma_1(0, q)}{\sigma_1'(0, q)}$ then from \refe{sigmaq12}
\[\sigma_2(x,q)=(q-1)\tau'(0, q)x^2+\left[q\sigma_1'(0, q)+(q-1)\tau(0, q)\right]x-q\sigma_1'(0, q)a_1(q)\]
where $\tau'(0,q)\neq0$ by hypothesis.
Then the $q$-Pearson equation \refe{openqpearson} takes the form
\begin{equation}\label{linearqpearson}
f(x, q):=\frac{\rho(qx, q)}{\rho(x, q)}=
\frac{q^{-1}\sigma_2(x, q)}{\sigma_1(qx, q)}=\frac{(1-q^{-1})\frac{\tau'(0, q)}{\sigma_1'(0, q)}[x-a_2(q)][x-b_2(q)]}{qx-a_1(q)}
\end{equation}
provided that the discriminant denoted by $ \Delta_q$, 
$$
 \Delta_q:=\left[1+\frac{(1-q^{-1})\tau(0, q)}{\sigma_1'(0, q)}\right]^2+4a_1(q)(1-q^{-1})\frac{\tau'(0, q)}{\sigma_1'(0, q)}
$$
of the quadratic polynomial in the nominator of $f$ in \refe{linearqpearson}
is non-zero. Note that here $x=a_2$ and $x=b_2$ are roots of $f$
which are constant multiplies of the roots of $\sigma_2$.
Moreover, $x=q^{-1}a_1$ is the vertical asymptote of $f$
and $y=1$ is its $y$-intercept since $\sigma_2(0, q)=q\sigma_1(0,q)$.
On the other hand, the locations of the zeros of $f$
are introduced by the following straightforward lemma.

\begin{lemma}
Let
$
\Lambda_q=\frac{\tau'(0, q)}{\sigma'_1(0, q)}\neq0.
$
Then, we have the following cases for the roots of the equation $f(x, q)=0$.

\medskip \noindent {\bf Case 1.} If $\Lambda_q$ and $a_1(q)$ 
have opposite signs, then there are two real distinct roots with opposite signs. 

\medskip \noindent {\bf Case 2.} If $\Lambda_q$ and $a_1(q)$ 
have same signs, then there exist three possibilities

\medskip {\bf (a)} if $\Delta_q>0$, $f$ has two real roots with same signs

\medskip {\bf (b)} if $\Delta_q=0$, $f$ has a double root

\medskip {\bf (c)} if $\Delta_q<0$, $f$ has a pair of complex conjugate roots.

\end{lemma}

\begin{figure}[!htp]
\centering
\includegraphics{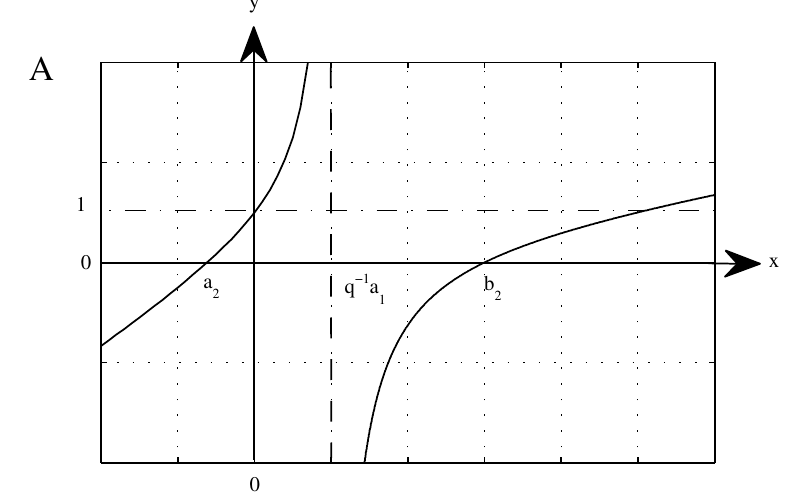}
\hfill
\includegraphics{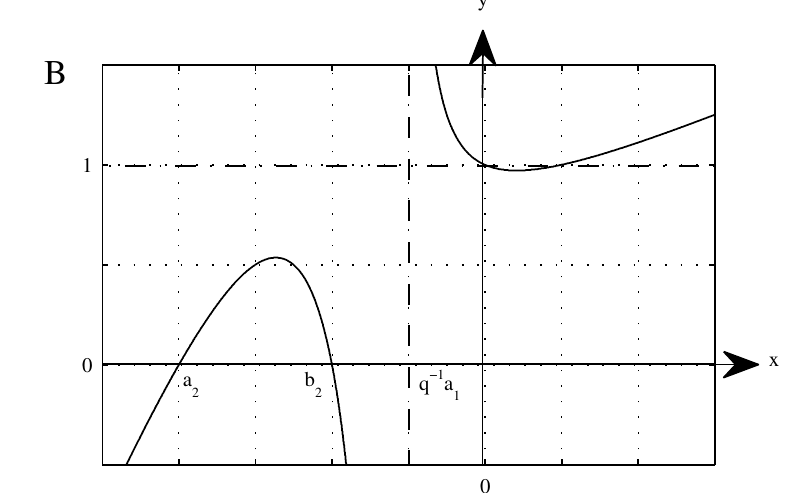}
\caption{The graph of $f(x, q)$. In A, we have \textbf{Case 1} with $\Lambda_q<0$ and $a_2<0<q^{-1}a_1<b_2$, 
and in B, we have \textbf{Case 2(a)} with $\Lambda_q<0$ and $a_2<b_2<q^{-1}a_1<0$.}
\label{linearsigma12}
\end{figure}

In Figure \ref{linearsigma12}A,
we first start with positivity condition of $q$-weight function
which allows us to exclude the intervals $(-\infty, a_2)$ and $(q^{-1}a_1, b_2)$.
Moreover, due to \textbf{PIII} 
$(a_2, q^{-1}a_1)$ can not be used.
On the other hand, the interval $(b_2, \infty)$
coincides with the 5{\it th} case of Theorem \ref{thm*}.
Notice that since  $\rho(qx, q)/\rho(x, q)=1$ at $x_0=-\tau(0, q)/\tau'(0, q)>b_2$,
$\rho$ is decreasing on $(x_0, \infty)$. Moreover, Since $\rho(qx, q)/\rho(x, q)$ has an infinite limit 
as $x\rightarrow+\infty$, we have $\rho\to0$ as $x\to\infty$. 
However, since it is infinite interval, we should check that 
$\sigma_1(x, q)\rho(x, q)x^k\to0$ as $x\to\infty$ by using
extended $q$-Pearson equation (\ref{openexpqpearson1}). 
The graph of the function $g$ defined in \refe{openexpqpearson1} looks like
the one for $f$. Then the analysis of the extended $q$-Pearson equation 
leads to $\sigma_1(x, q)\rho(x, q)x^k\to0$ as $x\to\infty$. Therefore, 
there exists a suitable $\rho$ supported on the set of points $\{q^{-k}b_2\}_{k\in\mathbb{N}_0}$. Thus, we have the following theorem.
\begin{theorem}\label{linearcase1A}
Let $a_2<0<a_1<b_2$ and $\Lambda_q<0$. 
Let $a=b_2(q)$ be the zero of $\sigma_2(x, q)$ and $b\to\infty$.
Then,  there exists a sequence  of polynomials  $(P_n)_n$ for $n\in\mathbb{N}_0$ 
orthogonal w.r.t. the weight function (see the 2nd expression of the
$\emptyset$-Jacobi/Laguerre case in Table \ref{tab2})
\begin{equation}\label{exp-linearcase1A}
\rho(x, q)=\left|x\right|^{\alpha}\sqrt{x^{\log_qx-1}}\frac{(qa_2/x, qa/x; q)_{\infty}}{(a_1/x; q)_{\infty}},
\quad q^{\alpha}=\frac{q^{-2}\hlf\sigma_2''(0, q)}{\sigma_1'(0, q)}  
\end{equation}
supported on $\{q^{-k}a\}_{k\in\mathbb{N}_0}$
(see \refe{qortho8} of Theorem \ref{thm*}-\textbf{5}).
\end{theorem}
The OPS in Theorem \ref{linearcase1A} coincides with
the case IIa2 in Chapter 11 of \cite[pages 337 and 358]{KLS}. 
An example of this family is the $q$-Meixner polynomials $M_n(x; b, c; q)$
satisfying the $q$-EHT with the coefficients
\[\sigma_1(x, q)=cq^{-2}(x-a_1), \quad \sigma_2(x, q)=(x-a_2)(x-b_2),\]
\begin{equation}\label{coeffqmeixner}
\tau(x, q)=-\frac{1}{1-q}x+\frac{cq^{-1}-bc+1}{1-q}\quad {\rm and} \quad \lambda_n(q)=\frac{[n]_q}{1-q}
\end{equation}
where $a_1=bq$, $a_2=-bc$ and $b_2=1$.
The conditions $\Lambda_q<0$ and $a_2<0<a_1<b_2$ give us the 
known restrictions $c>0$ and $0<b<q^{-1}$ on the parameters of 
$M_n(x; b, c; q)$ with orthogonality 
on $\{1, q^{-1},q^{-2},...\}$ in the sense \refe{qortho8} where
\begin{equation*}
d_n^2=(q^{-1}-1)c^{2n}q^{-n(2n+1)}(q, -c^{-1}q, bq; q)_n\frac{(q, -c; q)_{\infty}}{(bq; q)_{\infty}}.
\end{equation*}
In the literature, this relation can be found as an infinite sum \cite[page 360]{KLS}.

In Figure \ref{linearsigma12}B,
the only possible interval is $(q^{-1}a_1, \infty)$
which is the one identified in Theorem \ref{thm*}-\textbf{4}.  Notice that  
$\rho(qx, q)/\rho(x, q)=1$ at $x_0=-\tau(0, q)/\tau'(0, q)>q^{-1}a_1$,
then $\rho$ is increasing on $(q^{-1}a_1, x_0)$ and decreasing on $(x_0, \infty)$ which leads to
$\rho\to0$ as $x\to\infty$ since $\rho(qx, q)/\rho(x, q)\to\infty$.
But we still need to show
$\sigma_1(x, q)\rho(x, q)x^k\to0$ as $x\to\infty$ by using 
the extended $q$-Pearson equation (\ref{openexpqpearson1}).
By applying the same procedure to the extended $q$-Pearson
equation (\ref{openexpqpearson1}) 
whose graph looks like the one for $f$, 
we get $\sigma_1(x, q)\rho(x, q)x^k\to0$ as $x\to\infty$.
Consequently, we have a suitable $\rho$ on the interval 
$[a_1, \infty)$ supported on the set of points 
$\{a_1q^k\}_{k\in\mathbb{N}_0}\bigcup \{q^{\mp k}\}_{k\in\mathbb{N}_0}$. 

\begin{theorem}\label{linearcase2aD}
Let $a_2\leq b_2<a_1<0$, $\Lambda_q<0$. 
Let $a=a_1$ be a zero of $\sigma_1(x, q)$ and $b\to\infty$.
Then,  there exists a sequence  of polynomials  $(P_n)_n$ for $n\in\mathbb{N}_0$ 
orthogonal w.r.t. the weight function (see the 1st expression of the
$\emptyset$-Jacobi/Laguerre case in Table \ref{tab2})
\begin{equation*}
\rho(x, q)=\frac{(a^{-1}qx; q)_{\infty}}{(a_2^{-1}x, b_2^{-1}x; q)_{\infty}}
\end{equation*}
supported on the set of points $\{q^ka\}_{k\in\mathbb{N}_0}\bigcup
\{q^{\pm k}\alpha\}_{k\in\mathbb{N}_0}$ for arbitrary $\alpha>0$
in the sense \refe{qortho7} of Theorem \ref{thm*}-\textbf{4} with
\begin{equation}\label{dna1}
d_n^2=(1-q)q^{-n(2n-1)}\left(a_2b_2a_1^{-1}\right)^{2n}
(q, a_2^{-1}a_1, b_2^{-1}a_1; q)_n
\frac{(q, a_1, qa_1^{-1}, a_2^{-1}b_2^{-1}a_1, qa_2b_2a_1^{-1}; q)_{\infty}}
{(a_2^{-1}a_1, b_2^{-1}a_1, a_2^{-1}, b_2^{-1}, qa_2, qb_2; q)_{\infty}}.
\end{equation}
\end{theorem}
The OPS in Theorem \ref{linearcase2aD} coincides with 
the case VIa2 in Chapter 10 of \cite[pages 285 and 315]{KLS}.

\begin{figure}[!ht]
\centering
\includegraphics{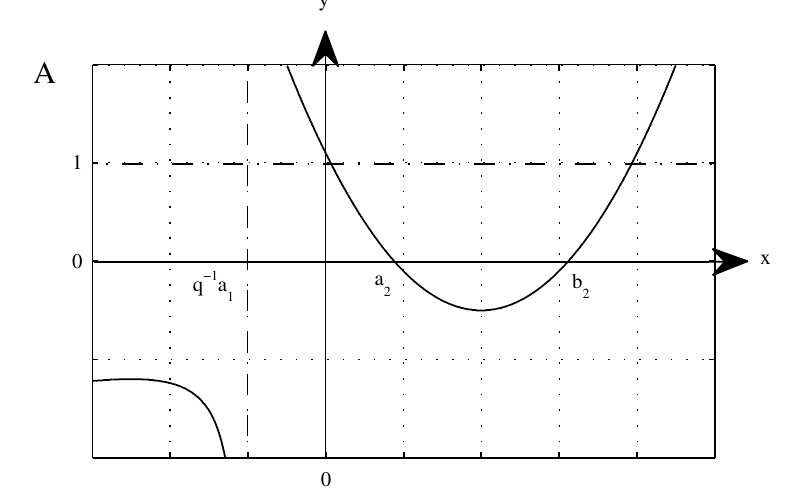}
\hfill
\includegraphics{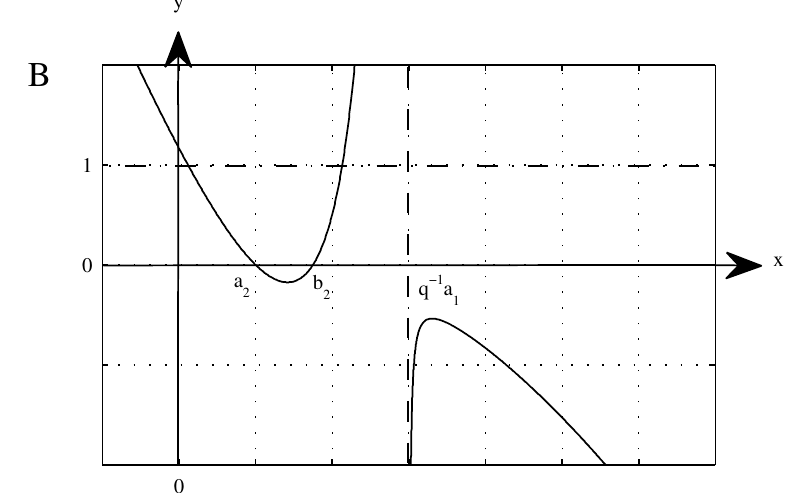}
\caption{The graph of $f(x, q)$ in \textbf{Case 2(a)}. In A, we have $\Lambda_q<0$ and $q^{-1}a_1<0<a_2<b_2$
and in B, $\Lambda_q>0$ and $0<a_2<b_2<q^{-1}a_1$.}\label{linearsigma34}
\end{figure}

In Figure \ref{linearsigma34}A,
the only possible interval is $(b_2, \infty)$. An analogous analysis
as the one that has been done for the case in Figure \ref{linearsigma12}A yields $\rho\to0$ as $x\to\infty$. 
Moreover, since from (\ref{openexpqpearson1}) $\sigma_1(x, q)\rho(x, q)x^k\to0$ as 
$x\to\infty$ for $k\in\mathbb{N}_0$, then there exists a $q$-weight function on $[b_2, \infty)$
supported at the points $q^{-k}b_2$ for $k\in\mathbb{N}_0$. Thus we have the following result.

\begin{theorem}\label{linearcase2aA}
Let $a_1<0<a_2\leq b_2$ and $\Lambda_q<0$. 
Let $a=b_2$ be the zero of $\sigma_2(x, q)$ and $b\to\infty$.
Then,  there exists a sequence  of polynomials  $(P_n)_n$ for 
$n\in\mathbb{N}_0$ orthogonal w.r.t. 
the weight function \refe{exp-linearcase1A} (see Theorem \ref{linearcase1A})
supported on $\{q^{-k}a\}_{k\in\mathbb{N}_0}$
(see \refe{qortho8} of Theorem \ref{thm*}-\textbf{5}).
\end{theorem}

A typical example of this family is again the $q$-Meixner polynomials 
orthogonal on $\{1, q^{-1},q^{-2},...\}$. 
They satisfy the $q$-EHT with the coefficients \refe{coeffqmeixner}
where $a_1=bq$, $a_2=-bc$ and $b_2=1$.
This set of $q$-Meixner polynomials corresponds to the case IIa2
in Chapter 11 of \cite[pages 337 and 358]{KLS} and their orthogonality 
relation is valid in a larger set of
parameters. In fact the conditions $a_1<0<a_2\leq b_2$ and $\Lambda_q<0$ yield
$c>0$, $b<0$ and $0<-bc\leq1$. This was not reported in \cite{KLS}.

In Figure \ref{linearsigma34}B, the only possible interval is 
$[b_2, q^{-1}a_1)$ which coincides with 3{\it th} case of Theorem \ref{thm*}.
In fact, $\rho(qx, q)/\rho(x, q)=1$ at $b_2<x_0=-\tau(0, q)/\tau'(0, q)<q^{-1}a_1$,
then $\rho$ is increasing on $(b_2, x_0)$ and decreasing on $(x_0, q^{-1}a_1)$. Moreover,
$\rho(qb_2, q)=0$ and $\rho(q^{-1}a_1, q)=0$ since $\rho(qb_2, q)/\rho(b_2, q)=0$
and $\rho(qx, q)/\rho(x, q)\to\infty$ as $x\to q^{-1}a_1^-$. Therefore, 
there is an OPS on $[b_2, q^{-1}a_1)$ w.r.t. a weight function supported on the set of points $\{q^{-k}b_2\}_{k=0}^N$ where $q^{-N-1}b_2=q^{-1}a_1$.

\begin{theorem}\label{linearcase2aB}
Let $0<a_2\leq b_2<a_1$ and $\Lambda_q>0$. 
Let $a=b_2$ be the zero of $\sigma_2(x, q)$ and $b=q^{-1}a_1$ of $\sigma_1(qx, q)$.
Then there exists a finite family of polynomials $(P_n)_n$
orthogonal w.r.t. the weight function  
(see the 3th expression  of $\emptyset$-Jacobi/Laguerre in Table \ref{tab2})
\begin{equation*}
\rho(x, q)=\left|x\right|^{\alpha}x^{\log_qx}(b^{-1}x, qa_2/x, a/x; q)_{\infty}, \quad
q^{\alpha}=-\frac{q^{-2}\hlf \sigma_2''(0, q)}{\sigma_1'(0, q)b}
\end{equation*}
supported on the set of points $\{q^{-k}a\}_{k=0}^N$
where $q^{-N-1}a=b$
(see \refe{qortho4-1} of Theorem \ref{thm*}-\textbf{3}).
\end{theorem}
The OPS in Theorem \ref{linearcase2aB} coincides with 
the case IIb1 in Chapter 11 of \cite[pages 337 and 361]{KLS}.  
An example of this family is the quantum $q$-Kravchuk polynomials 
$K_m^{qtm}(x; p, N; q)$ satisfying the $q$-EHT with the coefficients
$$\sigma_1(x, q)=-q^{-2}(x-a_1), \quad \sigma_2(x, q)=p(x-a_2)(x-b_2),$$
$$\tau(x, q)=-\frac{p}{1-q}x+\frac{p-q^{-1}+q^{-N-1}}{1-q} \quad 
{\rm and} \quad \lambda_n(q)=\frac{p}{1-q}[n]_q$$
where $a_1=q^{-N}$, $a_2=p^{-1}q^{-N-1}$ and $b_2=1$.
The conditions $\Lambda_q>0$ and $0<a_2\leq b_2<a_1$ give the 
constrain $p\geq q^{-N-1}$ on the parameter of
$K_m^{qtm}(x; p, N; q)$ with orthogonality 
on $\{1,q^{-1},q^{-2},...q^{-N}\}$ in the sense \refe{qortho4-1} where
\begin{equation*}
d_n^2=(q^{-1}-1)\frac{1}{(p^{-1}q^{-N}; q)_N}p^{-2n}q^{-n(2n+1)}(q, pq, q^{-N}; q)_n(q, p^{-1}q^{-N}, q^{N+1}; q)_{\infty}.
\end{equation*}
In the literature, this relation is usually written as a finite sum \cite[page 362]{KLS}.

\begin{figure}[!htp]
\centering
\includegraphics{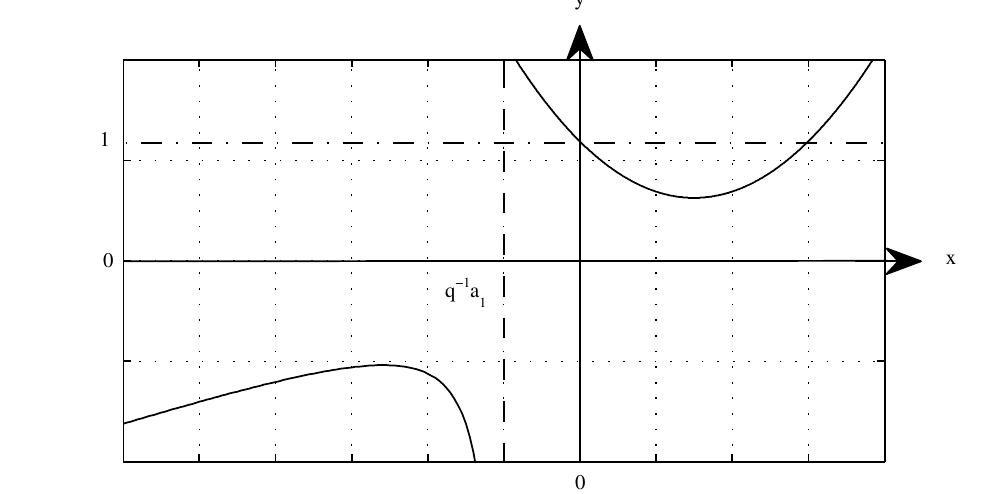}
\caption{The graph of $f(x, q)$ in \textbf{Case 2(c)} with $\Lambda_q<0$ and $a_1<0$, $a_2, b_2\in\mathbb{C}$.}\label{linearsigma8}
\end{figure}

In Figure \ref{linearsigma8}, $(q^{-1}a_1, \infty)$ is the only
interval where $f$ is positive.  Notice that the graphs of $f$ in the
interval $(q^{-1}a_1, \infty)$ in Figures \ref{linearsigma8} and \ref{linearsigma12}B 
have the same behaviour. Then, the analysis of Figure \ref{linearsigma12}B
is valid for this case and therefore there exists a suitable $\rho$ on $(q^{-1}a_1, \infty)$.
Thus, we have the following theorem.
\begin{theorem}\label{linearcase2c}
Let $a_1<0$, $a_2, b_2\in\mathbb{C}$ and $\Lambda_q<0$. 
Let $a=a_1$ be the zero of $\sigma_1(x, q)$ and $b\to\infty$.
Then, there exists a sequence  of polynomials  $(P_n)_n$ for 
$n\in\mathbb{N}_0$ orthogonal w.r.t. the weight function,  
given in Theorem \ref{linearcase2aD}, 
supported on the set of points $\{q^ka\}_{k\in\mathbb{N}_0}\bigcup
\{q^{\pm k}\alpha\}_{k\in\mathbb{N}_0}$ for arbitrary $\alpha>0$
in the sense \refe{qortho7} of Theorem \ref{thm*}-\textbf{4} with
$d_n^2$ defined by \refe{dna1}
\end{theorem}
The OPS in Theorem \ref{linearcase2c} coincides with 
the case VIa1 in Chapter 10 of \cite[pages 285 and 315]{KLS}.
The orthogonality relation of this OPS has the same form as in the
previous Case 2(a) defined in Theorem \ref{linearcase2aD} 
but now the zeros of $\sigma_2$  are complex.

For the two cases listed below the OPS fails to exist.

\medskip \noindent \textbf{Case 1.} $\Lambda_q>0$ and $q^{-1}a_1(q)<a_2(q)<0<b_2(q)$ and 
\textbf{Case 2(a).} $\Lambda_q>0$ and $0<a_2(q)<q^{-1}a_1(q)<b_2(q)$.

\subsection{$q$-Classical $\emptyset$-Jacobi/Hermite Polynomials}
Let the coefficients $\sigma_2$ and $\sigma_1$ be 
quadratic and constant polynomials in $x$, respectively, such that $\sigma_1(0, q)\sigma_2(0, q)\neq0$.
If $\sigma_1(x, q)=\sigma_1(0, q)\neq0$ then, from \refe{sigmaq12},
\[\sigma_2(x, q)=q\left[\sigma_1(x, q)+(1-q^{-1})x\tau(x, q)\right]=
(q-1)\tau'(0, q)x^2+(q-1)\tau(0, q)x+q\sigma_1(0, q)\]
where $\tau'(0, q)\neq0$ by hypothesis. Then the 
$q$-Pearson equation (\ref{openqpearson}) takes the form  
\begin{equation}\label{constantqpearson2}
f(x, q):=\frac{\rho(qx, q)}{\rho(x, q)}=
\frac{q^{-1}\sigma_2(x, q)}{\sigma_1(qx, q)}
=(1-q^{-1})\frac{\tau'(0, q)}{\sigma_1(0, q)}[x-a_2(q)][x-b_2(q)]
\end{equation}
provided that the discriminant denoted by $\Delta_q$, 
\[
 \Delta_q:=\left[(1-q^{-1})\frac{\tau(0, q)}{\sigma_1(0, q)}\right]^2-4(1-q^{-1})\frac{\tau'(0, q)}{\sigma_1(0, q)}
\]
of $f$ in \refe{constantqpearson2} is non-zero.
Notice that $y$-intercept  of $f$ 
is $y=1$ since $\sigma_2(0, q)=q\sigma_1(0, q)$. Moreover, 
$x=a_2$ and $x=b_2$ indicate its zeros which are constant
multiples of the roots of $\sigma_2$.
The following straightforward lemma allows us to determine the
locations of the zeros of $f$.

\begin{lemma}
Let $\Lambda_q=\frac{\tau'(0, q)}{\sigma_1(0, q)}\neq0$. 
Then we encounter the following cases for the roots of the equation $f(x, q)=0$. 

\medskip \noindent {\bf Case 1.} If $\Lambda_q>0$, $f$ has two real distinct roots with opposite signs. 

\medskip \noindent {\bf Case 2.} If $\Lambda_q<0$, there exist three possibilities

\medskip {\bf(a)} if $\Delta_q>0$, 
$f$ has two real roots with same signs

\medskip {\bf(b)} if $\Delta_q=0$, $f$ has a double root

\medskip {\bf(c)} if $\Delta_q<0$, $f$ has a pair of complex conjugate roots.
\end{lemma}

The next step is sketching roughly all graphs of $f$ by taking into account
all possible relative positions of the zeros of $f$ in question. As a result of 
analysis of the graphs of $f$, we determine a suitable $\rho>0$ satisfying the 
$q$-Pearson equation \refe{openqpearson} with BCs \refe{bc1}, \refe{bc3}.

\begin{figure}[!htb]
\centering
\includegraphics{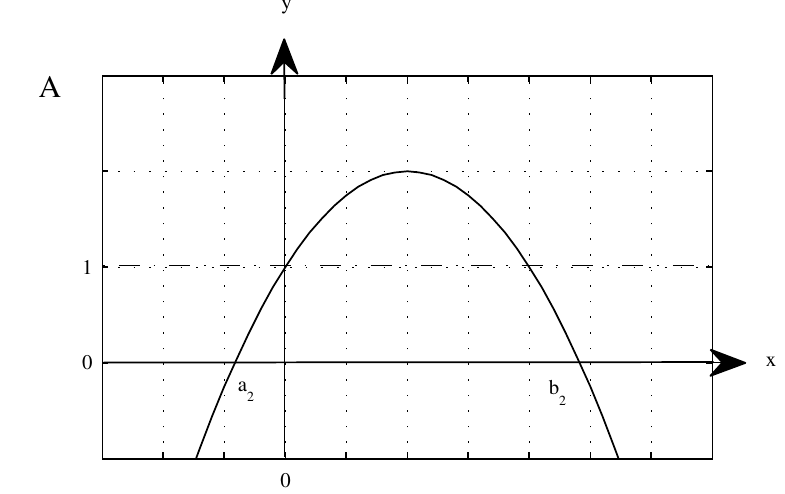}
\hfill
\includegraphics{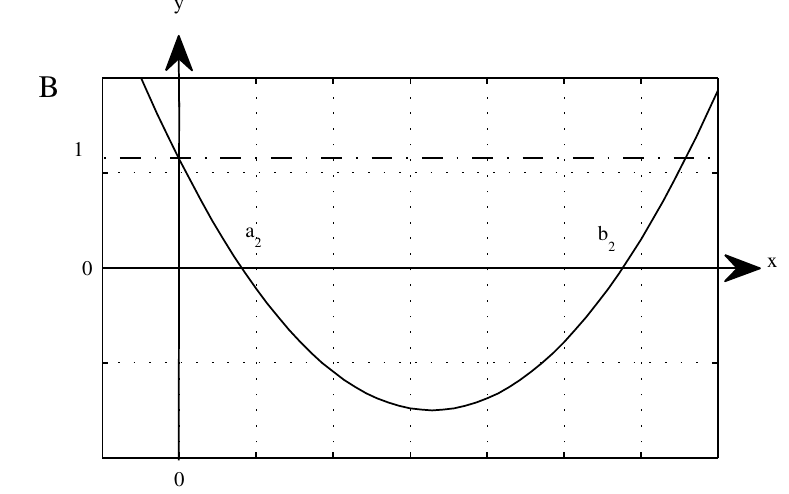}
\caption{The graph of $f(x, q)$. In A, we have \textbf{Case 1} with $\Lambda_q>0$ and $a_2<0<b_2$,
and in B, \textbf{Case 2(a)} with $\Lambda_q<0$ and $0<a_2<b_2$.}\label{constantsigma2}
\end{figure}
In Figure \ref{constantsigma2}A, let us consider the possible intervals in which we can have
a suitable weight function $\rho$ which are defined 
by the zeros of the polynomials $\sigma_1$ and $\sigma_2$. First of all,  notice that since
$\rho$ should be a positive weight function and $f$ is negative in the intervals
$(-\infty, a_2)$ and $(b_2, \infty)$, they are not suitable. On the other hand,
the interval $(a_2, b_2)$ is also eliminated in which $\rho=0$ due to \textbf{PII}.
As a result, an OPS fails to exist.

\begin{figure}[!htb]
\centering
\includegraphics{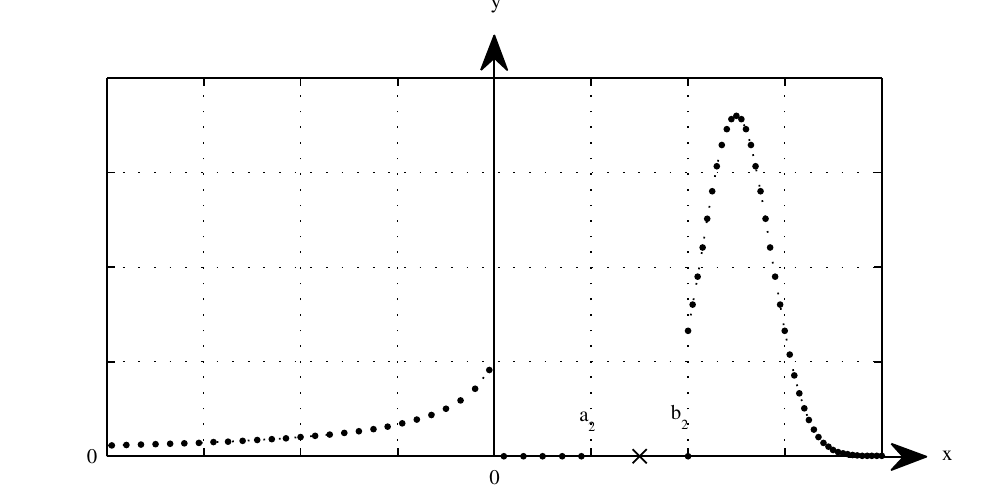}
\caption{The graph of $\rho(x, q)$ associated with the case in
Figure \ref{constantsigma2}B.}\label{constantsigma2rho}
\end{figure}

Let us analyse the case in Figure \ref{constantsigma2}B.
The positivity of $\rho$ implies that
the interval $(a_2, b_2)$ should be eliminated. On the other hand,
$(-\infty, a_2)$ is not suitable since $\rho=0$ in $(0, a_2)$ 
(this situation is similar to the one described in \textbf{PVI}.
The interval $(b_2, \infty)$ coincides with the 
5{\it th} case of Theorem \ref{thm*}. Notice that $\rho(qx, q)/\rho(x, q)=1$ at $x_0=-\tau(0, q)/\tau'(0, q)>b_2$,
then $\rho$ is decreasing on $(x_0, \infty)$. Since $f$ has  infinite limit 
as $x\rightarrow+\infty$, $\rho\to0$ as $x\to\infty$. As a result, the typical shape
of $\rho$ is constructed in Figure \ref{constantsigma2rho} assuming a positive initial
value of $\rho$ in each subinterval.

\begin{figure}[htb!]
\centering
\includegraphics{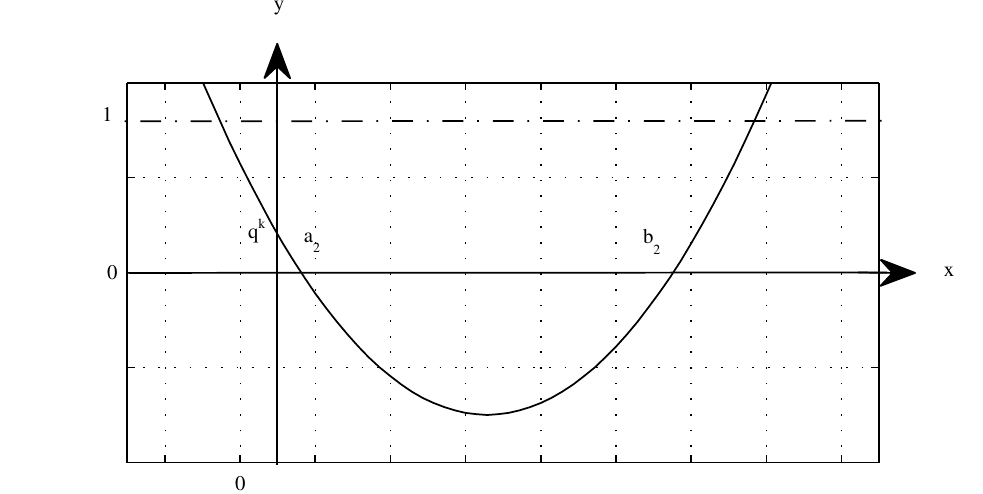}
\caption{The graph of $g(x, q)$ corresponding to Figure \ref{constantsigma2}B.}
\label{constantsigma2expandedqpearson}
\end{figure}

However, it is
not enough to assure that $\rho$ satisfies the BC at $+\infty$. 
In fact, even if $\rho\to0$ as $x\to\infty$ we should check that $\sigma_1(x, q)\rho(x, q)x^k\to0$ as $x\to\infty$
by using the extended $q$-Pearson equation (\ref{openexpqpearson1}), 
which is represented in Figure \ref{constantsigma2expandedqpearson}
for some $0<q<1$, where $k$ is large enough.

If we now provide a similar anaysis for $g$ in \refe{openexpqpearson1}, we see from Figure \ref{constantsigma2expandedqpearson} that, $g$ has the same property with
$f$. Therefore,  
$\sigma_1(x, q)\rho(x, q)x^k\rightarrow0$ as $x\to\infty$. 
That is, an OPS, to be stated in Theorem \ref{constantcase2abAB}, exists 
on the supporting set of points $\{b_2q^{-k}\}_{k\in\mathbb{N}_0}$.

\begin{theorem}\label{constantcase2abAB}
Let $0<a_2\leq b_2$ and $\Lambda_q<0$. 
Let $a=b_2(q)$ be a zero of $\sigma_2(x,q)$ and $b\rightarrow\infty$.
Then,  there exists a sequence  of polynomials  $(P_n)_n$ for $n\in\mathbb{N}_0$
orthogonal w.r.t. the weight function
 (see expression 2 for the $\emptyset$-Jacobi/Hermite case in Table \ref{tab2})
\begin{equation*}
\rho(x, q)=\left|x\right|^{\alpha}x^{\log_qx-1}(qa_2/x, qa/x; q)_{\infty},
\quad q^{\alpha}=\frac{q^{-1}\hlf\sigma_2''(0, q)}{\sigma_1(0, q)}
\end{equation*}
supported on the set of points $\{q^{-k}a\}_{k\in\mathbb{N}_0}$
(see \refe{qortho8} of Theorem \ref{thm*}-\textbf{5}).
\end{theorem}
The OPS in Theorem \ref{constantcase2abAB} coincides with 
the case Ia1 in Chapter 11 of \cite[pages 335 and 355-357]{KLS}.
In fact, a typical example of this family is the Al-Salam-Carlitz II 
polynomials $V_n^{(\alpha)}(x; q)$
satisfying the $q$-EHT with the coefficients
$$\sigma_1(x, q)=aq^{-1}, \quad \sigma_2(x, q)=(x-a_2)(x-b_2),$$
$$\tau(x, q)=\frac{1}{q-1}x-\frac{1+a}{q-1} \quad {\rm and} \quad \lambda_n(q)=\frac{1}{1-q}[n]_q$$
where $a_2=a$, $b_2=1$. 
The conditions $\Lambda_q<0$ and $0<a_2\leq b_2$ give the 
constrain $0<a\leq 1$ on the parameter of $V_n^{(\alpha)}(x; q)$
with orthogonality on $\{1,q^{-1},q^{-2},...\}$ in the sense \refe{qortho8} where
\begin{equation*}
d_n^2=(q^{-1}-1)a^nq^{-n^2}(q; q)_n(q; q)_{\infty}.
\end{equation*}
In the literature, this relation is usually written as an infinite sum \cite [page 357]{KLS}.

\begin{figure}[htp]
\centering
\includegraphics{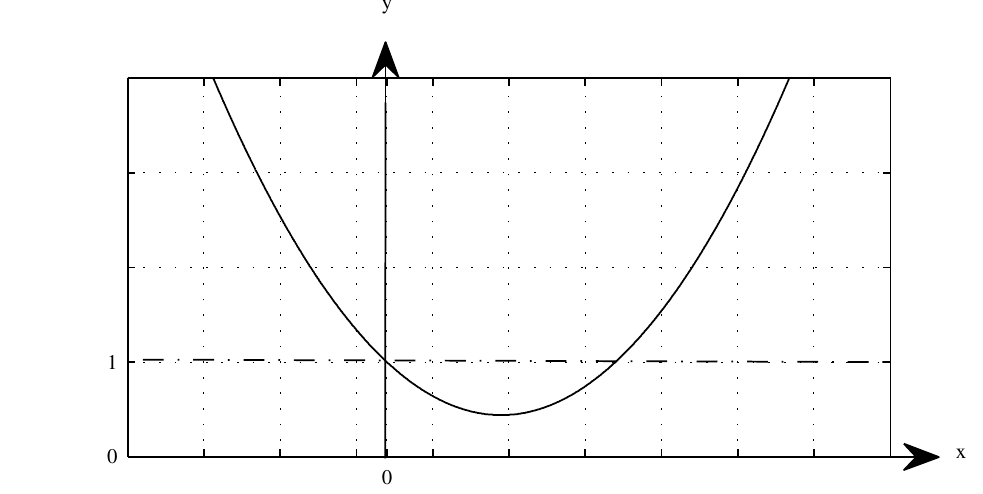}
\caption{The graph of $f(x, q)$ in \textbf{Case2(c)} with $\Lambda_q<0$ and $a_2(q), b_2(q) \in\mathbb{C}.$}
\label{constantsigma4}
\end{figure}

In Figure \ref{constantsigma4}, 
the only interval is $(-\infty, \infty)$ which corresponds to the 
7{\it th} case of Theorem \ref{thm*}.
Notice that $\rho(qx, q)/\rho(x, q)=1$ at $x_0=-\tau(0, q)/\tau'(0, q)$,
then it follows that $\rho$ is
increasing on $(-\infty, x_0)$ and decreasing
on $(x_0, \infty)$. Moreover, $\rho\rightarrow0$ as $x\rightarrow\mp\infty$
since $\rho(qx, q)/\rho(x, q)\rightarrow\infty$. 
Then, an OPS can exist on $\{\pm q^{\mp k}\}_{k\in\mathbb{N}_0}$.
But we should analyse the extended
$q$-Pearson equation (\ref{openexpqpearson1}) to check  
$\sigma_1(x, q)\rho(x, q)x^k\to0$ as $x\to\mp\infty$ which leads to 
similar figure as Figure \ref{constantsigma4}.
Then $\rho$ and $\sigma_1(x, q)\rho(x, q)x^k$ have same property that 
$q\sigma_1(x, q)\rho(x, q)x^k=\sigma_2(q^{-1}x, q)\rho(q^{-1}x, q)x^k\rightarrow 0$
as $x\rightarrow\mp\infty$ for $k\in\mathbb{N}_0$. 
Thus we can find a suitable $\rho$ 
supported on the set of points $\{\pm q^{\mp k}\}_{k\in\mathbb{N}_0}$. Therefore,
we have the following theorem.

\begin{theorem}\label{constantcase2c}
Let $\Lambda_q<0$ and $a_2, b_2 \in\mathbb{C}$.
Let $a\to-\infty$ and $b\to\infty$. 
Then, there exists a sequence  of polynomials  $(P_n)_n$ 
for $n\in\mathbb{N}_0$ orthogonal
w.r.t. the weight function  (see expression 1 for the $\emptyset$-Jacobi/Hermite in Table \ref{tab2})
\begin{equation*}
\rho(x, q)=\frac{1}{(a_2^{-1}x, b_2^{-1}x; q)_{\infty}}
\end{equation*}
supported on the set of points $\{\mp q^{\pm k}\alpha\}_{k\in\mathbb{N}_0}$
for arbitrary $\alpha>0$ (see \refe{qortho10} of Theorem \ref{thm*}-\textbf{7}).
\end{theorem}
The OPS in Theorem \ref{constantcase2c} corresponds to the case Ia1 in Chapter 11
and case Va2 in chapter 10 of \cite[pages 335, 355-356, 283 and 314-315]{KLS}.
An example of this family is the discrete $q$-Hermite II polynomials $\widetilde{h}_n(x; q)$
satisfying the $q$-EHT with the coefficients
$$\sigma_1(x, q)=q^{-1}, \quad \sigma_2(x, q)=(x-a_2)(x-b_2),$$
$$\tau(x, q)=\frac{1}{q-1}x \quad {\rm and} \quad \lambda_n(q)=\frac{1}{1-q}[n]_q$$
where $a_2=-i, b_2=i\in \mathbb{C}$.
Discrete $q$-Hermite II polynomials are orthogonal w.r.t. a measure
supported on the set of points $\{\pm q^{\mp k}\}_{k\in\mathbb{N}_0}$ with
\begin{equation*}
d_n^2=(1-q)q^{-n^2}(q; q)_n\frac{(q, -q, -1, -1, -q; q)_{\infty}}{(i, -i, -iq, iq, -i, i, iq, -iq; q)_{\infty}}
\end{equation*}
and the conditions $\Lambda_q<0$ and $a_2, b_2\in \mathbb{C}$ hold.

\subsection{$q$-Classical $\emptyset$-Laguerre/Jacobi Polynomials}
Let the coefficients $\sigma_2$ and $\sigma_1$ be linear and quadratic
polynomials in $x$, respectively, such that $\sigma_1(0,q)\sigma_2(0,q)\neq0$.
If $\sigma_1$ is written in terms of its roots, i.e., 
$\sigma_1(x, q)=\hlf\sigma_1''(0, q)[x-a_1(q)][x-b_1(q)]$,
then from (\ref{sigmaq12}) $\sigma_2(x, q)=\sigma_2'(0, q)x+\sigma_2(0, q)$ where
$$
\sigma_2'(0, q)=-q\Big[\hlf\sigma_1''(0, q)[a_1(q)+b_1(q)]-
(1-q^{-1})\tau(0, q)\Big]\neq0 \quad {\rm and} \quad \sigma_2(0, q)=q\hlf\sigma_1''(0, q)a_1(q)b_1(q)\neq0$$ 
provided that $\tau'(0, q)=-\frac{\hlf\sigma_1''(0, q)}{(1-q^{-1})}$.
Therefore, the $q$-Pearson equation (\ref{openqpearson}) takes the form
\begin{equation*}
f(x,q):=\frac{\rho(qx, q)}{\rho(x, q)}=\frac{-\Big[a_1(q)+b_1(q)-\frac{(1-q^{-1})\tau(0, q)}
{\hlf\sigma_1''(0, q)}\Big][x-a_2(q)]}{[qx-a_1(q)][qx-b_1(q)]}
\end{equation*}
where $\Big[a_1(q)+b_1(q)-\frac{(1-q^{-1})\tau(0, q)}
{\hlf\sigma_1''(0, q)}\Big]a_2(q)=a_1(q)b_1(q)$.
Let us point out that $f(x,q)$ intersects
the $y$-axis at the point $y=1$ since $\sigma_2(0, q)=q\sigma_1(0, q)$.
On the other hand, we consider the cases depending on
the signs of zeros of $\sigma_1$ and $\Lambda_q$ defined by
$$\Lambda_q:=\Big[a_1(q)+b_1(q)-\frac{(1-q^{-1})\tau(0, q)}
{\hlf\sigma_1''(0, q)}\Big].
$$

\medskip \noindent \textbf{Case 1.} $\Lambda_q<0$ with $a_1<0<b_1$, 
\medskip \noindent \textbf{Case 2.} $\Lambda_q>0$ with $0<a_1<b_1$,
\medskip \noindent \textbf{Case 3.} $\Lambda_q<0$ with $0<a_1<b_1$.

\begin{figure}[!htp]
\centering
\includegraphics{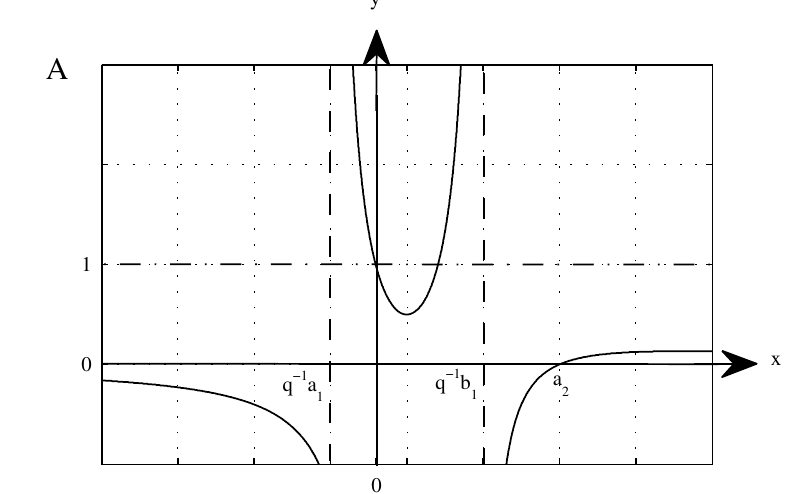}
\hfill
\includegraphics{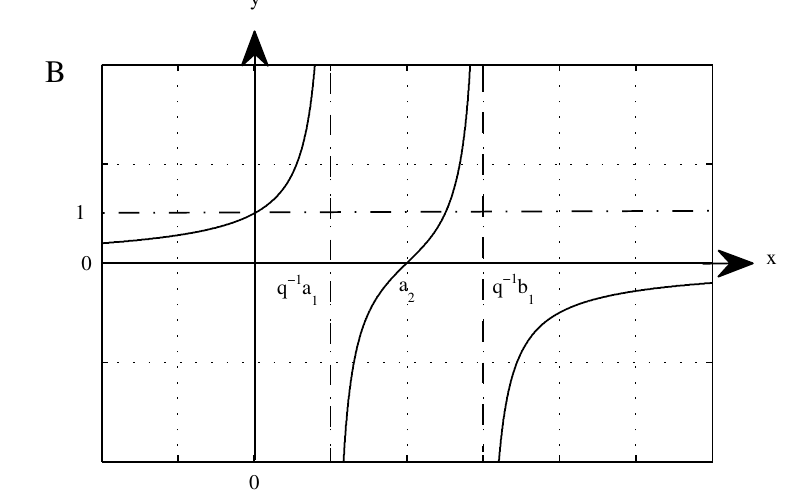}
\caption{The graph of $f(x, q)$. 
In A, we have \textbf{Case 1} with $\Lambda_q<0$ and $q^{-1}a_1<0<q^{-1}b_1<a_2$ and in B, 
\textbf{Case 2} with $\Lambda_q>0$ and $0<q^{-1}a_1<a_2<q^{-1}b_1$.}
\label{JLcase1AB}
\end{figure}

In Figure \ref{JLcase1AB}A, the only possible interval is $(q^{-1}a_1, q^{-1}b_1)$
which is the one described in Theorem \ref{thm*}-\textbf{1}. 
In fact, $\rho(qx, q)/\rho(x, q)=1$ at $q^{-1}a_1<x_0=-\tau(0, q)/\tau'(0, q)<q^{-1}b_1$,
Then, $\rho$ is increasing on $(q^{-1}a_1, x_0)$ and decreasing on
$(x_0, q^{-1}b_1)$. Moreover, $\rho\to0$ as $x\to q^{-1}a_1^+$ and $x\to q^{-1}b_1^-$
since $\rho(qx, q)/\rho(x, q)\to\infty$. Then there exists an OPS, 
to be stated in Theorem \ref{theoremJLcase1A},  
w.r.t. a $\rho$ supported at the points $x=q^ka_1$ and $x=q^kb_1$ for
$k\in\mathbb{N}_0$. 

\begin{theorem}\label{theoremJLcase1A}
Let $a_1<0<b_1<a_2$ and 
$\Lambda_q<0$. 
Let $a=a_1$ and $b=b_1$ be the zeros of $\sigma_1(x,q)$.
Then, there exists a sequence  of polynomials  $(P_n)_n$ for 
$n\in\mathbb{N}_0$ w.r.t.  weight function (see the 1st expression of the $\emptyset$-Laguerre/Jacobi 
case in Table \ref{tab2})
\begin{equation*}
\rho(x, q)=\frac{(qa^{-1}x, qb^{-1}x; q)_{\infty}}{(a_2^{-1}x; q)_{\infty}}
\end{equation*}
 supported on $\{q^ka\}_{k\in\mathbb{N}_0}
\bigcup\{q^kb\}_{k\in\mathbb{N}_0}$
(see \refe{qortho1} of Theorem \ref{thm*}-\textbf{1}).
\end{theorem}
The OPS in Theorem \ref{theoremJLcase1A} coincides with 
the case VIIa1 in Chapter 10 of \cite[pages 292 and 318]{KLS}.  
A typical example of this family is the big $q$-Laguerre polynomials $P_n(x; a, b; q)$
satisfying the $q$-EHT with the coefficients 
$$\sigma_1(x, q)=q^{-2}(x-a_1)(x-b_1), \quad \sigma_2(x, q)=-abq(x-a_2),$$
$$\tau(x, q)=-\frac{q^{-1}}{q-1}x+\frac{a+b-abq}{q-1}\quad {\rm and} \quad 
\lambda_n(q)=\frac{q^{-n}}{q-1}[n]_q$$
where $a_1=bq$, $b_1=aq$ and $a_2=1$. The conditions
$\Lambda_q<0$ and $a_1<0<b_1<a_2$ give the restrictions $b<0$ and $0<a<q^{-1}$
on the parameters of $P_n(x; a, b; q)$ with orthogonality on $\{bq,bq^2,bq^3,...\}
\bigcup\{...,aq^{3},aq^2,aq\}$ in the sense \refe{qortho1} where
\begin{eqnarray*}
d_n^2=(a-b)q(1-q)(-ab)^nq^{n(n+3)/2}(q; q)_n(aq, bq; q)_n
\frac{(q, a^{-1}bq, ab^{-1}q; q)_{\infty}}{(aq, bq; q)_{\infty}}.
\end{eqnarray*}

In Figure \ref{JLcase1AB}B,
the only possible interval is $[a_2, q^{-1}b_1)$ which coincides with the one described by 
Theorem \ref{thm*}-\textbf{3}.
Notice that $\rho(qx, q)/\rho(x, q)=1$ at $a_2<x_0=-\tau(0, q)/\tau'(0, q)<q^{-1}b_1$. 
Thus, $\rho$ is increasing on $(a_2, x_0)$ and decreasing on 
$(x_0, q^{-1}b_1)$. Moreover, $\rho(qa_2, q)=0$ 
and $\rho\to0$ as $x\to q^{-1}b_1^-$ since $\rho(qa_2, q)/\rho(a_2, q)=0$
and $\rho(qx, q)/\rho(x, q)\to\infty$ as $x\to q^{-1}b_1^-$.
Therefore, $[a_2, q^{-1}b_1)$ is suitable interval in which 
we have a positive $\rho$ supported on the set of points
$\{q^{-k}a_2\}_{k=0}^N$ where $q^{-N-1}a_2=q^{-1}b_1$. 

\begin{theorem}\label{theoremJLcase2E}
Let $0<a_1<a_2<b_1$ and $\Lambda_q>0$. 
Let $a=a_2$ be the zero of $\sigma_2(x,q)$ and $b=q^{-1}b_1$ of $\sigma_1(qx,q)$.
Then, there exists a finite family of polynomials  $(P_n)_n$  orthogonal 
w.r.t. the weight function (see the 2nd expression of the $\emptyset$-Laguerre/Jacobi case in Table \ref{tab2})
\begin{equation*}
\rho(x, q)=\left|x\right|^{\alpha}\frac{(a/x, qb^{-1}x; q)_{\infty}}{(a_1(q)/x; q)_{\infty}},
\quad q^{\alpha}=-\frac{q^{-2}\sigma_2'(0, q)}{\hlf\sigma_1''(0, q)b}
\end{equation*}
supported on $\{q^{-k}a_2\}_{k=0}^N$ where $q^{-N-1}a_2=q^{-1}b_1$. 
(see \refe{qortho4-1} of Theorem \ref{thm*}-\textbf{3}).
\end{theorem}

The OPS in Theorem \ref{theoremJLcase2E} coincides with 
the case IIIb3 in Chapter 11 of \cite[pages 343 and 363]{KLS}. 
An example of this family is the affine $q$-Kravchuk polynomials $K^{Aff}_n(x; p, N; q)$
satisfying the $q$-EHT with the coefficients
$$\sigma_1(x, q)=q^{-1}(x-a_1)(x-b_1), \quad \sigma_2(x, q)=-pq^{1-N}(x-a_2),$$
$$\tau(x, q)=\frac{1}{1-q}x-\frac{pq+q^{-N}-pq^{1-N}}{1-q}\quad {\rm and} \quad 
\lambda_n(q)=\frac{1}{q-1}[n]_{q^{-1}}$$
where $a_1=pq$, $b_1=q^{-N}$ and $a_2=1$.
The conditions $\Lambda_q>0$ and $0<a_1<a_2<b_1$ give the 
constrain $0<p<q^{-1}$ on the parameter of
$K^{Aff}_n(x; p, N; q)$ with orthogonality 
$\{1,q^{-1},q^{-2},...,q^{-N}\}$ in the sense \refe{qortho4-1} where
\begin{eqnarray*}
d_n^2=(-1)^np^{n-N}(q^{-1}-1)q^{-N(n+1)}q^{n(n+1)/2}(q, pq, q^{-N}; q)_n
\frac{(q, q^{N+1}; q)_{\infty}}{(pq; q)_{\infty}}.
\end{eqnarray*}
In the literature, this relation is usually written as a finite sum \cite[page 364]{KLS}.

The following four cases listed below fail to define an OPS.

\medskip \noindent \textbf{Case 1.} $\Lambda_q<0$ and $q^{-1}a_1<0<a_2<q^{-1}b_1$, 
\textbf{Case 2.} $\Lambda_q>0$ and $0<q^{-1}a_1<q^{-1}b_1<a_2$, 
\textbf{Case 2.} $\Lambda_q>0$ and $0<a_2<q^{-1}a_1<q^{-1}b_1$ and
\textbf{Case 3.} $\Lambda_q<0$ and $a_2<0<q^{-1}a_1<q^{-1}b_1$.

\subsection{$q$-Classical $\emptyset$-Hermite/Jacobi Polynomials}
Let the coefficients $\sigma_2$ and $\sigma_1$ be 
constant and quadratic polynomials in $x$, respectively, such that $\sigma_1(0,q)\sigma_2(0,q)\neq0$.
If $\sigma_1$ can be written in terms of its roots, i.e., 
$\sigma_1(x, q)=\hlf\sigma_1''(0, q)[x-a_1(q)][x-b_1(q)]$,
then, from (\ref{sigmaq12}) 
$$\sigma_2(x, q)=\sigma_2(0, q)=q\hlf\sigma_1''(0, q)a_1(q)b_1(q)$$
provided that $(1-q^{-1})\tau'(0, q)=-\hlf\sigma_1''(0, q)$ and
$(1-q^{-1})\tau(0, q)=\hlf\sigma_1''(0, q)[a_1(q)+b_1(q)]$.
Therefore, the $q$-Pearson equation (\ref{openqpearson}) becomes
\begin{equation*}
f(x,q):=\frac{\rho(qx, q)}{\rho(x, q)}=\frac{a_1(q)b_1(q)}
{[qx-a_1(q)][qx-b_1(q)]}.
\end{equation*}
Notice that the point $y=1$ is $y$-intercept of $f$.
In a similar fashion as before, we introduce the following two cases.

\medskip \noindent \textbf{Case1.} $a_1(q)<0<b_1(q)$. \textbf{Case 2.} $0<a_1(q)<b_1(q)$. 

\begin{figure}[!htp]
\centering
\includegraphics{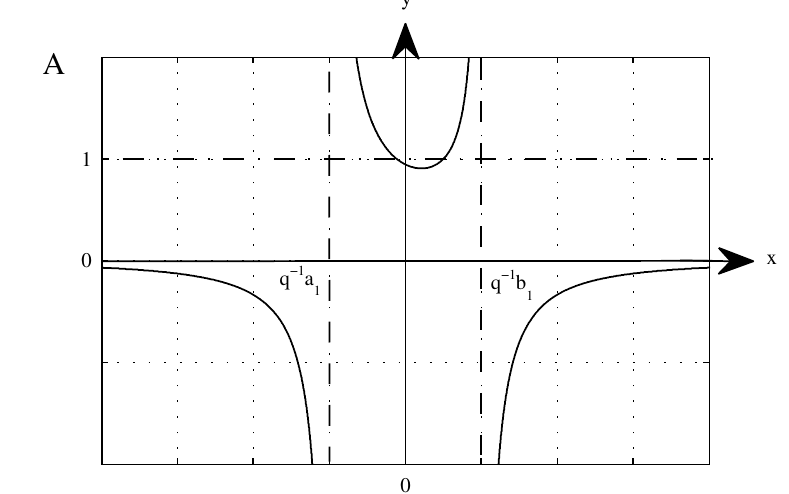}
\hfill
\includegraphics{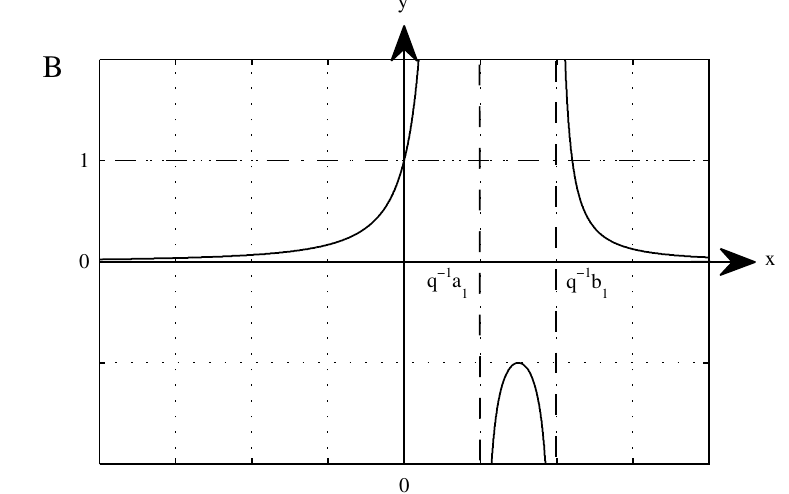}
\caption{The graph $f(x, q)$ in A, we have Case 1. and in B, Case 2.}
\label{JHcase1AB}
\end{figure}

In Figure \ref{JHcase1AB}A,
the only possible interval is $(q^{-1}a_1, q^{-1}b_1)$ which coincides with the 
Theorem \ref{thm*}-\textbf{1}. 
Notice that $\rho(qx, q)/\rho(x, q)=1$ at 
$q^{-1}a_1<x_0=-\tau(0, q)/\tau'(0, q)<q^{-1}b_1$.
Then, $\rho$ is increasing on $(q^{-1}a_1, x_0)$ and decreasing on $(x_0, q^{-1}b_1)$. 
Moreover, $\rho\to0$ as $x\to q^{-1}a_1^+$ and $x\to q^{-1}b_1^-$ 
since $\rho(qx, q)/\rho(x, q)\to\infty$.
It is obvious that 
BC holds at $x=a_1(q)$ and $x=b_1(q)$. Then 
there exists an OPS with positive $q$-weight function supported on
$\{q^ka_1\}_{k\in\mathbb{N}_0}\bigcup \{q^kb_1\}_{k\in\mathbb{N}_0}$, as it is stated 
in the following theorem.

\begin{theorem}\label{theoremJHcase1A}
Let $a_1<0<b_1$. 
Let $a=a_1$  and $b=b_1$ be the zeros of $\sigma_1(x,q)$.
Then,  there exists a sequence  of polynomials  $(P_n)_n$ 
for $n\in\mathbb{N}_0$ orthogonal 
w.r.t. the weight function (see the $\emptyset$-Hermite/Jacobi case in Table \ref{tab2})
\begin{equation*}
\rho(x, q)=(qa^{-1}x, qb^{-1}x; q)_{\infty}>0,
\, x\in(a, b)
\end{equation*}
supported on $\{q^ka_1\}_{k\in\mathbb{N}_0}\bigcup \{q^kb_1\}_{k\in\mathbb{N}_0}$  
(see \refe{qortho1} of Theorem \ref{thm*}-\textbf{1}).
\end{theorem}
The OPS in Theorem \ref{theoremJHcase1A} coincides with 
the case VIIa1 in Chapter 10 of \cite[pages 292 and 318-320]{KLS}
An example of this family is Al-Salam-Carlitz I polynomials $U^{(a)}_n(x; q)$
satisfying the $q$-EHT with the coefficients
$$\sigma_1(x, q)=q^{-1}(x-a_1)(x-b_1), \quad \sigma_2(x, q)=a,$$
$$\tau(x, q)=\frac{1}{1-q}x-\frac{1+a}{1-q} \quad {\rm and} \quad 
\lambda_n(q)=\frac{q^{1-n}}{q-1}[n]_{q}$$
where $a_1=a$ and $b_1=1$. 
The condition $a_1<0<b_1$ gives the restriction $a<0$ on the parameter of
$U^{(a)}_n(x; q)$ with orthogonality on $\{a,qa,q^2a,...\}\bigcup \{...,q^2,q,1\}$  
 in the sense \refe{qortho1} where
\begin{equation*}
d_n^2=(1-a)(-a)^{n}q^{(^n_2)}(1-q)(q; q)_n(q, aq, a^{-1}q; q)_{\infty}.
\end{equation*}

Another example of this family is the discrete $q$-Hermite I polynomials which are 
special case of Al-Salam-Carlitz I polynomials 
(see \cite[page 320]{KLS} for further details).
Finally, let us mention that the case represented in Figure \ref{JHcase1AB}B is 
inappropriate to define an OPS.

\subsection{$q$-Classical $0$-Jacobi/Laguerre Polynomials}
Let $\sigma_2$ and $\sigma_1$ be quadratic and linear polynomials in $x$, respectively,
such that $\sigma_2(0,q)=\sigma_1(0,q)=0$.
If $\sigma_1(x, q)=\sigma_1'(0, q)x$,
then from \refe{sigmaq12}, $\sigma_2(x, q)=\hlf\sigma_2''(0, q)x^2+\sigma_2'(0, q)x$ where
$$
\hlf\sigma_2''(0, q)=q(1-q^{-1})\tau'(0, q)\neq0 \quad {\rm and} \quad
\sigma_2'(0, q)=q[\sigma_1'(0, q)+(1-q^{-1})\tau(0, q)]\neq0
$$
provided that $(1-q^{-1})\tau(0, q)\neq-\sigma_1'(0, q)$.
For this case the $q$-Pearson equation reads
\begin{equation}\label{zerolinearqpearson1}
f(x,q):=\frac{\rho(qx, q)}{\rho(x, q)}=q^{-1}(1-q^{-1})\frac{\tau'(0, q)}{\sigma_1'(0, q)}[x-a_2(q)] 
\end{equation}
where $-(1-q^{-1})\frac{\tau'(0, q)}{\sigma_1'(0, q)}a_2(q)=1+\frac{(1-q^{-1})\tau(0, q)}{\sigma_1'(0, q)}$.
Let us point out that $f$ intersects the $y$-axis at the point 
$$y:=y_0=q^{-1}\left[1+\frac{(1-q^{-1})\tau(0, q)}{\sigma_1'(0, q)}\right].$$

Notice that for the zero cases one of the boundary of $(a,b)$ interval could be zero. 
This requires to find the behaviour of $\rho$ at the origin. 
\begin{lemma}\label{linearzeropoint1}
If $0<y_0<1$, then $\rho(z, q)\to0$ as $z\to0$. Otherwise it diverges to $\mp\infty$.
\end{lemma}
\noindent \textbf{Proof:}
From \refe{zerolinearqpearson1} it follows that 
\begin{eqnarray*}
\rho(q^kx, q)&=&q^{-k}\left[1+\frac{(1-q^{-1})\tau(0, q)}{\sigma_1'(0, q)}\right]^k(x/a_2(q); q)_k\rho(x, q)
\end{eqnarray*}
from where the result follows.\qed

\medskip

Again we identify the cases depending on $\sigma_2$, 
$\Lambda_q:=\frac{\tau'(0, q)}{\sigma_1'(0, q)}$ and 
$y_0$.

\medskip \noindent \textbf{Case 1.} $\Lambda_q>0$, $a_2>0$ and $y_0>1$, 
\textbf{Case 2.} $\Lambda_q<0$, $a_2<0$ and $0<y_0<1$, 
\textbf{Case 3.} $\Lambda_q<0$, $a_2>0$ and $y_0<0$. 

\begin{figure}[!htp]
\centering
\includegraphics{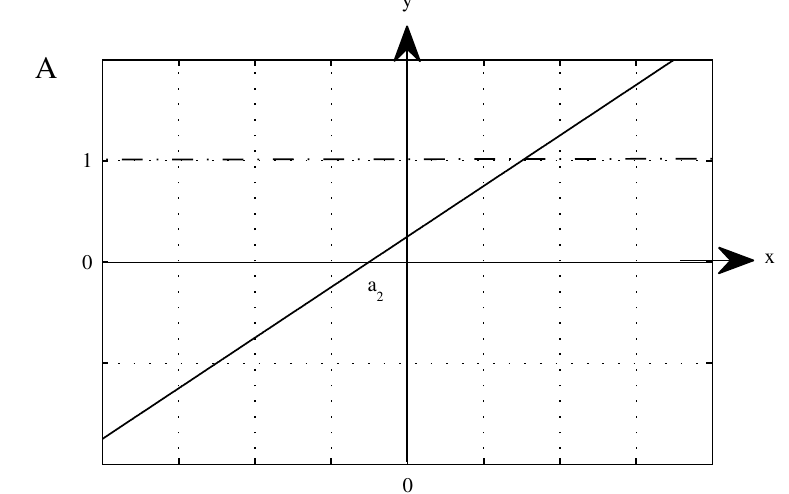}
\hfill
\includegraphics{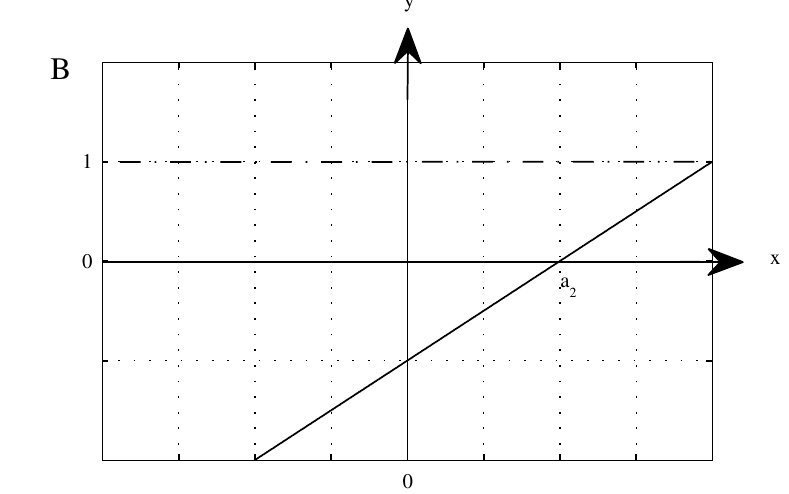}
\caption{The graph of $f(x, q)$ in A, we have
Case 2. and in B, Case 3.}
\label{zeroL-Jcase23AB}
\end{figure}

The \textbf{Case 1}, do not lead to any OPS. 
\textbf{Case 2-3} are introduced in Figure \ref{zeroL-Jcase23AB}.
In Figure \ref{zeroL-Jcase23AB}A,
the only possible interval is
$(0, \infty)$ which coincides with 6{\it th} case of Theorem \ref{thm*}. 
Notice that $\rho(qx, q)/\rho(x, q)=1$ at $x_0=-\tau(0, q)/\tau'(0, q)>0$. 
Then $\rho$ is increasing on $(0, x_0)$ and decreasing on $(x_0, \infty)$.
Furthermore, $\rho\to0$ as $x\to0^+$
by Lemma \ref{linearzeropoint1} since $0<y_0<1$ 
and $\rho\to0$ as $x\to\infty$ since $\rho(qx, q)/\rho(x, q)\to\infty$. 
Therefore, it could be possible to have 
a suitable $\rho$ on $(0, \infty)$. But we need to check $\sigma_1(x, q)\rho(x, q)x^k\to0$
as $x\to\infty$ for $k\in\mathbb{N}_0$ by using extended $q$-Pearson equation (\ref{openexpqpearson1}).
It is clear from (\ref{openexpqpearson1}) that graph of the function $g$ defined in \refe{openexpqpearson1} 
looks like the one represented in 
Figure \ref{zeroL-Jcase23AB}A with $y$-intercept, $0<q^{k+1}y_0<1$, $k\in\mathbb{N}_0$. Thus
$\sigma_1(x, q)\rho(x, q)x^k\to0$ as $x\to\infty$ for $k\in\mathbb{N}_0$ and therefore, 
there exists an OPS supported on $\{q^{\pm k}\}_{k\in\mathbb{N}_0}$ which is established in the next theorem.
\begin{theorem}\label{theoremzeroL-Jcase1B}
Let $\Lambda_q<0$, $a_2<0$ and  $0<qy_0<1$.
Let $a=0$  and $b\to\infty$. Then, there exists a sequence  
of polynomials $(P_n)_n$ for $n\in\mathbb{N}_0$ orthogonal 
w.r.t. the weight function (see the 1st expression of the $0$-Jacobi/Laguerre case in Table \ref{tab2})
\begin{equation*}
\rho(x, q)=\left|x\right|^{\alpha}\frac{1}{(a_2^{-1}x; q)_{\infty}},
 \quad q^{\alpha}=-\frac{q^{-2}\hlf\sigma_2''(0, q)a_2}{\sigma_1'(0, q)}
\end{equation*}
supported on $\{q^{\pm k}\}_{k\in\mathbb{N}_0}$
(see \refe{qortho9} of Theorem \ref{thm*}-\textbf{6}).
\end{theorem}
The OPS in Theorem \ref{theoremzeroL-Jcase1B} coincides with 
the case IIIa2 in Chapter 10 of \cite[pages 272 and 309]{KLS}.
An example of this family is the $q$-Laguerre polynomials $L^{(\alpha)}_n(x; q)$
satisfying the $q$-EHT with the coefficients
$$\sigma_1(x, q)=q^{-2}x, \quad \sigma_2(x, q)=q^{\alpha}x(x-a_2),$$
$$\tau(x, q)=-\frac{q^{\alpha}}{1-q}x+\frac{q^{-1}-q^{\alpha}}{1-q}
\quad {\rm and} \quad \lambda_n(q)=[n]_{q}\frac{q^{\alpha}}{1-q}$$
where $a_2=-1$. The conditions 
$\Lambda_q<0$, $a_2<0$ and $0<qy_0<1$ give the 
constrain $\alpha>-1$ on the parameter of $L^{(\alpha)}_n(x; q)$
with orthogonality on $\{q^{\pm k}\}_{k\in\mathbb{N}_0}$ in the sense \refe{qortho9} where
\begin{eqnarray*}
d_n^2= q^{-n}(1-q)\frac{(q^{\alpha+1}; q)_n}{(q; q)_n}
\frac{(q, -q^{\alpha+1}, -q^{-\alpha}; q)_{\infty}}{(q^{\alpha+1}, -q, -q; q)_{\infty}}.
\end{eqnarray*}
 In Figure \ref{zeroL-Jcase23AB}B,
the positivity of $\rho$ enables us to skip the intervals $(-\infty, 0)$ and $(0, a_2)$.
So the only interval is $(a_2, \infty)$ which is the one
described in Theorem \ref{thm*}-\textbf{5}. 
Notice that $\rho(qx, q)/\rho(x, q)=1$ at $x_0=-\tau(0, q)/\tau'(0, q)>a_2$.
Therefore, $\rho$ is increasing on $(a_2, x_0)$ and decreasing on $(x_0, \infty)$.
Moreover, $\rho(qa_2, q)=0$ since $\rho(qa_2, q)/\rho(a_2, q)=0$ and
$\rho\to0$ as $x\to\infty$ since $\rho(qx, q)/\rho(x, q)\to\infty$.
Furthermore, since the graph of the function $g$ defined in \refe{openexpqpearson1} looks like 
the one represented in Figure \ref{zeroL-Jcase23AB}B
one can conclude that $\sigma_1(x, q)\rho(x, q)x^k\to0$ as $x\to\infty$ for $k\in\mathbb{N}_0$
and therefore we have the following theorem.

\begin{theorem}\label{theoremzeroL-Jcase3C}
Let $\Lambda_q<0$, $a_2>0$ and  $qy_0<0$.
Let $a=a_2$ be the zero of $\sigma_2(x, q)$  and $b\to\infty$.
Then, there exists a sequence of polynomials $(P_n)_n$ 
for $n\in\mathbb{N}_0$ orthogonal w.r.t. the weight function
(see the 2nd expression of the $0$-Jacobi/Laguerre case in Table \ref{tab2})
\begin{equation*}
\rho(x, q)=\left|x\right|^{\alpha}\sqrt{x^{\log_qx-1}}(qa/x; q)_{\infty},
\quad q^{\alpha}=\frac{q^{-2}\hlf\sigma_2''(0, q)}{\sigma_1'(0, q)}
\end{equation*}
supported on the set of points $\{q^{-k}a_2\}_{k\in\mathbb{N}_0}$
(see \refe{qortho8} of Theorem \ref{thm*}-\textbf{5}).
\end{theorem}
The OPS in Theorem \ref{theoremzeroL-Jcase3C} coincides with 
the case IIa2 in Chapter 11 of \cite[pages 337 and 358]{KLS}.
An example of this family is the $q$-Charlier polynomials $C_n(x; a; q)$ 
satisfying the $q$-EHT with the coefficients
$$\sigma_1(x, q)=aq^{-2}x, \quad \sigma_2(x, q)=x(x-a_2),$$
$$\tau(x, q)=-\frac{1}{1-q}x+\frac{a+q}{(1-q)q}\quad {\rm and} \quad \lambda_n(q)=[n]_{q}\frac{1}{1-q}$$
where $a_2=1$. 
The conditions $\Lambda_q<0$, $a_2>0$ and $qy_0<0$ give the 
restriction $a>0$ on the parameter of $C_n(x; a; q)$
with orthogonality 
on $\{1,q^{-1},q^{-2},...\}$  in the sense \refe{qortho8} where
\begin{eqnarray*}
d_n^2=a^{2n}q^{-n(2n+1)}(-a^{-1}q, q; q)_n(-a, q; q)_{\infty}.
\end{eqnarray*}
In the literature, this relation is usually written as an infinite sum \cite[page 360]{KLS}.

\subsection{$q$-Classical $0$-Bessel/Jacobi Polynomials}
Let $\sigma_2$ and $\sigma_1$ be quadratic polynomials in $x$, respectively, such that
$\sigma_2'(0,q)=0$ and $\sigma_2(0,q)=\sigma_1(0,q)=0$.
If $\sigma_1(x, q)=\hlf\sigma''_1(0, q)x[x-a_1(q)]$,
$\frac{\tau'(0, q)}{\hlf\sigma_1''(0, q)}\neq
-\frac{1}{(1-q^{-1})}$ and $\frac{\tau(0, q)}{\hlf\sigma_1''(0, q)}=
\frac{a_1(q)}{(1-q^{-1})}$, then from (\ref{sigmaq12}) we have 
$\sigma_2(x, q)=\hlf\sigma_2''(0, q)x^2
=q\left[\hlf\sigma_1''(0, q)+(1-q^{-1})\tau'(0, q)\right]x^2$.
As a result, the $q$-Pearson equation (\ref{openqpearson}) becomes
\begin{eqnarray*}
f(x,q):=\frac{\rho(qx, q)}{\rho(x, q)}&=&
\frac{\left[1+\frac{(1-q^{-1})\tau'(0, q)}{\hlf\sigma_1''(0, q)}\right]x}{q[qx-a_1(q]}.
\end{eqnarray*}
Let us point out that
$f(x, q)$ passes through the origin and the line
$y= \Lambda_q:=q^{-2}\left[1+\frac{(1-q^{-1})\tau'(0, q)}{\hlf\sigma_1''(0, q)}\right]\neq0$
is its horizontal asymptote.
Hence, we have the following two cases: 

\medskip \noindent \textbf{Case 1.} $\Lambda_q<0$ and $a_1>0$ and 
\textbf{Case 2.} $\Lambda_q>0$ and $a_1>0$. 

\medskip

\begin{figure}[!htp]
\centering
\includegraphics{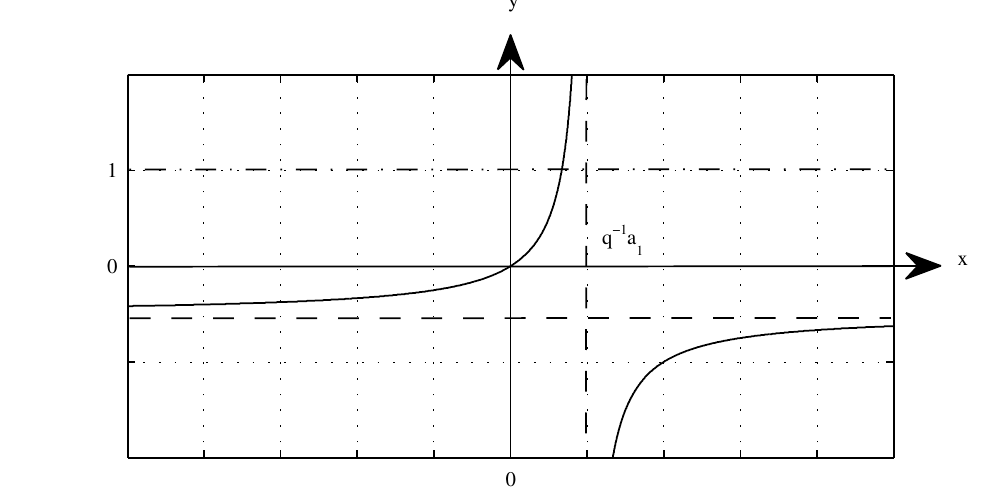}
\caption{The graph of $f(x, q)$ in Case 1.}
\label{zeroJ-Bcase2}
\end{figure}

The \textbf{Case 2} with $\Lambda_q>1$ and $0<\Lambda_q<1$ do not lead to any OPS.
The \textbf{Case 1} is represented in Figure \ref{zeroJ-Bcase2} from where it 
follows that the only possible interval is $(0, q^{-1}a_1)$ which is the one defined 
in Theorem \ref{thm*}-\textbf{2}.
Notice also that $\rho(qx, q)/(x, q)=1$ at $0<x_0=-\tau(0, q)/\tau'(0, q)<q^{-1}a_1$. Then,
$\rho$ is increasing on $(0, x_0)$ and decreasing on $(x_0, q^{-1}a_1)$.
Moreover, $\rho\to0$ as $x\to0^+$ and $x\to q^{-1}a_1^-$ since
$\rho(qx, q)/(x, q)\to0$ and $\rho(qx, q)/(x, q)\to\infty$, respectively.
Then, there exists an OPS with a suitable $\rho$ defined on $(0, a_1]$ supported at the points
$a_1q^k$ for $k\in\mathbb{N}_0$ and the following theorem holds.

\begin{theorem}\label{theoremzeroJ-Bcase2C}
Let $q^2\Lambda_q<0$ and $a_1>0$. 
Let $a=0$ and $b=a_1$ be the zeros of $\sigma_1(x,q)$.
Then, there exists a sequence  of polynomials  
$(P_n)_n$ for $n\in\mathbb{N}_0$ orthogonal 
w.r.t. the weight function (see the $0$-Bessel/Jacobi case in Table \ref{tab2})
\begin{equation*}
\rho(x, q)=\left|x\right|^{\alpha}\sqrt{x^{\log_qx-1}}(b^{-1}qx; q)_{\infty},
\quad q^{\alpha}=-\frac{q^{-2}\hlf \sigma_2''(0, q)}{\hlf \sigma_1''(0, q)b}
\end{equation*}
supported on the set of points $\{q^ka_1\}_{k\in\mathbb{N}_0}$
(see \refe{qortho2} of Theorem \ref{thm*}-\textbf{2}).
\end{theorem}
The OPS in Theorem \ref{theoremzeroJ-Bcase2C} coincides with 
the case IVa5 in Chapter 10 of \cite[pages 278 and 313]{KLS}. 
An example of this family is the Alternative $q$-Charlier ($q$-Bessel) polynomials 
$K_n(x; a; q)$ satisfying the $q$-EHT with the coefficients
$$\sigma_1(x, q)=-q^{-2}x(x-a_1), \quad \sigma_2(x, q)=ax^2,$$
$$\tau(x, q)=-\frac{1+aq}{(1-q)q}x+\frac{1}{(1-q)q}\quad {\rm and} \quad \lambda_n(q)=q^{-n}[n]_{q}\frac{1+aq^n}{1-q}$$
where $a_1=1$. 
The conditions $q^2\Lambda_q<0$ and $a_1>0$ give the 
constrain $a>0$ on the parameter of $K_n(x; a; q)$
with orthogonality on $\{...,q^2,q, 1\}$ in the sense \refe{qortho2} where
\begin{equation*}
d_n^2=a^nq^{n(3n-1)/2}(-aq, q; q)_{\infty}\frac{(q, -a; q)_n}{(-a, -aq; q)_{2n}}.
\end{equation*}
In the literature, this relation can be found as an infinite sum \cite[page 314]{KLS}.

\subsection{$q$-Classical $0$-Bessel/Laguerre Polynomials}
Let $\sigma_2$ and $\sigma_1$ be quadratic and linear polynomials in $x$, respectively,
such that $\sigma_2'(0,q)=0$ and $\sigma_2(0,q)=\sigma_1(0,q)=0$.
If $\sigma_1(x, q)=\sigma_1'(0, q)x$,
then, from \refe{sigmaq12} $\sigma_2(x, q)=\hlf\sigma_2''(0, q)x^2=q(1-q^{-1})\tau'(0, q)x^2$
provided that $(1-q^{-1})\tau(0, q)=-\sigma_1'(0, q)$.
So the $q$-Pearson equation is now
\begin{equation*}
f(x,q):=\frac{\rho(qx, q)}{\rho(x, q)}=
q^{-1}(1-q^{-1})\frac{\tau'(0, q)}{\sigma_1'(0, q)}x.
\end{equation*}
Clearly, $f$ passes through the origin.
According to the sign of $\Lambda_q:=\frac{\tau'(0, q)}{\sigma_1'(0, q)}$
we have only one possible case.
\begin{figure}[!hb]
\centering
\includegraphics{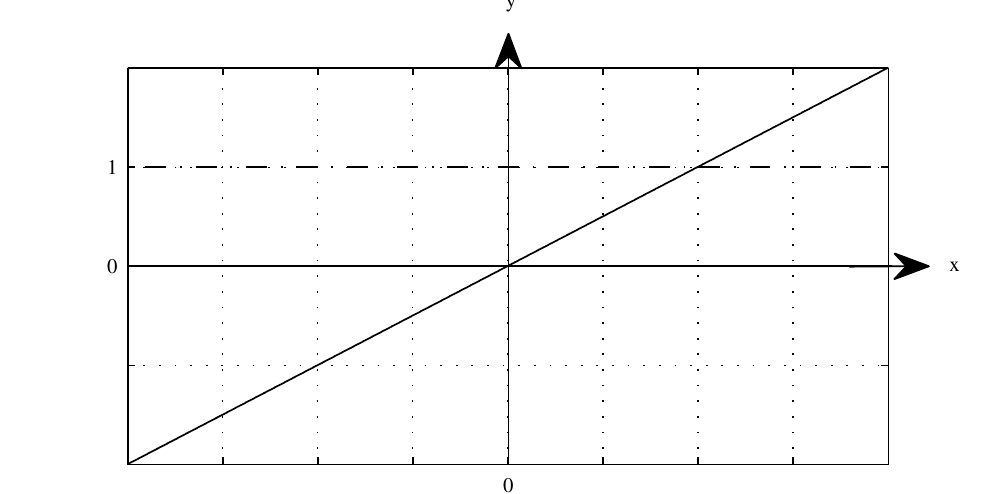}
\caption{The graph of $f(x, q)$ with $\Lambda_q<0$, $a_2=0$.}
\label{zeroL-B}
\end{figure}

From Figure \ref{zeroL-B} it follows that $(0, \infty)$ is the only possible interval and it
coincides with the one described in Theorem \ref{thm*}-\textbf{6}. 
Notice that $\rho(qx, q)/\rho(x, q)=1$ at $x_0=-\tau(0, q)/\tau'(0, q)>0$.
Then, $\rho$ is increasing on $(0, x_0)$ and decreasing on $(x_0, \infty)$. Moreover,
by use of the extended $q$-Pearson equation (\ref{openexpqpearson1}) it is straightforward to see
that $\sigma_1(x, q)\rho(x, q)x^k\to0$ as $x\to+\infty$. Thus, the following theorem holds.

\begin{theorem}\label{theoremzeroL-Bcase1A}
Let $\Lambda_q<0$, $a_2=0$ and  $qy_0=0$.
Let $a=0$  and $b\to\infty$.
Then, there exists a sequence of polynomials
$(P_n)_n$ for $n\in\mathbb{N}_0$ orthogonal w.r.t. the weight function (see the $0$-Bessel/Laguerre case in Table \ref{tab2})  
\begin{equation*}
\rho(x, q)=\left|x\right|^{\alpha}\sqrt{x^{\log_qx-1}},
\quad q^{\alpha}=\frac{q^{-2}\hlf\sigma_2''(0, q)}{\sigma_1'(0, q)}
\end{equation*}
supported on $\{q^{\mp k}\}_{k\in\mathbb{N}_0}$
(see \refe{qortho9} of Theorem \ref{thm*}-\textbf{6}).
\end{theorem}
The OPS in Theorem \ref{theoremzeroL-Bcase1A} coincides with 
the case IIIa2 in Chapter 10 of \cite[pages 272 and 309]{KLS}.
An example of this family is Stieltjes-Wigert polynomials $S_n(x; q)$ 
satisfying the $q$-EHT with the coefficients
$$\sigma_1(x, q)=q^{-2}x, \quad \sigma_2(x, q)=x^2,$$
$$\tau(x, q)=-\frac{1}{1-q}x+\frac{1}{(1-q)q}\quad {\rm and} \quad \lambda_n(q)=[n]_{q}\frac{1}{1-q}.$$
The conditions $\Lambda_q<0$, $a_2=0$ and $qy_0=0$ are satisfied for
$S_n(x; q)$ and they are orthogonal w.r.t. a measure 
supported on $\{q^{\mp k}\}_{k\in\mathbb{N}_0}$ in the sense \refe{qortho9} with
\begin{eqnarray*}
d_n^2=q^{-n}(1-q)\frac{(-tq, -1/t, q; q)_{\infty}}{(q^2; q)_n}.
\end{eqnarray*}

\subsection{$q$-Classical $0$-Laguerre/Jacobi Polynomials}
Let $\sigma_2$ and $\sigma_1$ be linear and quadratic polynomials in $x$, respectively,
such that $\sigma_2(0,q)=\sigma_1(0,q)=0$.
If $\sigma_1(x, q)=\hlf\sigma''_1(0, q)x[x-a_1(q)]$
and $\frac{\tau'(0, q)}{\hlf\sigma_1''(0, q)}=
-\frac{1}{(1-q^{-1})}$, then from (\ref{sigmaq12}) we get
$\sigma_2(x, q)=\sigma_2'(0, q)x=q\left[(1-q^{-1})\tau(0, q)-\hlf\sigma_1''(0, q)a_1(q)\right]x$.
Therefore, the $q$-Pearson equation has the form
\begin{equation*}
f(x,q):=\frac{\rho(qx, q)}{\rho(x, q)}=\frac{(1-q^{-1})\frac{\tau(0, q)}{\hlf\sigma_1''(0, q)}-a_1(q)}{q[qx-a_1(q]}.
\end{equation*}
Notice that $y=0$ is the horizontal asymptote of $f(x, q)$, and
its $y$-intercept is
$$y:=y_0=q^{-1}\left[1-\frac{(1-q^{-1})}{a_1(q)}\frac{\tau(0, q)}{\hlf\sigma_1''(0, q)}\right].$$

We have the following two cases: \textbf{Case 1.} $y_0>0$ and $a_1>0$, \textbf{Case 2.} $y_0<0$ and $a_1>0$. 

\medskip

\begin{figure}[!htp]
\centering
\includegraphics{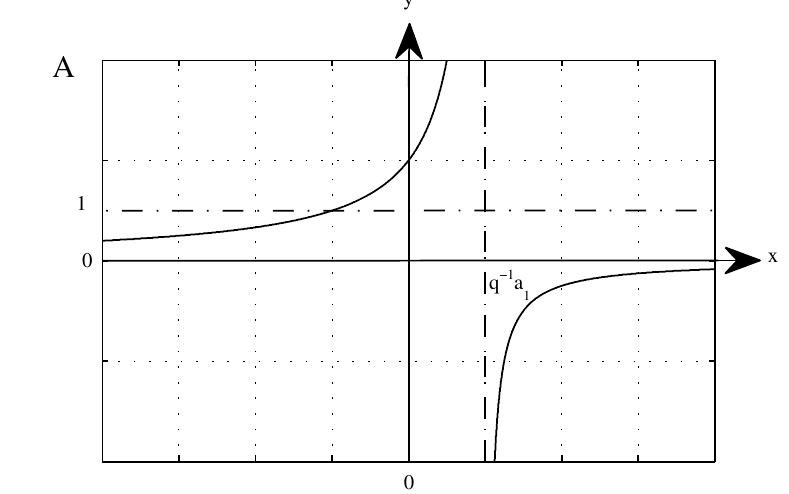}
\hfill
\includegraphics{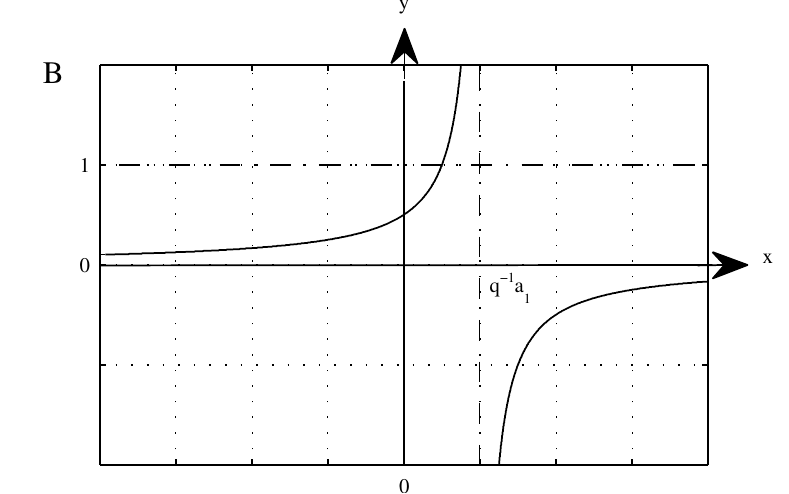}
\caption{The graph of $f(x, q)$ in \textbf{Case 1}. In A, we have 
$y_0>1$ and $a_1>0$ and in B, $0<y_0<1$ and $a_1>0$.}
\label{zeroJ-Lcase1}
\end{figure}
The Case 1 represented in Figure \ref{zeroJ-Lcase1}A as well as the \textbf{Case 2} do not yield any OPS.
From Figure \ref{zeroJ-Lcase1}B, it follows that the only possible interval is $(0, q^{-1}a_1)$
which coincides with the 2\textit{nd} case of Theorem \ref{thm*}.
A completely similar analysis as the one
done in the previous case allows us to conclude that in $(0, a_1]$ an OPS can be defined 
which is orthogonal w.r.t. a suitable $\rho$ supported on the set of points $\{q^ka_1\}_{k\in\mathbb{N}_0}$.
Then, we have the following Theorem.

\begin{theorem}\label{theoremzeroJ-Lcase1iiB}
Let $a_1>0$ and $0<qy_0<1$. 
Let $a=0$ and $b=a_1$ be the zeros of $\sigma_1(x,q)$.
Then, there exists a sequence  of polynomials  
$(P_n)_n$ for $n\in\mathbb{N}_0$ orthogonal 
w.r.t. the weight function (see the $0$-Laguerre/Jacobi case in Table \ref{tab2})
\begin{equation*}
\rho(x, q)=\left|x\right|^{\alpha}(b^{-1}qx; q)_{\infty},
\quad q^{\alpha}=-\frac{q^{-2}\hlf \sigma_2''(0, q)}{\hlf \sigma_1''(0, q)b}
\end{equation*}
supported on the set of points $\{q^ka_1\}_{k\in\mathbb{N}_0}$.
(see \refe{qortho2} of Theorem \ref{thm*}-\textbf{2}).
\end{theorem}
The OPS in Theorem \ref{theoremzeroJ-Lcase1iiB} coincides with
the case IVa4 in Chapter 10 of \cite[pages 278 and 312]{KLS}. 
An example of this family is the little $q$-Laguerre (Wall) 
polynomials $P_n(x; \alpha|q)$
satisfying the $q$-EHT with the coefficients
$$\sigma_1(x, q)=q^{-2}x(a_1-x), \quad \sigma_2(x, q)=ax,$$
$$\tau(x, q)=-\frac{1}{(1-q)q}x+\frac{1-aq}{(1-q)q}\quad {\rm and} \quad \lambda_n(q)=\frac{q^{-n}}{1-q}[n]_q$$
where $a_1=1$. 
The conditions $0<qy_0<1$ and $a_1>0$ give the 
restriction $0<a<q^{-1}$ on the parameter of $P_n(x; \alpha|q)$
with orthogonality on 
$\{...,q^2,q,1\}$ in the sense \refe{qortho2} where
\begin{equation*}
d_n^2=a^nq^{n^2}\frac{(q; q)_{\infty}}{(aq; q)_{\infty}}(q, aq; q)_n.
\end{equation*}
In the literature, this relation can be found as an infinite sum \cite[page 312]{KLS}.


\section{Concluding remarks}
The $q$-polynomials of the Hahn class have been revisited
by use of a direct and very simple geometrical approach based on the 
qualitative analysis of solutions of the $q$-Pearson \refe{openqpearson} 
and the \textit{extended} $q$-Pearson \refe{openexpqpearson1} equations.
By this way, it is shown that it is possible to introduce in a 
unified manner all orthogonal polynomial solutions of the $q$-EHT, 
which are orthogonal w.r.t. a measure supported 
on some set of points in certain intervals.
In this review article we are able to extend the well known orthogonality relations
for the big $q$-Jacobi polynomials (see Theorem \ref{case3aA} and Theorem \ref{theoremcase3cD}), 
$q$-Hahn polynomials (see Theorem \ref{case3aB} and Theorem \ref{theoremcase4C}), and for the 
$q$-Meixner polynomials (see Theorem \ref{linearcase2aA}) to a larger set 
of their parameters.


\subsection*{Acknowledgements}  This work was partially supported by MTM2009-12740-C03-02
(Ministerio de Econom\'\i a y Competitividad), FQM-262, FQM-4643, FQM-7276 (Junta de Andaluc\'\i a),
Feder Funds (European Union), and METU OYP program (RSA). The second author (RSA) thanks the
Departamento de An\'alisis Matem\'atico and IMUS for their kind hospitality during her stay in
Sevilla.


\end{document}